\documentclass[a4paper,10pt]{article}

\usepackage{latexsym}
\usepackage{amsmath}
\usepackage{amsthm}
\usepackage{amssymb}
\usepackage{amscd}
\usepackage{bm}
\usepackage{graphicx}
\usepackage{wrapfig}
\usepackage{fancybox}
\usepackage{enumerate}
\usepackage{url}

\newcommand{\bbc}{{\mathbb C}}

\newcommand{\bbn}{{\mathbb N}}

\newcommand{\bbr}{{\mathbb R}}
\newcommand{\bbs}{{\mathbb S}}
\newcommand{\bbt}{{\mathbb T}}

\newcommand{\bbz}{{\mathbb Z}}

\newcommand{\al}{{\alpha}}

\newcommand{\gam}{{\gamma}}
\newcommand{\Gam}{{\Gamma}}
\newcommand{\del}{{\delta}}
\newcommand{\Del}{{\Delta}}
\newcommand{\ep}{{\epsilon}}
\newcommand{\vep}{{\varepsilon}}
\newcommand{\vph}{{\varphi}}

\newcommand{\lam}{{\lambda}}

\newcommand{\sig}{{\sigma}}
\newcommand{\Sig}{{\Sigma}}

\newcommand{\im}{{\operatorname {Im}}}
\newcommand{\re}{{\operatorname {Re}}}

\newcommand{\R}{\bbr}
\newcommand{\C}{\bbc}

\newcommand{\ovl}[1]{{\overline{#1}}}
\newcommand{\bk}{\backslash}

\newcommand{\str}{{\operatorname{Str}}}

\newcommand{\rd}{{\partial}}
\newcommand{\nab}{{\nabla}}

\newcommand{\cae}{{\mathcal E}}
\newcommand{\caf}{{\mathcal F}}
\newcommand{\cao}{{\mathcal O}}
\newcommand{\cas}{{\mathcal S}}
\newcommand{\mct}{{\mathcal T}}

\newcommand{\loc}{{\operatorname{loc}}}

\newcommand{\nor}[2]{\left\| {#1} \right\|_{#2}}
\newcommand{\inp}[2]{\langle {#1} , {#2} \rangle }

\newcommand{\dist}{\operatorname{dist}}
\newcommand{\rad}{\operatorname{rad}}
\newcommand{\supp}{\operatorname{supp}}

\newcommand{\eh}{\hat{e}}
\newcommand{\fh}{\hat{f}}
\newcommand{\tq}{\tilde{q}}
\newcommand{\tu}{\tilde{u}}
\newcommand{\bj}{\boldsymbol{j}}
\newcommand{\us}{\underline{s}}

\newcommand{\sy}{4\hspace{-0.5pt} /\hspace{-0.5pt} 3}
\newcommand{\sh}{8\hspace{-0.5pt} /\hspace{-0.5pt} 3}
\newcommand{\Na}{3\hspace{-0.5pt} /\hspace{-0.5pt} 2}
\newcommand{\Nb}{1\hspace{-0.5pt} /\hspace{-0.5pt} 2}
\newcommand{\Nc}{3\hspace{-0.5pt} /\hspace{-0.5pt} 4}

\newcommand{\La}{L^{\sy}_t L^{\sy}_x}
\newcommand{\Lb}{L^{\sh}_t L^8_x}
\newcommand{\Lc}{L^{\infty}_t L^2_x}
\newcommand{\Ld}{L^4_t L^4_x}

\newcommand{\scr}{Schr\"{o}dinger }

\newcommand{\dhe}{{\dot{H}_{e}}}

\makeatletter

  \@addtoreset{equation}{section}
\makeatother

\begin{document}

\newtheorem{dfn}{Definition}[section]
\newtheorem{thm}{Theorem}[section]
\newtheorem{lem}{Lemma}[section]
\newtheorem{prop}{Proposition}[section]
\newtheorem{rmk}{Remark}[section]

\title{REMARKS ON LOCAL THEORY FOR SCHR\"ODINGER MAPS NEAR HARMONIC MAPS
\footnote{Running head: LOCAL WELL-POSEDNESS OF SCHR\"ODINGER MAPS}
}
\author{{\sc Ikkei Shimizu}}
\date{\hfill}
\maketitle

\abstract{
We consider the initial-value problem for the equivariant Schr\"odinger maps near a family of harmonic maps. 
We provide some supplemental arguments for the proof of local well-posedness result by Gustafson, Kang and Tsai in [Duke Math. J. 145(3) 537--583, 2008]. %, and prove that their assertions are actually true except uniqueness. 
We also prove that the solution near harmonic maps is unique in $C(I;\dot{H}^1(\R^2)\cap \dot{H}^2(\R^2))$ for time interval $I$. In the proof, we give a justification of the derivation of 
the modified Schro\"odinger map equation in low regularity settings without smallness of energy. 
}

\section{Introduction}

We consider the initial value problem for the Schr\"odinger map equation (or Schr\"odinger flow) 
from $\R^n$ to  a sphere $\bbs^2$, which is given by
\begin{equation}\label{1}
\rd_t u = u\times \Del u,\quad
u(x,0)= u_0(x),
\end{equation}
where $u=u(x,t)$ is unknown function from $\R^n\times \R$ to a sphere
\begin{equation}\label{2}
\bbs^2 = \left\{ y\in \R^3 : |y|=1  \right\} \subset \R^3,
\end{equation}
and $\times$ denotes the vector product of vectors in $\R^3$. 
This equation arises in various ways in physics; we refer, for example, to \cite{D}, \cite{KIK} 
for details. 
%%%%%%%%%%%%%%%%%%%%%%%reference
The equation (\ref{1}) admits the following conserved energy
\begin{equation}\label{3}
\cae (u) = \frac{1}{2} \int_{\R^2} | \nabla u|^2 dx,
\end{equation}
and (\ref{1}) has the scale invariance 
\begin{equation}\label{4}
u(x,t) \mapsto u(\frac{x}{\lam} , \frac{t}{\lam^2}) \text{ for } \lam >0.
\end{equation}
In this work, we restrict ourselves to the case $n=2$. 
%In this paper, we provide some supplemental argument of the paper by Gustafson et al. []. 
Our aim of the present paper is to supplement arguments 
concerning the regularity, which is used without proof in the paper by Gustafson, Kang and Tsai. \cite{GKT2}. 
%To prepare the background of their work, we first review 
\bigskip\par %%%%%%%%%%referance
%We first review the background of their work. 
We first recall the background of the problem. For $m\in\bbn$, a map $u:\R^2 \to \bbs^2$ is 
said to be $m$-\textit{equivariant} if 
$u$ has the form 
\begin{equation}\label{5}
u(x) = e^{m\theta R}v(r),
\end{equation}
where $(r,\theta)$ is the polar coordinates of $x$, $v={}^t (v_1,v_2,v_3)$ is a function from $(0,\infty)$ to $\R^3$, 
and $R$ is the matrix 
$
R=
\begin{pmatrix}
0 & -1 & 0\\
1 & 0 & 0\\
0 & 0 & 0
\end{pmatrix}
$. 
Note that $e^{\al R}$ represents a rotation in angle $\alpha$ around the $u_3$-axis for $\al\in\R$. 
We can observe that the Schr\"odinger map equation (\ref{1}) formally preserves the $m$-equivariance; i.e., 
if $u_0$ is $m$-equivariant, then the solution $u(t)$ to (\ref{1}) is $m$-equivariant for all $t$. 
Thus, it makes sense to restrict our function space to the $m$-equivariant class. \bigskip\par
If an $m$-equivariant map $u=e^{m\theta R}v(r)$ satisfies $\cae (u)<\infty$, the followings hold: 
$v:(0,\infty) \to \R^3$ is continuous, and both limits $v(0) :=\lim_{r\to 0} v(r)$ and $v(\infty):= \lim_{r\to \infty} v(r)$
 exist and are equal to either of $\pm \vec{k}$, where $\vec k = {}^t(0,0,1)$. 
%For the proof of this fact, see [BT] for example. 
(For the proof, see Section 4 below.) 
By the rotational symmetry of (\ref{1}), we may assume $v(0)=-\vec{k}$ without loss of generality. 
Here, we have two choices of $v(\infty)$:
$$
v(\infty)= -\vec{k} \quad \text{ or }\quad v(\infty)= \vec{k}.
$$
Each choice corresponds to different homotopy type of maps. 
When $v(0)=v(\infty)= -\vec{k}$ and $\cae (u) <8\pi m$, 
the image of $u$ never covers the whole sphere $\bbs^2$, which implies that $u$ is homotopic to 
a constant map $Q$ for some $Q\in\bbs^2$ (see \cite{BIKT2}). 
%We can say that this case is less geometric. 
On the other hand, when $v(0)=-\vec k$ and $v(\infty) =\vec k$, 
the image of $u$ must cover $\bbs^2$, thus $u$ is \textit{not} homotopic to a constant map. 
In this paper, we focus on the latter case.\bigskip\par
%%%%%%%%%%%%%%%%%%%%%%%%reference
%In what follows, we assume $m\ge 1$. 
Our function space is 
\begin{equation}\label{6}
\Sigma_m := \left\{ u=e^{m\theta R}v(r)\ |\ u\in\dot H^1 , v(0)=-\vec{k}, v(\infty)=\vec{k}  \right\} .
\end{equation}
Then $\Sig_m$ is complete metric space with metric 
$d(u,\tilde{u}) = \nor{u-\tilde{u}}{\dot H^1}$. 
%(It may seem to be strange since $\nor{\cdot}{\dot{H}^1}$ is not a norm. 
%However, $\nor{\cdot}{\dot{H}^1}$ turns to be a complete metric under restriction to $m$-equivariant class.) 
For $u\in \Sig_m$, we have 
\begin{equation}\label{7}
\cae(u) = \pi \int_0^\infty \left( |\rd_r u |^2 + m^2 \frac{{u_1}^2 +{u_2}^2}{r^2}  \right) r\, dr.
\end{equation}
Then we can write (\ref{7}) as 
\begin{equation}
\cae(u) = \pi \int_0^\infty \left| v_r - \frac{m}{r} J^v Rv \right|^2 rdr + 4\pi m
\end{equation}
for $u\in \Sig_m$, where we define $J^v:= v\times \cdot$. (See \cite{GKT} for details.) 
Thus, we have $\cae (u) \ge 4\pi m$ for all $u\in\Sig_m$, 
and $u$ minimizes the energy if and only if $v_r - \frac{m}{r}J^vRv =0$ for almost all $r\in (0,\infty)$. 
By solving this ODE, it turns out that the minimizing set can explicitly be written as 
\begin{equation}\label{8}
\cao_m= \left\{ e^{\al R}Q(\frac{\cdot}{s}) \ |\ s>0, \al\in\R \right\},
\end{equation}
where $Q := e^{m\theta R}h(r)$, 
$h(r)={}^t ( h_1(r), 0, h_3(r) )$, 
$h_1(r)= \frac{2r^m}{r^{2m} +1}$, 
$h_3(r)= \frac{r^{2m}-1}{r^{2m} +1}$.
We call $\cao_m$ the family of \textit{harmonic maps}. 
Note that any element in $\cao_m$ is a stationary solution to (\ref{1}): $u\times \Del u =0$.\bigskip\par
%In \cite{GKT}, Gustafson et al. gives a geometric description of $\cao_m$. 
We recall a geometric description of $\cao_m$ obtained in \cite{GKT}.

\begin{prop}\label{P1} (\cite{GKT})
There exist $\del_0 >0$ and $C_0, C_1>0$ such that 
for $u\in \Sig_m$ with $\cae (u) <4\pi m + {\del_0}^2$, the following hold:
\begin{enumerate}[(i)]
\item There exist unique $s_*=s_*(u)\in (0,\infty)$ and $\al_*=\al_*(u) \in \bbt^1$ such that 
\begin{equation}\label{9}
\dist_{\dot{H}^1} (u ,\cao_m ) = \nor{u- e^{\al_* R} Q(\frac{\cdot}{s_*})}{\dot{H}^1}. 
\end{equation}
\item The map $u\mapsto (s_*(u), \al_*(u))$ is continuous.
\item $C_0 \dist_{\dot{H}^1} (u ,\cao_m ) \le \sqrt{\cae (u) - 4\pi m} \le C_1 \dist_{\dot{H}^1} (u ,\cao_m ) $.
\end{enumerate}
\end{prop}

The above proposition ensures the unique existence of $\dot{H}^1$-closest harmonic map for each $u\in \Sig_m$ with $\cae (u) - 4\pi m \ll 1$. 
By the scaling pair $(s_*(u), \al_*(u))$, we can get precise information on the position of $u$ relative to $Q$ along the harmonic map family.

%%%%%%%%%%%%%%%%%%
%%%%%%%%%%%%%%%%%%
%%%%%%%%%%%%%%%%%%

\section{The Paper by Gustafson et al. and Our Main Result}

With the aim of studying the stability of $\cao_m$, 
Gustafson, Kang, and Tsai \cite{GKT2} considers local problems for (\ref{1}) 
near the family $\cao_m$ in the class $\Sig_m$. 
To present the statement of their results, we introduce the notion of weak solution of (\ref{1}). 
We first note that the equation (\ref{1}) can be written by divergence form as follows:
\begin{equation}\label{10}
\rd_t u= %u\times \Del u = 
\sum_{j=1}^2 \rd_{x_j} (u \times (\rd_{x_j} u) )
\end{equation}
where $x_j$ is the $j$-th spatial coordinate. %This form enables us a weak formulation of (\ref{1}) as follows. 
Considering (\ref{10}), we define the weak solution to (\ref{1}) in the following way.

\begin{dfn}
For interval $I$, $u(t) \in L^\infty_{\loc} (I;\Sig_m)$ is said to be a weak solution 
if $u(t)$ satisfies 
\begin{equation*}
\int_{I\times \R^2} u\rd_t \phi\  dxdt = \int_{I\times \R^2} \sum_{j=1}^2(u\times \rd_{x_j} u) 
\rd_{x_j} \phi\ dxdt 
\end{equation*}
for all $\phi\in C^\infty_0 (I\times \R^2)$.
\end{dfn}

The pioneering paper \cite{GKT2} by Gustafson et al. contains the following results (Theorem 1.4, page 543): 

\begin{enumerate}
\item[(LP1)] (Existence) There exist $\del_0>0$, $\sig >0$, and $C>0$ such that the following holds: 
If $u_0\in\Sig_m$ satisfies $\del :=\sqrt{\cae(u_0) -4\pi m} <\del_0$, 
then (\ref{1}) has a weak solution $u(t)\in C(I; \Sig_m)$, $I=[0,\sig {s_0}^2]$, where $s_0:= s_*(u_0)$. 
\item[(LP2)] (Uniqueness) The above solution is unique in $C(I;\Sig_m)$; i.e., 
if $\tilde{u}(t) \in C(I';\Sig_m)$ satisfies (\ref{1}) for $I'=[0,T]$ with $T>0$, 
then $\tilde{u}(t)=u(t)$ for all $t\in I\cap I'$.
\item[(LP3)] (Energy conservation) The above solution conserves the energy, that is, $\cae (u(t)) =\cae(u_0)$ for all $t\in I$.
\item[(LP4)] (Regularity) If we further assume that $u_0\in \dot{H}^2$, then the above solution $u(t)$ 
is in $C(I; \Sig_m\cap \dot{H}^2)$. 
\item[(LP5)] (Continuous dependence) The map $\{u\in \Sig_m : \cae (u) < 4\pi m +{\del_0}^2\} \ni u_0 \mapsto u(t) \in C(I;\Sig_m)$ is continuous. 
\end{enumerate}

These assertions play an important role in the investigation of  
global behavior of the solution to (\ref{1}) near the harmonic map family $\cao_m$. 
Indeed, the ensured existence time $\sig s_0^2$ in (LP1) implies that 
the possible finite time blow-up senario for (\ref{1}) is $s_*(u(t)) \to 0$. 
%namely, the \textit{energy concentration of closest harmonic map}. 
See \cite{GKT2} for more details. (See also \cite{GNT}, \cite{MRR} and \cite{P}.) \bigskip\par
%Among their great works, there are three points of lack of argument in their proof. 
%Although there works are great, there are three main points which need to be fixed or 
%supplemented. 
In the present paper, we mainly focus on the following three points which are not explicitly mentioned in their paper. 
The first one is concerned with the limiting argument in their proof. 
Their way to show (LP1) is to reduce the problem to a PDE-ODE system (\ref{13}) and (\ref{22}) defined below. 
For the construction of weak solution, they first 
construct a solution $(q(t),s(t),\al(t))$ of (\ref{13}) and (\ref{22}), then reconstruct the original map $u(t)$ from it. %$(q(t),s(t),\al(t))$. 
Then, %it is necessary to prove that $u(t)$ is actually a weak solution to (). 
they claim that this $u(t)$ is actually a weak solution. 
To show that, they approximate $u(t)$ by smooth solutions. 
%In order to work this argument well,
In the argument, they implicitly use the fact that the maximal existence time of each element of approximating sequence of solutions $\{ u^k(t)\}_{k=1}^\infty$ is bounded from below uniformly in $k$. 
%Therefore, we provide a proof of this fact in the present paper. \par
Our first aim is to provide a proof of this fact.\par 
The second point is related to the regularity persistence stated in (LP4). 
% one is related to regularity propagation (LP4). 
%In their paper, they provides a priori bound for $\nor{\Del (e^{i(m+1)\theta}q)}{L^2}$, which leads to 
%a priori bound for $\nor{u}{\dot{H}^2}$. 
%and this implies that the solution $u(t)$ is in $L^\infty (I;\Sig_m\cap\dot{H}^2)$. %-bounded with respect to time $t$. 
%but this does not mean that (LP4) immediately. 
%However, it is not trivial that $u(t)$ is continuous as an $\dot{H}^2$-valued function, namely $u(t)\in C(I;\dot{H}^2)$. 
%In this paper, we give a proof of the fact that the a priori bound implies $u(t)\in C(I;\dot{H}^2)$. 
%To this end, we need to check the continuity of the map $H^1\times \R^+\times\bbt^1\ni (q,s,\al)\mapsto u \in \Sig_m\cap \dot{H}^2$.
In their argument, there is no explicit mention of how we ensure the continuity of the map $u:I\to \dot{H}^2$. Hence, we give a proof of this fact in the present paper.
\par
The third one is concerned with the uniqueness of solutions. 
They implicitly reduce the problem to the modified system (\ref{13}) and (\ref{22}), and then 
show the uniqueness for these equations. 
%this does not mean that of original equation (\ref{1}). 
%In order to deduce (LP2) to their result, 
%For the reduction, 
For sufficiently smooth solutions (more precisely, when $u(t)\in C(I;\dot{H}^1\cap \dot{H}^3)$), the reduction is justified since it is known that 
the corresponding $(q(t),s(t),\al(t))$ actually satisfies (\ref{13}) and (\ref{22}). 
However, 
%the derivation has not been justified %only in the case that the solution is sufficiently smooth; 
%in the case that $u(t)$ is a general weak solution, 
%this justification has not been done. 
such kind of justification is not given for weak solutions. 
%more precisely, 
%for the class $C(I; \dot{H}^1\cap \dot{H}^3)$, which is given in \cite{BT}. 
Hence, we attempt to give a new justification of the derivation of (\ref{13}) and 
(\ref{22}) in a larger class of solutions. \bigskip\par
In this paper, we reproduce the proof of (LP1), (LP3), and (LP4). Moreover, we show 
the restated propositions (LP2)' and (LP5)' as follows. 
%For uniqueness, we assert that the solutions are unique in the class $C(I;\Sig_m\cap \dot{H}^2)$. 
%The main result is the following: 

\begin{thm}\label{T1}
The propositions (LP1), (LP3), and (LP4) hold. Moreover, 
the restated propositions (LP2)' and (LP5)'  hold.
\begin{itemize}
\item[(LP2)'] %In the class $C(I; \Sig_m\cap \dot{H}^2)$, the solutions to (\ref{1}) are unique. 
For $u_0\in \{u\in \Sig_m : \cae (u) < 4\pi m +{\del_0}^2\}\cap \dot{H}^2$, the solutions to (\ref{1}) 
are unique in $C(I;\Sig_m\cap \dot{H}^2)$. 
\item[(LP5)'] The map $\{u\in \Sig_m : \cae (u) < 4\pi m +{\del_0}^2\}\cap \dot{H}^2 \ni u_0 \mapsto u(t) \in C(I;\Sig_m)$ is continuous with $\dot{H}^1$-topology. 
Moreover, the above map can be uniquely extended to $\{u\in \Sig_m : \cae (u) < 4\pi m +{\del_0}^2\}$ 
as a limit of $C(I;\Sig_m\cap\dot{H}^2)$-solutions, and coincide with 
the solutions constructed in (LP1).
\end{itemize}
\end{thm}

\begin{rmk}
The uniqueness stated in (LP2)' is more restricted than that in (LP2), 
but stronger than that ensured in \cite{M}. 
The statement in (LP5)' is essentially unchanged from that in (LP5), 
since %the reason why
we just make the definition of solution map clearer according to (LP2)'.
\end{rmk}

We briefly explain how we supplement the points mentioned above. 
For the first one, we establish a priori estimate for second derivative of $q(t)$ (see (\ref{c31})), 
which leads to a priori bound for third order derivative of $u(t)$. 
Hence, by McGahagan \cite{M}, we can ensure that 
$u^k(t)$ continues to exist as long as the corresponding solution to (\ref{13}) and (\ref{22}) exists. \par
%the lower bound of maximal existence time of approximating sequence $\{ u^k(t) \}$. \par
The essential step for the second point is the continuity of reconstruction stated in Lemma \ref{L10}, 
in which we claim that the map 
$H^1\times \R^+\times\bbt^1\ni (q,s,\al)\mapsto u \in \Sig_m\cap \dot{H}^2$ is continuous. \par
%In terms of the third problem, %weakening their claim, we assert tha
To prove (LP2)', which is concerned with the third point, 
we make a new justification of the derivation of (\ref{13}) and (\ref{22}) for the solutions $u(t)\in 
C(I;\dot{H}^1)\cap L^\infty (I;\dot{H}^2)$. 
This class is the lowest regularity ever. 
%This enable us to reduce the problem to the uniqueness of PDE-ODE system (\ref{13}) and (\ref{22}), which has been already established by \cite{GKT2}. 
This immediately leads to (LP2)' by the uniqueness of the system (\ref{13}) and (\ref{22}) 
established in \cite{GKT2}. 
The main difficulty is that the calculation needs to be performed in the distributional class, 
while we have to use the polar coordinates essentially. 
This is why we introduce a new function class $H^{-1}_e$ defined below. 
In this space, several kinds of calculation related to polar coordinates are justified 
in a larger class than $L^2$. \par
The rest of the proof is essentially a reproduction of the argument of \cite{GKT2}, while we make small modifications.\bigskip\par
Here, we make a few remarks on the preceding results concerning the well-posedness for Schr\"odinger maps. 
%The local and global well-posedness 
%of Schr\"odinger maps (\ref{1}) have been extensively studied. 
%The oldest result is established 
In \cite{SSB}, the existence of global weak solution 
$u(t)\in L^\infty (\R ; \dot{H}^1)$ is established, which, however, says nothing about the singularities or uniqueness. 
The local well-posedness for large data has been studied, for example, by \cite{SSB}, \cite{DW}, \cite{M}, and \cite{KLPST}. The lowest regularity is the work by McGahagan \cite{M}, in which the 
existence and uniqueness of solutions is established in $L^\infty (I;\dot{H}^1\cap \dot{H}^3)$. 
(We essentially use this  result in the present paper.) 
%It is worth noting that \cite{DW}, \cite{M}, and \cite{KLPST} consider some kinds of \textit{perturbed equation}, and construct solutions by taking limit so that 
%the perturbing term tends to $0$. 
%The common weak point to these method is that we need high regularity for function spaces. 
For small data, the local and global well-posedness have been extensively studied, 
and the cutting-edge result is \cite{BIKT}. 
%(see also \cite{IK}, \cite{IK2}, \cite{BIK}, \cite{B}, \cite{B2}, \cite{BIKT2}, and \cite{Dod}, for example).
The propositions (LP1)--(LP5) by Gustafson, Kang and Tsai \cite{GKT2}, which is the main subject of the present work, is 
the first local-in-time result for rough data near harmonic maps under the restriction to equivariance. 
%and includes geometrical implication. 
\bigskip\par
%Based on the characterization of blow-up dynamics of the solution near harmonic maps in \cite{GKT2}, 
%The global behavior of the solutions near harmonic maps has been explored in the last decade. 
%In \cite{GKT}, Gustafson, Kand, and Tsai shows that the harmonic map family $\cao_m$ is orbitally stable for $m\ge 1$. 
%In the subsequent paper \cite{GKT2} and the paper by Gustafson, Nakanishi, and Tsai \cite{GNT}, it is shown that $\cao_m$ is asymptotically 
%stable in energy space when $m\ge 3$. 
%On the other hand, for $1$-equivariant case, Bejenaru and Tataru \cite{BT} shows that $Q(x)= e^{\theta R} h(r)$ is not stable. After that, 
%Merle, Rapha\"el, and Rodnianski \cite{MRR} and Perelman \cite{P} construct a family of initial data in $\Sig_1$ which intersects 
%arbitrary small neighborhood of $\cao_1$, and the corresponding solutions blows up in finite time. 
%(They also obtain the blow-up profile of these solutions.) 
%This indicates that $\cao_1$ is not %orbitally stable, but not
%asymptotically stable. 
%(Energy concentration can occur.)\bigskip\par

The organization of this paper is as follows. 
In Section 3, we provide a proof of Theorem \ref{T1}, reproducing 
the argument in \cite{GKT2}. The subsequent sections is devoted to 
the proof of technical lemmas. 
In detail, we provide a justification of derivation of the modified system (\ref{13}) and (\ref{22}) 
in Section 4. 
In Section 5, we derive the a priori estimates (\ref{c30}) and (\ref{c31}). 
In Section 6, a detailed proof of the properties of scaling pair $(s,\al)$ is given. 
%Especially, we take care of its regularity properties. 
In Section 7, the continuity of the map $H^1\times \R^+\times\bbt^1\ni (q,s,\al)\mapsto u \in \Sig_m\cap \dot{H}^2$ is shown. In the same section, we provide a proof of a lemma 
concerning approximation by smooth maps.\par
We close this section with introducing notations used in the present paper. 
We set 
$\bbn:= \{n\in \bbz : n\ge 1 \}$ and % is the set of all positive integers. 
$\R^+ := \{ s\in\R : s>0 \}$.% is the set of all positive numbers. 
We use the letter $C$ in many times to indicate a constant, and the representing quantity 
varies from each situation, if there is no risk of mathematical validity. 
For $p,q\in [1,\infty]$ and for interval $I\subset\R$, we sometimes abbreviate $L^p(I;L^q(\R^2))$ as $L^p_tL^q_x (I)$, 
or $L^p_t L^q_x$. 
We define $L^2_{\rad} \equiv L^2_{\rad} (\R^2) := \{ f\in L^2(\R^2) : f \text{ is radially symmetric.} \}$. %is the set of all radially symmetric functions in $L^2(\R^2)$. 
For a Banach space $X$, $\inp{\cdot}{\cdot}_{X^*,X}$ denotes the 
coupling of the elements in $X^*$ and in $X$. And for a Hilbert space $H$, 
$\inp{\cdot}{\cdot}_{H}$ denotes the inner product of $H$.
%\begin{equation*}
%\nor{f}{L^p_t L^q_x}
%\end{equation*}

%%%%%%%%%%%%%%%%%%%%
%%%% Section 3 %%%%%
%%%%%%%%%%%%%%%%%%%%

\section{Proof of Theorem \ref{T1}}

\subsection{Coulomb Gauge and Modified Schr\"odinger Map}
%In what follows, we assume that $m>0$. The case $m<0$ is led from $-m$ 
%by change of variables. \bigskip\par
We begin with a proposition concerning %the \textit{Coulomb gauge}
the choice of frame, which we show 
in more general setting in Section 4 (see Lemma \ref{L1}). 

\begin{prop} (\cite{BT})
Let $u=e^{m\theta R}v(r)\in \Sig_m$. Then, there exists $\eh (r):(0,\infty)\to \R^3$ such that 
\begin{enumerate}[(i)]
\item $\eh$ is absolutely continuous on any closed subinterval of $(0,\infty)$.
\item $\lim_{r\to \infty} v(r)= {}^t (1,0,0)$.
\item $\rd_r\eh (r) = - (\eh(r)\cdot v_r(r)) v(r)$ for $r\in (0,\infty)$.
\end{enumerate}
\end{prop}
From the properties of $\eh$, it follows that 
$$
|\eh(r)|\equiv 1,\quad \eh (r) \cdot v(r) =0 \text{ for all } r\in (0,\infty).
$$
Therefore, %$\{\eh (r), J^v(r) \eh(r) \}$ forms an orthonormal frame of $T_{v(r)}\bbs^2$, 
%where $T_p \bbs^2$ denotes the tangent space of $\bbs^2$ at $p\in\bbs^2$. 
$\{e^{m\theta R}\eh, J^{u} e^{m\theta R}\eh \}$ forms an orthonormal frame of 
$T_{u}\bbs^2$, where $T_p \bbs^2$ denotes the tangent space of $\bbs^2$ at $p\in\bbs^2$. 
(In other words, this is an orthonormal frame of the tangent bundle $u^{-1}T\bbs^2$) 
%when we consider its component as an element of its section $\Gam (u^{-1}T\bbs^2)$.) 
This choice of frame is called the \textit{Coulomb gauge} (or \textit{Coulomb frame}). \bigskip\par
For $u\in\Sig_m$, 
define $q=q(u)\in L^2_{\rad}$ and $\nu = \nu(u)\in L^\infty (\R^2)$ by 
\begin{equation}\label{11}
q := (v_r - \frac{m}{r}J^v Rv)\cdot(\eh + iJ^v\eh),\quad \nu := J^v Rv\cdot(\eh + iJ^v\eh).
\end{equation}
%This quantity corresponds to the representation of $v_r - \frac{m}{r}J^v Rv$ 
%in complex coordinates associated with the frame ${\eh ,J^v\eh}$. 
%We also define 
%\begin{equation}\label{12}
%a
%\end{equation}
%Let $u(t)\in C(I;\Sig_m)\cap C(I;\dot H^2)$ be a solution to (\ref{1}) for some interval $I\subset \R$. 
Let $I\subset \R$ be an interval. 
For $u(t)\in C(I;\Sig_m)$, we write $q(t):=q(u(t))$, and set $\tq (t)=e^{i(m+1)\theta}q(t)$. 
Then, the following holds.

\begin{prop}\label{P3}
Let $u(t)\in C(I;\Sig_m)\cap C(I;\dot H^2)$ be a solution to (\ref{1}). 
Then, $\tq (t)\in C(I; H^1)\cap C^1(I; H^{-1})$, and $\tq (t)$ satisfies %the following equation:
\begin{equation}\label{13}
i \tilde{q}_t +\Del \tilde{q} = \frac{m(1+v_3)(mv_3-m-2)}{r^2}\tq +\frac{mv_{3r}}{r}\tq 
+\tq N(q)
\end{equation}
where
\begin{equation}\label{14}
N(q) := \re \int_{r}^{\infty} \left( \bar{q} +\frac{m}{r}\bar{\nu} \right)
 \left( q_r +\frac{1-mv_3}{r} q \right) dr.
\end{equation}
Moreover, if a solution $u(t)$ is in $C(I;\Sig_m)\cap L^{\infty} (I;\dot{H}^2)$, then 
$\tq (t)\in C(I;L^2)\cap L^\infty(I;H^1)\cap W^{1,\infty} (I;H^{-1})$, and (\ref{13}) holds.
\end{prop}

The equation (\ref{13}) is called the \textit{modified Schr\"odinger map} equation.
The equation (\ref{13}) is first obtained by Chang et al. \cite{CSU} by formal calculations. 
Recently, Bejenaru and Tataru \cite{BT} shows that the calculation is rigorous when $u(t)\in C(I;\Sig_m\cap \dot{H}^3)$. 
Our claim is that the calculation can also be justified when $u(t)\in C(I;\Sig_m)\cap L^\infty(I;\dot{H}^2)$. 
Moreover, our proof of the proposition does not need any condition of smallness of energy, 
although this improvement brings no benefit on our main argument. 
We prove this proposition in Section 4. 
%Note that under the regularity condition $u(t)\in C(I;\Sig_m)\cap C(I;\dot H^2)$,
% the conservation of energy can be obtained by direct differentiation with respect to $t$. 
%Hence, the condition for energy is 

\subsection{Scaling}
%Since we have obtained a relatively accessible equation (\ref{13}), %(\ref{13}) one kind of nonlinear Schr\"odinger equation, 
%it is natural to tempt to reduce the original problem (\ref{1}) to (\ref{13}). 
%However, 
The quantity $q$ does not have sufficient information about the original map $u$. 
Indeed, we have 
$$
q=0 \Longleftrightarrow u\in\cao_m,
$$
which %tells us that for each $q\in L^2_{\rad}$, it remains 
indicates the scale indefiniteness 
of original map. 
(Note that the situation is different in the case when $v(0)=v(\infty)=-\vec{k}$, 
where $u$ can be completely reconstructed from $q$. See \cite{BIKT2} for more details.) 
From this observation, we need to consider the information about the position of $u$ along $\cao_m$. 
To this end, it seems to be natural to see $(s_*(u), \al_*(u))$ defined in Proposition \ref{P1}. 
%Nevertheless, Gustafson et al. \cite{GKT2} points out that 
%there is a better choice of scaling for their analysis, as follows.
However, we make a different choice of scaling by following \cite{GKT2}, instead of $(s_*(u),\al_*(u))$. \bigskip\par
\vspace{-5pt}
We make here two preparations. 
First, we introduce Hilbert space $\dhe^1$ as follows:
\begin{equation}\label{15}
\dhe^1 := \{ f:(0,\infty) \to \C \ |\ \nor{f}{\dhe^1}<\infty  \},
\end{equation}
\begin{equation}\label{16}
\inp{f}{g}_{\dot{H}^1_e} := \int_{0}^{\infty} \left( f_r \ovl{g_r} +\frac{m^2}{r^2} f\ovl{g} \right) r\,dr .
\end{equation}
%Note that $v_1,v_2\in \dhe^1$ for $u\in\Sig_m$. 
%Thus $\dhe^1$ roughly corresponds to the $\dot{H}^1$-norm for $m$-equivariant maps. 
%Moreover, we have 
%$\nor{f}{L^\infty} \le C \nor{f}{\dhe^1}$ for some constant $C$ depending only on $m$. 
%(This is one of the great advantages of equivariance.) 
Some properties of this space are observed in Section 4. 
The notation $\dot{H}^1_e$ is adopted from \cite{BIKT2}, 
while \cite{GKT} and \cite{GKT2} use $X$ instead.\par
%the notation in \cite{BT}.\par
Next, we set 
\begin{equation}\label{17}
\boldsymbol{j}:= {}^t (0,1,0)
%\begin{pmatrix}
%0\\
%1\\
%0
%\end{pmatrix}
,\qquad 
J^{h}\boldsymbol{j} := {}^t (-h_3,0,h_1)
%\begin{pmatrix}
%-h_3\\
%0\\
%h_1
%\end{pmatrix}
.
\end{equation}
Then, $\{ \bj, J^h \bj, h \}$ forms an orthonormal basis of $\R^3$. 
Using (\ref{17}), for given $s>0$, $\al\in\bbt^1$, and for a map $u\in\Sig_m$, 
we can decompose  
\begin{equation}\label{18}
e^{-\al R}v(s\cdot ) = z_1 \bj + z_2 J^h\bj + (1+\gam) h, 
\end{equation}
and we set $z= z_1 +iz_2$. %Then, if $\cae(u) -4\pi m\ll 1$, we have
%
%aaaaaaaaaaaaa

\begin{prop}\label{P4}
(\cite{GKT2})
There exist $\del_0>0$, $C>0$ such that the followings hold:
\begin{enumerate}[(i)]
\item For $u\in\Sig_m$ satisfying $\del:= \sqrt{\cae (u) - 4\pi m} <\del_0$, 
there exist $(s,\al)=(s(u),\al(u))\in \R_{>0}\times\bbt^1$ such that 
\begin{equation}\label{19}
\inp{z}{h_1}_{\dhe^1} =0
\end{equation}
\begin{equation}\label{20}
\left| \frac{s}{s_*(u)} -1 \right| + \left|\al - \al_* (u)\right| \le C\del
\end{equation}
\item For the $u$ above, if $(\tilde{s},\tilde{\al})$ satisfies (\ref{19}) and 
%\begin{equation}\label{21}
$
| {s_*(u)}^{-1} s - 1| +| \al - \al_* (u)| \le C\del_0,
$
%\end{equation} 
then $(\tilde{s},\tilde{\al})=(s,\al)$. 
\item If $u(t)\in C( I ;\Sig_m)\cap C^1(I;L^2(\R^2))$ for some open interval $I\subset\R$, then $s(u(t)),\al(u(t))$ are $C^1$. 
If $u(t)\in C(I;\Sig_m)\cap W^{1,\infty}(I; L^2)$ for some open interval $I\subset \R$, then 
$s(u(t)),\al(u(t))\in W^{1,\infty}(I;\R)$. 
\end{enumerate}
\end{prop}

%We make a few comments on this proposition. 
%Compared with the choice $(s_*(u) ,\al_*(u))$, the scaling in Proposition \ref{P4} is more compatible to 
%the inner product $\inp{\cdot}{\cdot}_{\dhe^1}$. Hence, we can obtain relatively simpler equations of  $(s,\al)$. 
%In \cite{GKT2}, 
%it is asserted with no proof that the choice of $(s,\al)$ satisfying (\ref{19}) is unique without 
%the restriction (\ref{21}). However, this might not be always the case. 
%Nevertheless, it is sufficient for our purpose that the uniqueness holds among such $(s,\al)$ as (\ref{21}). 
%We also note that we make a close observation of regularity, which is not shown  explicitly in \cite{GKT2}. 
%Concerning the above, 
A proof of Proposition \ref{P4} is given in Section 6.\bigskip\par
\vspace{-4pt}
For $u\in\Sig_m$, 
we have extracted three quantities $q(u)\in L^2_{\rad}$, $s(u)>0$, and $\al \in \bbt^1$. 
Conversely, 
we can show that $(q,s,\al)\in L^2_{\rad}\times \R^+\times \bbt^1$ possesses enough information to reconstruct the original map $u\in\Sig_m$ completely, which is stated in the following proposition.\\[-5pt] 

\begin{prop}{(\cite{GKT2})}\label{P5} % (\cite{GKT2})
There exists $\del_0>0$ such that for $(q,s,\al)\in L^2_{\rad}\times \R^+\times \bbt^1$ 
with $\del := \nor{q}{L^2} <\del_0$, there is a unique $u\in\Sig_m$ which satisfies 
$(q,s,\al)=(q(u),s(u),\al(u))$. 
Moreover, the map 
$
L^2_{\rm{rad}}\times\R^+ \times\bbt^1 \ni (q,s,\al) \mapsto u \in \Sig_m
$
 is continuous.
\end{prop}
The proof of Proposition \ref{P5} can be found in \cite{GKT2}, Lemma A.2.%, but 
%for reader's convenience, we provide a proof in Appendix. 
\bigskip\par
Now, 
let $u(t)\in C(I;\Sig_m)\cap L^\infty(I;\dot{H}^2)$ be a solution to (\ref{1}), 
and set $s(t):=s(u(t))$, $\al(t) :=\al(u(t))$. 
Direct calculations yield the equation which $s(t)$ and $\al(t)$ satisfy. 

\begin{prop}{(\cite{GKT2})}\label{P6} %(\cite{GKT2})
There exists $\del_0>0$ such that the following holds: 
If $u(t)\in C(I;\Sig_m)\cap L^\infty(I;\dot{H}^2)$ is a solution to (\ref{1}) which satisfies 
$\del:=\sqrt{\cae (u) - 4\pi m} <\del_0$, then 
\begin{equation}\label{22}
\begin{pmatrix}
s(t)\\
\al (t)
\end{pmatrix}
=
\begin{pmatrix}
s(0)\\
\al(0)
\end{pmatrix}
+
\int_{0}^{t}
\left\{
\begin{pmatrix}
0 & -(ms)^{-1}\\
s^{-2} & 0
\end{pmatrix}
(\nor{h_1}{\dot{H}^1_e}^2 I + A)^{-1} \vec{G_2}
\right\}
(\tau) d\tau,
\end{equation}
where 
\begin{equation}
A
=
\begin{pmatrix}
\inp{h_1}{\gamma h_1 -z_2 h_3}_{\dhe^1},
 & 
\frac{1}{m} \inp{rN_0h_1}{z_{1r}}_{L^2} \\
\inp{h_1}{z_1h_3}_{\dhe^1}
 & 
\inp{h_1}{\gamma h_1}_{\dhe^1} + \frac{1}{m}\inp{rN_0h_1}{z_{2r}}_{L^2}
\end{pmatrix}
,
\end{equation}
\begin{equation}
N_0 = -\rd_r^2 -\frac{1}{r}\rd_r +\frac{m^2}{r^2},\quad
L_0 = \rd_r + \frac{m}{r}h_3 ,
\end{equation}
\begin{equation}\label{25}
I=
\begin{pmatrix}
1 & 0\\
0 & 1
\end{pmatrix}
,\quad
\vec{G}_2 = 
\begin{pmatrix}
\inp{L_0N_0h_1}{L_0 z_2}_{L^2} \\
-\inp{L_0N_0h_1}{L_0 z_1}_{L^2}
\end{pmatrix}
+
\begin{pmatrix}
\re G_1 \\
\im G_1
\end{pmatrix}
,
\end{equation}
\begin{equation}\label{26}
\begin{aligned}
G_1 = 
\int_{0}^{\infty} 
(
ig_r ( &-\gam z_r +z\gam_r ) 
+ \frac{m}{r} h_1g (-2\gam_r - iz_2z_{1r} + iz_1z_{2r}) \\
& + \frac{m}{r} (h_1g)_r (\gam^2 -iz_2z) 
+ i\frac{m^2}{r^2}(2h_1^2-1)g\gam z\\
& - i \frac{2m^2}{r^2}h_1h_3 g z_2z
)
r\,dr.
\end{aligned}
\end{equation}
\end{prop}

For the proof of this proposition, see \cite{GKT2}, Section A.2. 
%The differentiability of $s(t)$ and $\al(t)$ is ensured by Proposition \ref{P4} (iii). 
%which leads the differentiability of $z(t)$. 

\subsection{Local well-posedness of the PDE-ODE System}
We have seen that 
if $u(t)\in C(I;\Sig_m)\cap L^\infty(I;\dot{H}^2)$ is a solution to (\ref{1}), 
then $(\tq (t),s(t),\al(t))$ must satisfy the system of equations (\ref{13}) and (\ref{22}). 
Note that this is a closed system. 
Indeed, the quantities such as $v$ and $z$ in (\ref{13}) and (\ref{22}) can be reconstructed from $(q,s,\al)$ 
by Proposition \ref{P5}. 
%, thanks to the reconstruction in Proposition \ref{P5}. 
In the converse direction, %we consider the local-wellposedness of the PDE-ODE system (\ref{13}), (\ref{22}). 
it is reasonable to expect that if $(\tq (t),s(t),\al(t))$ is a solution to the system (\ref{13}) and (\ref{22}), 
then the reconstructed map $u(t)$ is a weak solution to (\ref{1}). 
Hence, we now consider the local-wellposedness of the PDE-ODE system (\ref{13}) and (\ref{22}) 
as in \cite{GKT2}. \bigskip\par
\vspace{-5pt}
%The initial data here is $(\tq (0), s(0),\al(0)) = (e^{i(m+1)i}q(u_0), s(u_0), \al(u_0))$, 
%but we may assume that $s(u_0)=1$. Indeed, since 
%$\tilde{u}_0:=u(s(u_0)\cdot )$ satisfies $s(\tilde{u}_0)=1$ and $\cae (\tilde{u}_0)=\cae (u_0)$, 
%if we succeed in constructing a solution $\tilde{u}$ to (\ref{1}) with initial data $\tilde{u}_0$, then 
%\begin{equation*}
%u(x,t) := \tilde{u}(\frac{x}{s(u_0)}, \frac{t}{{s(u_0)}^2})
%\end{equation*}
%satisfies (\ref{1}) with $u(0)=u_0$.\par
For time interval $I$, we set $\str (I) := L^\infty_t L^2_x \cap L^4_tL^4_x \cap L^{\sh}_t L^8_x (I)$. 
In \cite{GKT2}, the following proposition is established. 

\begin{prop}\label{P7} (\cite{GKT2})
(i) There exist $\del_0>0$, $\sig>0$, and $C >0$ such that 
the following holds: 
For $(q_0, s_0,\al_0)\in L^2_{\rad}\times \R^+\times \R$ 
%$\del = \sqrt{\cae(u_0)-4\pi m} <\del_0$, 
%then for initial data $(\tq (0), s(0),\al(0)) = (e^{i(m+1)\theta}q(u_0), 1, \al(u_0))$, 
satisfying $\nor{q_0}{L^2}\le \del_0$, 
there exists a solution to the system (\ref{13}) and (\ref{22}); $(q (t), s(t), \al(t))$ on the 
interval $I=[0,\sig s_0^2]$ which satisfies the following properties:
\begin{itemize}
\item $(q (t), s(t), \al(t)) \in C(I; L^2_{\rad}\times \R^+\times\R)$,
\item $(q (0), s(0), \al(0)) = (q_0, s_0, \al_0)$, 
%\item $q(t) = e^{-i(m+1)\theta} \tq (t) \in L^2_{\rm{rad}}$ for all $t\in I$,
\item $\tq\in \str (I)$ and $\nor{\tq}{\str (I)} \le C \nor{q_0}{L^2}$,
\item $s_0^{-1}s(t) \in [0.5,\ 1.5]$ for all $t\in I$.
\end{itemize}
(ii) The solution to (\ref{13}) and (\ref{22}) is unique in 
$(C(I; L^2_{\rad}) \cap \str (I) ) \times C(\R^+\times\R)$.\\
%Namely, if $(q' (t), s'(t), \al'(t)) \in C([0,T']; L^2(\R^2)\times \R^+\times\R)$ is a 
%solution to (\ref{13}) and (\ref{22}) with $T'>0$, then $(\tq (t), s(t), \al(t))=(\tq' (t), s'(t), \al'(t))$ in $[0,\min \{\sig s_0^2, T' \}]$.\\
(iii) There exists $C>0$ such that if $(\tq^i (t), s^i(t), \al^i(t)) \in C(\tilde{I}; L^2(\R^2)\times \R^+\times\R)$ are 
solutions with initial data $(q^i_0,s^i_0,\al^i_0)$ as in (i) for $i=1,2$ 
for $\tilde{I}=[0,\tilde{T}] \subset [0, \min_{i=1,2} \sig ({s_0^i})^2]$, then 
the following difference estimate holds: 
\begin{equation}\label{27}
\begin{aligned}
\nor{q^1-q^2}{\str (\tilde{I})} + &\nor{s^1-s^2}{L^\infty(\tilde{I})} +\nor{\al^1-\al^2}{L^\infty(\tilde{I})} \\
&\hspace{15pt} \le C \left( 
\nor{q^1_0-q^2_0}{L^2} + |s^1_0-s^2_0| + |\al^1_0-\al^2_0| \right) .
\end{aligned}
\end{equation}
\end{prop}

\begin{rmk}\label{R1}
The nonlocal term $N(q)\tq$ can be written as 
\begin{equation}\label{c28}
N(q)\tq = \left( -V(r) + \int_{r}^\infty \frac{2}{r'} V(r') dr\right) \tq ,\quad V(r)= \frac{|q|^2}{2} +\re \frac{m}{r}
\ovl{\nu}  q
\end{equation}
Thus (\ref{13}) makes sense for $q\in \str (I)$ via Duhamel formula. 
\end{rmk}
For the proof of Proposition \ref{P7}, see \cite{GKT2}, Section A.3.
%%%%%%%%%%%%%%%%%%%%%%%Shoumei sita houga ii
%%%%%%%%%%%%%%%%%%%%%%%Detail ha [GKT2] ni mawasite iikara

\subsection{The Proof of (LP1)}

Let us see the proof of (LP1). 
%As mentioned above, we may assume that $s(u_0)=1$. 
%The initial data here is $(\tq (0), s(0),\al(0)) = (e^{i(m+1)i}q(u_0), s(u_0), \al(u_0))$, 
%but 
We may assume that $s(u_0)=1$ by rescaling. %Indeed, since 
%$\tilde{u}_0:=u(s(u_0)\cdot )$ satisfies $s(\tilde{u}_0)=1$ and $\cae (\tilde{u}_0)=\cae (u_0)$, 
%if we succeed in constructing a weak solution $\tilde{u}$ to (\ref{1}) with initial data $\tilde{u}_0$, then 
%\begin{equation*}
%u(x,t) := \tilde{u}(\frac{x}{s(u_0)}, \frac{t}{{s(u_0)}^2})
%\end{equation*}
%is a weak solution to (\ref{1}) with $u(0)=u_0$.\par
By Proposition \ref{P7}, for initial data $(q(u_0), 1, \al(u_0))$, 
there is a unique solution $(\tq (t), s(t), \al(t))$ %\in C(I; L^2(\R^2)\times \R^+\times\R)$ 
to (\ref{13}) and (\ref{22}) for $I=[0,\sig]$ if $\del<\del_0$. 
By Proposition \ref{P5}, we reconstruct $u(t)\in \Sig_m$ from $(\tq (t), s(t), \al(t))$ for each $t\in I$. 
By the continuity of reconstruction, $u(t)\in C(I;\Sig_m)$. 
It suffices to show:\bigskip\par
\vspace{-5pt}
{\bf Claim 1.} $u(t)$ is a weak solution to (\ref{1}). \bigskip\par
\vspace{-5pt}
We prove this claim by approximation by smooth solutions. 
First, 
by Lemma \ref{L11}, %%%%%%%%%%%%%%%%%%%%%%%
we take the following sequence:
\begin{equation*}
\{ u_0^k\}_{k=1}^\infty\subset \Sig_m\cap \dot{H}^3 ,\quad \lim_{k\to\infty} u_0^k = u_0  
\text{ in } \dot{H}^1.
\end{equation*}
Since $s(u_0^k)\to s(u_0)=1$, we may assume that $s(u_0^k)=1$ for all $k\in\bbn$ by rescaling. 
By McGahagan's theorem \cite{M}, there exists a unique solution $u^k(t)\in L^\infty ([0,T_k) 
;\Sig_m\cap \dot{H}^3)$ for each $k\in \bbn$, 
and $\lim_{t\to T_k} \nor{\nabla u^k}{H^2}=\infty$ if $T_k<\infty$. \par
For each $k\in \bbn$ and $t\in [0,T_k)$, 
we denote $(q^k(t),s^k(t),\al^k(t)) = (q(u^k(t)),s(u^k(t)),\al(u^k(t)))$. 
Then Propositions \ref{P3} and \ref{P6} imply that 
$(q^k(t),s^k(t),\al^k(t))$ satisfies (\ref{13}) and (\ref{22}). 
On the other hand, 
applying Proposition \ref{P7} to the initial data $(q^k(0),s^k(0),\al^k(0))= (q(u_0^k), 1, \al(u_0^k))$, 
we have another solution $(q^k_W(t),s^k_W(t),\al^k_W(t))$. %\in C([0,\sig ];L^2\times \R^+\times \R)$. 
%By uniqueness, 
Since both of these belong to the space as in Proposition \ref{P7} (ii), 
we have $(q^k(t),s^k(t),\al^k(t))=(q^k_W(t),s^k_W(t),\al^k_W(t))$ on $t\in [0,\min \{ T_k ,\sig\}]$. Here, we claim that\bigskip\par
\vspace{-5pt}
{\bf Claim 2.} $T_k \ge\sig$. \bigskip\par
\vspace{-5pt}
We first prove Claim 1 provided that Claim 2 holds. 
By Claim 2, $(q^k(t),s^k(t),\al^k(t))$ satisfies all the properties in Proposition \ref{P7}. 
Therefore, by the difference estimate in Proposition \ref{P7} (iii), we have
\begin{equation}
\begin{aligned}
&\nor{q^k-q}{\str (I)} +  \nor{s^k-s}{L^\infty (I)} +\nor{\al^k-\al}{L^\infty (I)} \\
&\hspace{40pt} \lesssim 
\nor{q^k(0)-q_0}{L^2} + |s^k(0)-s(0)| + |\al^k(0)-\al(0)| \to 0
%& \to 0 \text{ as } k\to\infty .
\end{aligned}
\end{equation}
as $k\to \infty$. 
Hence, by the continuity of reconstruction (see Proposition \ref{P5}), 
%\begin{equation}\label{28}
$
\nor{u^k (t) - u (t)}{\dot{H}^1} \to 0$ as $k\to\infty$ for all $t\in I=[0,\sig] 
$. 
%\end{equation}
Since all $u^k(t)$ satisfy (\ref{1}), $u(t)$ is a weak solution to (\ref{1}), which is the desired conclusion. \bigskip\par
We now return to Claim 2. 
While the above argument is established by \cite{GKT2}, 
the proof of Claim 2 is not explicitly described in \cite{GKT2}. 
It is essential to ensure that $T_k$ is bounded from below uniformly in $k$ so that 
the approximation works. 
%This is the first problem which I mentioned in Section 1. 
Hence, we provide a proof of this claim here. \par
Our main ingredients are a priori estimates for modified Schr\"odinger map equation (\ref{13}). 
More precisely, 
\begin{prop}\label{P8}
There exists $\del_0>0$ and $C>0$ such that the following holds: 
If $u(t)\in C(I;\Sig_m)\cap L^\infty(I;\dot{H}^2)$ is a solution to (\ref{1}) on interval $I=(\tau , \tau +\sig)$ 
for $\tau$, $\sig >0$ 
which satisfies $\del := \sqrt{\cae (u) -4\pi m} <\del_0$, then 
\begin{equation}\label{c30}
\begin{aligned}
\nor{\nabla \tilde{q}}{\str (I)} \le C ( \nor{\nabla \tq(\tau)}{L^2_x}
 + (\underline{s}^{-1}& \sig^{1/2} + \underline{s}^{-3/2} \sig^{3/4} \\
& +\nor{q}{L^4_t L^4_x \cap L^{8/3}_t L^8_x}^2) \nor{\nabla \tq}{L^\infty_t L^2_x\cap L^4_t L^4_x} ),
\end{aligned}
\end{equation}
where $\underline{s}= \inf_{t\in I} s(t)$. 
%and 
%$\str (I) = L^{\infty}_tL^2_x(I) \cap L^4_tL^4_x(I) \cap L^{8/3}_t L^8_x(I)$
Furthermore, if $u(t)\in L^\infty(I;\Sig_m\cap \dot{H}^3)$, then
\begin{equation}\label{c31}
\begin{aligned}
\nor{\Del\tq}{\str (I)}\le C ( 
\nor{\Del \tq (\tau )}{L^2_x} + 
 (& \underline{s}^{-3/2} \sig^{3/2} + \nor{q}{L^{8/3}_t L^8_x}^2  + \nor{q}{L^4_t L^4_x}^3 \\
&\hspace{40pt} + \sig^{1/2} \nor{\nabla \tq}{L^4_t L^4_x}
)
 \nor{\Del \tq}{L^\infty_t L^2_x} ).
\end{aligned}
\end{equation}

\end{prop}
The former estimate is obtained in \cite{GKT}, Lemma 3.1. 
What is new is the estimate (\ref{c31}), which is concerned with the second derivative of $\tq$. 
%which plays an important role in deriving a priori bound for 
%$\nor{u(t)}{\dot{H}^3}$. by combining the following estimate: 
This provides a priori bound if we use the following estimate:
\begin{prop}\label{P9}
There exists $\del_0>0$ and $C>0$ such that for $u\in \Sig_m\cap \dot{H}^3$ 
with $\del:=\sqrt{\cae (u) -4\pi m} <\del_0$, we have
\begin{equation}\label{c32}
\nor{u}{\dot{H}^3} \le C (\nor{\Del\tq}{L^2} +\nor{q}{L^6}^2 + s(u(t))^{-2} ).
\end{equation}
\end{prop}
The $\nor{u}{\dot{H}^2}$ counterpart to (\ref{c32}) is obtained in \cite{GKT}, Lemma 4.8. 
%Unlike the case of trivial homotopy type $v(\infty) =-\vec{k}$, %the corresponding estimate is also achieved by \cite{BIKT2}. 
%we have to use $s$ in order to bound the norms concerning $u$. (For the corresponding estimate 
%in the case where $v(\infty)=-k$, see \cite{BIKT2}.) 
%In our setting, since $q$ only indicates the difference between $u$ and the harmonic maps, 
%sole $q$ is not capable of bounding the norm of $u$ which is different from that in \cite{BIKT2}. 
A proof of Propositions \ref{P8} and \ref{P9} is given in Section 5. 
\bigskip\\
\textit{Proof of Claim 2}. 
We omit the index $k$ for simplicity. 
Suppose $T< \sig$. 
We take two numbers $0<\ep <\tau \ll 1$ and set the interval $J=[T-\tau ,T-\ep]$. 
Then, by (\ref{c30}) in Proposition \ref{P8}, we have %%%%%%%%%%%%%%%%%%%%%%%%%%
\begin{equation}\label{29}
\begin{aligned}
\nor{\nabla \tq}{\str (J)} \le C & \left( \nor{\nabla \tq(T-\tau)}{L^2_x} \right. \\
&\hspace{10pt} + (\tau^{1/2} + \tau^{3/4} +
\left. \nor{q}{L^4_t L^4_x \cap L^{8/3}_t L^8_x (J)}^2) \nor{\nabla \tq}{\str (J)} \right) .
\end{aligned}
\end{equation}
By absolute continuity of integral, there exists $\tau_0$ such that for $\tau\le\tau_0$, 
\begin{equation}\label{30}
C (\tau^{1/2} + \tau^{3/4} +
\nor{q}{L^4_t L^4_x \cap L^{8/3}_t L^8_x (J)}^2) < \frac{1}{2} .
\end{equation}
Therefore, by (\ref{29}) and (\ref{30}), for $\tau \le \tau_0$, we obtain 
\begin{equation}\label{31}
\nor{\nabla \tq}{\str (J)} \le 2 C \nor{\nabla \tq(T-\tau)}{L^2_x} .
\end{equation}
Next, by (\ref{c31}) in Proposition \ref{P8}, for $\tau\le\tau_0$, %%%%%%%%%%%%%%%%%%%%%%%%%%
\begin{equation}\label{32}
\begin{aligned}
\nor{\Del\tq}{L^\infty_t L^2_x (J)}\le C ( 
&\nor{\Del \tq (T-\tau )}{L^2_x} \\
&+  ( \tau^{3/2} + \nor{q}{L^{8/3}_t L^8_x(J)}^2  + \nor{q}{L^4_t L^4_x(J)}^3 \\
&\hspace{15pt} + 2 C \tau^{1/2} \nor{\nabla \tq(T-\tau_0)}{L^2_x} %\nor{\nabla \tq}{L^4_t L^4_x}
)
 \nor{\Del \tq}{L^\infty_t L^2_x} ).
\end{aligned}
\end{equation}
(Here, we have applied (\ref{31}) with $\tau=\tau_0$.) 
As in (\ref{31}), 
there exists $\tau_1\le\tau_0$ such that for $\tau\le\tau_1$, 
\begin{equation}\label{33}
C ( \tau^{3/2} + \nor{q}{L^{8/3}_t L^8_x(J)}^2  + \nor{q}{L^4_t L^4_x(J)}^3 
 + 2 C \tau^{1/2} \nor{\nabla \tq(T-\tau_0)}{L^2_x}
)
<\frac{1}{2} .
\end{equation}
By (\ref{32}) and (\ref{33}), for $\tau\le\tau_1$, we obtain 
\begin{equation}\label{34}
\nor{\Del\tq}{L^\infty_t L^2_x (J)}\le 2 C \nor{\Del \tq (T-\tau )}{L^2_x}.
\end{equation}
In particular, from (\ref{31}) and (\ref{34}), 
\begin{equation}\label{35}
\nor{\nabla \tq (T-\ep)}{L^2(\R^2)} \le 2 C \nor{\nabla \tq(T-\tau_1)}{L^2(\R^2)},
\end{equation}
\begin{equation}\label{36}
\nor{\Del\tq (T-\ep)}{L^2(\R^2)}\le 2 C \nor{\Del \tq (T-\tau_1 )}{L^2(\R^2)}.
\end{equation}
Hence, taking limit $\ep \to 0$, we have
\begin{equation}\label{37}
\limsup_{t\to T} \nor{\nabla \tq (t)}{L^2(\R^2)} \le 2 C \nor{\nabla \tq(T-\tau_1)}{L^2(\R^2)} <\infty,
\end{equation}
\begin{equation}\label{38}
\limsup_{t\to T} \nor{\Del\tq (t)}{L^2(\R^2)}\le 2 C \nor{\Del \tq (T-\tau_1 )}{L^2(\R^2)} <\infty.
\end{equation}
Here, Proposition \ref{P9} and Sobolev embedding imply 
\begin{equation}\label{39}
\nor{u(t)}{\dot{H}^3} \le C(\nor{\Del \tq(t)}{L^2} + \nor{\nabla \tq (t)}{H^1}^2 + 1),
\end{equation}
hence
\begin{equation}\label{40}
\limsup_{t\to T} \nor{u(t)}{\dot{H}^3}
\lesssim 
\limsup_{t\to T} \nor{\Del \tq (t)}{L^2} + \limsup_{t\to T}\nor{\nabla \tq(t)}{H^1}^2 + 1
<\infty,
\end{equation}
which leads to a contradiction. Therefore, we obtain $T\ge\sig$.\hfill $\square$

\subsection{Energy Conservation, Regularity, Uniqueness and Continuous Dependence}
We first note that if the solution $u(t)$ is in $C(I;\Sig_m)\cap L^\infty (I;H^2)$, 
then (LP3) can be obtained directly by differentiation with respect to $t$. 
Thus, the energy conservation (LP3) follows from approximation by smooth solution as in the previous subsection. \par
Next, we prove regularity propagation (LP4). 
%The proof is essentially due to \cite{GKT2} and we add some small supplemental arguments. \par
Let $u_0\in \Sig_m\cap \dot{H}^2$ with $\del <\del_0$. 
Without loss of generality, we may assume $s(u_0)=1$. 
Lemma \ref{L5} below implies $q(u_0)\in H^1(\R^2)$. %%%%%%%%%%%%%%%%%
As in the previous subsection, we take a sequence $\{ u_0^k \}_{k=1}^\infty\subset \Sig_m\cap \dot{H}^3$ such that $u_0^k\to u_0$ 
in $\Sig_m$. %(The existence of such sequence is ensured by Lemma \ref{L11} below.) 
We may assume $s(u_0^k)=1$. 
As shown above, for each $k$, there exists a unique solution $u^k(t) \in L^\infty ([0,\sig];\Sig_m\cap \dot{H}^3)$ with $u^k(0)=u_0^k$, where $\sig$ is independent of $k$. 
If we denote $(\tq^k(t), s^k(t), \al^k(t) ):= (\tq(u^k(t)), s(u^k(t)), \al(u^k(t)) )$, Proposition \ref{P3} implies $\tq^k(t) \in L^\infty (I;H^2)$ and 
$\tq^k$ satisfies (\ref{13}). 
Furthermore, %By uniqueness of (\ref{13}) obtained in 
Proposition \ref{P7} implies $\nor{\tq^k}{\str ([0,\sig])} \le C\del_0$ and $s^k(t) \in [0.5,1.5]$ for 
all $k$ (by restricting the sequence to sufficiently large $k$, if necessary). 
Combining these with the estimate (\ref{c30}), 
it follows that $\nor{\tq^k}{L^\infty ([0,\sig]; H^1)}$ is bounded uniformly in $k$. 
On the other hand, Proposition \ref{P7} (iii) implies that 
$\nor{\tq^k -\tq}{L^\infty ([0,\sig];L^2)} \to 0$ as $k\to \infty$. 
Thus, it follows that $\tq \in L^\infty ([0,\sig];H^1)$ (see Proposition 1.4.24 in \cite{CH}). 
By direct differentiation of the integral form of (\ref{13}), we have
\begin{equation*}
\nab \tq (t) = e^{-it\Del} \nab \tq_0 - i \int_{0}^{t} e^{-i (t-\tau) \Del} \nab \Gam (\tq) dt ,
\end{equation*}
where
\begin{equation*}
\Gam (\tq) = \frac{m(1+v_3)(mv_3-m-2)}{r^2}\tq +\frac{mv_{3r}}{r}\tq 
+\tq N(q).
\end{equation*}
Since $\nab \tq_0 \in L^2$ and $\nab \Gam (\tq) \in L^{\sy}_t L^{\sy}_x $, which is the 
consequence of the proof of \cite{GKT}, Lemma 3.1, 
the Strichartz estimates provides $\tq (t) \in C([0,\sig]; H^1)$. 
Hence Lemma \ref{L10}, which is shown later, implies that $u(t)\in C(I;\Sig_m\cap\dot{H}^2)$. 
This is the desired conclusion.\par
%\subsection{Uniqueness and Continuous Dependence}
%Now we prove the uniqueness (LP2)'. 
%We may assume $s(u_0) =1$. 
%For $k=1,2$, let $u^k(t)\in C(I; \Sig_m\cap \dot{H}^2)$ be two solutions to (\ref{1}) 
%with the same initial data $u^k(0)= u_0$ for some interval $I\subset [0,\sig ]$ including $0$. 
%If we construct $(q^k(t), s^k(t), \al^k(t))$ as above for each $k$, 
%both of these satisfy (\ref{13}) and (\ref{22}) by Propositions \ref{P3} and \ref{P6}. 
%Hence by Proposition \ref{P7} (ii), we have $q^1(t) =q^2(t)$ for all $t\in I$. 
%Thus, by uniqueness of reconstruction stated in Proposition \ref{P5}, 
%we obtain $u^1(t) = u^2(t)$ for all $t\in I$, which is the desired conclusion.
%\par
(LP2)' is now immediate from Propositions \ref{P3}, \ref{P5}, \ref{P6} and \ref{P7}. 
%Note that our method cannot be applied in 
%the general case (LP2); i.e., $u(t)\in C(I;\Sig_m)$. 
%This is because the justification of the derivation of (\ref{13}) and (\ref{22}) has not been attained for general solutions in this class. 
%It is worth emphasizing that corresponding regularity of $\tq(t)$ is $C(I;L^2(\R^2))$, 
%for which we are able to formulate (\ref{13}) as we remarked in Remark \ref{R1}.\par %\bigskip\par
(LP5)' also follows immediately from Proposition \ref{P7} (iii), and from Propositions \ref{P5}.
%The uniqueness in this class is quite difficult problem due to serious lack of regularity. 
%We finally show the continuous dependence (LP5)'. Let $u_0\in \Sig_m\cap\dot{H}^2$ 
%with $\del= \sqrt{\cae(u) -4\pi m} < \del_0$, and sequence $\{ u_0^k \}_{k=1}^{\infty}$ with 
%$u_0^k \to u_0$ in $\Sig_m$. Since the map $\Sig_m \ni u\to s_*(u) \in \R$ is continuous, 
%for any number $T< \sig [s_* (u_0)]^2$, we have $[0,T] \subset [0, \sig [s_*(u_0^k)]^2]$ 
%for sufficiently large $k$. 
%We denote the solutions with initial value $u_0$, $u_0^k$ by 
%$u(t)$, $u^k(t)$, respectively. 
%By regularity persistence (LP4), we have $u(t)$, $u^k(t) \in C([0,T]; \Sig_m\cap \dot{H}^2)$. 
%Moreover, when we denote $\tq (t):= \tq (u(t))$ and $\tq^k(t):= \tq (u^k(t))$, 
%$\tq(t)$ and $\tq^k(t)$ satisfies (\ref{13}) by Proposition \ref{P3}. 
%Thus, Proposition \ref{P7} (iii) implies that $\nor{\tq^k(t) -\tq(t)}{L^\infty([0,T] ;L^2)} \to 0$
%as $k\to \infty$. Hence, (LP5)' follows from Proposition \ref{P5}. 

%%%%%%%%%%%%%%%%%%%%%%%%%%%%%%%%%%%%%%%%%%%%%%%%%%%%%%%%%%%%%%%%%%%
%%%%%%%%%%%%%%%
%%%%%%%%%%%%             SECTION 4
%%%%%%%%%%%%%%%
%%%%%%%%%%%%%%%%%%%%%%%%%%%%%%%%%%%%%%%%%%%%%%%%%%%%%%%%%%%%%%%%%%%

\section{Derivation of Modified Schr\"odinger Map}
In this section, we prove Proposition \ref{P3} by 
following the argument in Bejenaru and Tataru \cite{BT}, Chapter 3.
%We emphasize that there are two improvements in our argument. 
%First, we give a justification of calculations when solutions $u(t)$ belong to $C(I;\dot{H}^1)\cap L^\infty (I;\dot{H}^2)$. 
%This regularity is lower than that in \cite{BT}, where they treat the solutions in $C(I;\dot{H}^1\cap\dot{H}^3)$. 
%Note that this improvement of regularity is the essential step in the proof of our uniqueness result (LP2)'. 
%Second, we remove the smallness condition of energy ($\cae(u)-4\pi m \ll 1$), although 
%this improvement brings no benefit on our well-posedness argument. 
We only prove in the case where $u(t)\in C(I;\Sig_m\cap\dot{H}^2)$, 
while the case where $u(t)\in C(I;\Sig_m)\cap L^{\infty} (I;\dot{H}^2)$ is 
almost parallel and achieved by small modifications. 
\par
%\vspace{-1pt}
We introduce some notations here. For a subset $I\subset (0,\infty)$, define 
$$
L^2_e(I) := \left\{ f:(0,\infty)\to \C\ |\ \nor{f}{L^2_e(I)} := \left[ \int_I |f|^2 rdr \right]^{1/2}<\infty \right\} ,
$$
and $L^2_e := L^2_e((0,\infty))$. Note that we can write 
$
\nor{f}{\dot{H}^1_e}^2 = \nor{\rd_r f}{L^2_e}^2 + m^2 \nor{\frac{f}{r}}{L^2_e}^2
$, where $\dot{H}^1_e$ is defined in Section 3. \par
Next, we make a few fundamental observations. For $m$-equivariant maps $u=e^{m\theta R}v(r)$, 
we have the equivalence
\begin{equation}\label{c45}
\nor{u}{\dot{H}^1} \sim \nor{v_1}{\dot{H}^1_e} +\nor{v_2}{\dot{H}^1_e} + \nor{\rd_r v_3}{L^2_e}.
\end{equation}
We can show that the norm $\dot{H}^1_e$ bounds the $L^\infty$-norm, i.e.,
\begin{equation}\label{c46}
\nor{f}{L^\infty} \le C \nor{f}{\dot{H}^1_e}
\end{equation}
for some constant $C>0$. Indeed, for $f\in \dot{H}^1_e$ and $0<r_1<r_2<\infty$, 
\begin{equation}\label{c47}
\begin{aligned}
|f^2 (r_1) -f^2(r_2)| &\le \frac{1}{2} \int_{r_1}^{r_2} |f(r)|\cdot |\rd_r f(r')| dr \\
&\le \frac{1}{2} \nor{\rd_r f}{L^2_e((r_1,r_2))}^2 \nor{\frac{f}{r}}{L^2_e((r_1,r_2))}^2, 
\end{aligned}
\end{equation}
which implies $\lim_{r\to\infty} f(r)$ exists, and it must be $0$ since $\frac{f}{r}\in L^2_e$. 
Taking the limit as $r_2\to \infty$, we have (\ref{c46}) with $C=2^{-1/2}$. \par
(\ref{c47}) also implies that for any $m$-equivariant map $u=e^{m\theta R}v(r) \in \dot{H}^1$, 
$\lim_{r\to 0} v(r)$ and $\lim_{r\to\infty} v(r)$ exist and in $\{\pm \vec{k} \}$, as we remarked before.\\[-5pt]
%%%%%%%%%%%%%%%%%%%%%%%%%%%%%%%%%%%%%%%%%%%%%%%%%%%%%%%%%%%%%%%%%%
%%%%%%%%%%%%%               Lemma 1               %%%%%%%%%%%%%%%%
%%%%%%%%%%%%%%%%%%%%%%%%%%%%%%%%%%%%%%%%%%%%%%%%%%%%%%%%%%%%%%%%%%
\begin{lem}\label{L1}
(i) Let $u(x)=e^{m\theta R}v(r)\in \dot{H}^1(\R^2; \bbs^2)$ and $v(\infty) =\vec{k}$. 
Then there exists a unique function $\eh (r)\in C((0,\infty) ;\R^3)$ such that 
\begin{itemize}
\item $\eh (r)$ is absolutely continuous on any closed subinterval of $(0,\infty)$.
\item $\lim_{r\to\infty} \eh(r) ={}^t(1,0,0)$.
\item $D_r\eh \equiv \rd_r \eh +(\eh\cdot\rd_r v)v =0$ for almost every $r\in (0,\infty)$. 
\end{itemize}
(ii) Let $u^{(i)}(x)=e^{m\theta R}v^{(i)}(r)\in \dot{H}^1(\R^2;\bbs^2)$, and suppose $v^{(i)}(\infty)=\vec{k}$ 
for $i=1,2$. And let $M>0$ satisfying $\| u^{(i)}\|_{\dot{H}^1}\le M$ for $i=1,2$. 
Then, there exists $C=C(M)$ such that
\begin{equation}\label{a41}
\nor{\eh^{(1)} -\eh^{(2)}}{\dot{H}^1_C} \le C(M) \nor{u^{(1)} -u^{(2)}}{\dot{H}^1} ,
\end{equation}
where for $f=(f_1,f_2,f_3) :(0,\infty)\to \R^3$, 
\begin{equation}\label{a42}
\nor{f}{\dot{H}^1_C} := \nor{\rd_r f}{L^2_e} +\nor{f}{L^\infty} +\nor{\frac{f_3}{r}}{L^2_e}. 
\end{equation}
(iii) Let $I\subset\R$ be an open interval, and let $u(t,r)=e^{m\theta R}v(t,r)\in C(I;\dot{H}^1(\R^2;\bbs^2))$ with $v(t_0,\infty) =\vec{k}$ for some $t_0\in I$. 
If $\eh (t)$ is the function as in (i) corresponding to $u(t)$, then we have 
$\eh (t)\in C(I;\dot{H}^1_C)$.
\end{lem}

\begin{rmk}
We can obtain the same result by replacing $\vec{k}$ by $-\vec{k}$. $D_r$ corresponds to the 
covariant derivative along the curve $v(\cdot) :(0,\infty)\to \bbs^2$. 
\end{rmk}

\noindent\textit{Proof of Lemma \ref{L1}.}
The proof of (i) is the same as that in \cite{BT}, Chapter 3. Since (iii) is immediate from (ii), 
we only prove (ii).\par
In the proof, we use the notation $\del\cdot$ to describe the difference between $j=1$ and $j=2$ 
(for example, $\del v= v^{(1)} -v^{(2)}$). 
We take a partition $0=b_0 <b_1<\cdots <b_{n-1} <b_n=\infty$, $n=n(M)$ satisfying 
\begin{equation}\label{a45}
\max_{i=1,2} \left\{ \nor{\rd_r v^{(i)}}{L^2_e(I_k)} + \nor{\frac{\rd_r v^{(i)}_1}{r}}{L^2_e(I_k)} + \nor{\frac{\rd_r v^{(i)}_2}{r}}{L^2_e(I_k)} \right\} \le \frac{1}{8}
\end{equation}
for $i=1,2$ and $k=1,\cdots ,n$ where $I_k= [b_{k-1}, b_k]$. 
Then for $I=I_k$ and $j=1,2$, we have
\begin{equation}\label{a46}
\begin{aligned}
& \nor{\rd_r \del \eh_j}{L^1(I, dr)}  \\
&\le \left\| 
-[\del\eh \cdot \rd_r v^{(1)}] v^{(1)}_j  
-[\eh^{(2)}\cdot \rd_r \del v] v^{(1)}_j -[\eh^{(2)}\cdot \rd_r v^{(2)}]\del v_j
 \right\| _{L^1(I,dr)} \\
&\le \nor{\del\eh}{L^\infty (I)} \nor{\rd_r v^{(1)}}{L^2_e(I)} \nor{\frac{v^{(1)}_j}{r}}{L^2_e(I)}\\
&\hspace{65pt}+ \nor{|\rd_r \del v|\cdot |v^{(1)}_j|}{L^1(I,dr)} 
+ \nor{|\rd_r v^{(2)}_j|\cdot |\del v_j|}{L^1(I,dr)}\\
&\le \frac{1}{16} \nor{\del \eh}{L^\infty (I)} +R(I) ,
\end{aligned}
\end{equation}
where we set
\begin{equation}\label{a47}
%\begin{aligned}
R(I')= \max_{j=1,2} \left\{ \nor{|\rd_r \del v |\cdot |v^{(1)}_j|}{L^1(I',dr)} %\right. \\
%&\hspace{80pt} \left. 
+ \nor{|\rd_r v^{(2)}_j|\cdot |\del v_j|}{L^1(I',dr)} \right\}
%\end{aligned}
\end{equation}
for interval $I'\subset\R$. 
Next, let $\fh :=J\eh$. 
Then $\fh$ satisfies 
$
\rd_r \fh =-(\fh \cdot \rd_r v) v
$, 
and $\lim_{r\to\infty}\fh (r) = {}^t(0,1,0)$. 
Therefore, $\fh$ possesses the same properties as $\eh$ (except boundary condition). 
Thus, (\ref{a46}) also holds when we replace $\eh$ by $\fh$. 
%\begin{equation}\label{a48}
%\nor{\rd_r \del \fh_j}{L^1(I, dr)} \le \frac{1}{16} \nor{\del \fh}{L^\infty (I)} +R(I) .
%\end{equation}
Here, we consider the interval $I=I_n$. From (\ref{a46}), 
\begin{equation}\label{a49}
\begin{aligned}
|\del \eh_j (r)| \le \int_r^\infty |\rd_r \del \eh_j(r')| dr' 
&\le \nor{\rd_r \del \eh_j }{L^1(I_n,dr)}\\
&\le \frac{1}{16} \nor{\del \eh}{L^\infty (I_n)} +R(I_n)
\end{aligned}
\end{equation}
for $r\in I_n$ and $j=1,2$. Hence we have
\begin{equation}\label{a50}
\nor{\del \eh_j}{L^\infty (I_n)} \le \frac{1}{16} \nor{\del \eh}{L^\infty (I_n)} +R(I_n)
\end{equation}
for $j=1,2$, and (\ref{a50}) also holds when $\eh$ is replaced by $\fh$. 
%Similarly, we obtain 
%\begin{equation}\label{a51}
%\nor{\del \fh_j}{L^\infty (I_n)} \le \frac{1}{16} \nor{\del \fh}{L^\infty (I_n)} +R(I_n).
%\end{equation}
%for $j=1,2$. 
To bound the third components, we use the simple relations 
\begin{equation}\label{a52}
\eh = \fh \times v,\qquad \fh = v\times \eh .
\end{equation}
In particular, we have $\eh_3 = \fh_1 v_2 - \fh_2 v_1$. Thus, from (\ref{a50}), 
\begin{equation}\label{a53}
\begin{aligned}
\nor{\del \eh_3}{L^\infty (I_n)} 
&= \left\|
(\del \fh_1 ) v_2^{(1)} 
+\fh_1^{(2)} (\del v_2) -(\del \fh_2) v_1^{(1)} -\fh_2^{(2)} (\del v_1) \right\|_{L^\infty (I_n)}\\
&\le \sum_{j=1}^2 \nor{\del \fh_j}{L^\infty(I_n)} 
+ \sum_{j=1}^2 \nor{\del v_j}{L^\infty(I_n)} \\
& \le \frac{1}{8} \nor{\del \fh}{L^\infty(I_n)} 
+ 2R(I_n) +2C_0 \nor{\del u}{\dot{H}^1(\R^2)}.
\end{aligned}
\end{equation}
where $C_0=\pi C$ and $C$ is the constant in (\ref{c46}). %%%%%%%%%%%%%%%%%%%%%%%%%%%%%%%%%%%%%%%%
Similarly, we have 
\begin{equation}
%\begin{aligned}
\nor{\del \fh_3}{L^\infty (I_n)} 
\le 
\frac{1}{8} \nor{\del \eh}{L^\infty(I_n)} 
+ 2R(I_n) +2C_0 \nor{\del u}{\dot{H}^1(\R^2)} .
%\end{aligned}
\end{equation}
Hence, 
\begin{equation}
\begin{aligned}
& \nor{\del \eh}{L^\infty(I_n)} + \nor{\del \fh}{L^\infty(I_n)}
 \le \sum_{j=1}^3 \left(
\nor{\del \eh_j}{L^\infty(I_n)} + \nor{\del \fh_j}{L^\infty(I_n)}
 \right)\\
& \le \frac{1}{4} \left(\nor{\del \eh}{L^\infty(I_n)} + \nor{\del \fh}{L^\infty(I_n)}\right) + 8R(I_n) +4C_0\nor{\del u}{\dot{H}^1}, 
\end{aligned}
\end{equation}
which implies 
\begin{equation}\label{a56}
\begin{aligned}
\nor{\del \eh}{L^\infty(I_n)} + \nor{\del \fh}{L^\infty(I_n)} \le \frac{32}{3} R(I_n) + \frac{16}{3} C_0 \nor{\del u}{\dot{H}^1} .
\end{aligned}
\end{equation}
Next, we consider $I= I_{n-1}$. 
For $r\in I_{n-1}$ and for $j=1,2$, 
\begin{equation}
\begin{aligned}
\left| \del \eh_j (r)\right| 
&\le \left| \del \eh_j (b_{n-1})\right| +\int_{r}^{b_{n-1}} \left| \rd_r \del \eh_j (r')\right| dr'\\
&\le \nor{\del \eh_j}{L^\infty (I_n)} + \nor{\rd_r \del \eh_j}{L^1(I_{n-1},dr)}\\
&\le\frac{32}{3} R(I_n) + \frac{16}{3} C_0 \nor{\del u}{\dot{H}^1} +\frac{1}{16} \nor{\del \eh}{L^\infty (I_{n-1})} +R(I_{n-1}), 
\end{aligned}
\end{equation}
where the last inequality comes from (\ref{a46}) and (\ref{a56}). Taking the supremum, we have 
\begin{equation}\label{a58}
%\begin{aligned}
\nor{\del \eh_j}{L^\infty (I_{n-1})} \le 
\frac{1}{16} \nor{\del \eh}{L^\infty (I_{n-1})} + 
\frac{32}{3} R(I_{n-1}\cup I_n) +
\frac{16}{3} C_0 \nor{\del u}{\dot{H}^1}
%\end{aligned}
\end{equation}
for $j=1,2$, and (\ref{a58}) also holds when we replace $\eh$ by $\fh$. 
%\begin{equation}
%\begin{aligned}
%\nor{\del \fh_j}{L^\infty (I_{n-1})} \le 
%\frac{1}{16} \nor{\del \fh}{L^\infty (I_{n-1})} + 
%\frac{32}{3} R(I_{n-1}\cup I_n) +
%\frac{16}{3} C_0 \nor{\del u}{\dot{H}^1}.
%\end{aligned}
%\end{equation}
By using (\ref{a52}),  
\begin{equation}
\begin{aligned}
& \nor{\del \eh_3}{L^\infty (I_{n-1})} \\
&\le \sum_{j=1}^2 \nor{\del \fh_j}{L^\infty(I_{n-1})} 
+ \sum_{j=1}^2 \nor{\del v_j}{L^\infty(I_{n-1})} \\
& \le \frac{1}{8} \nor{\del \fh}{L^\infty (I_{n-1})} 
+ \frac{64}{3} R(I_{n-1}\cup I_n) %\\ 
%&\hspace{120pt} 
+ \left( \frac{32}{3}+2 \right) C_0 \nor{\del u}{\dot{H}^1(\R^2)},
\end{aligned}
\end{equation}
and 
\begin{equation}
\begin{aligned}
&\nor{\del \fh_3}{L^\infty (I_{n-1})}
 \le
\frac{1}{8} \nor{\del \eh}{L^\infty (I_{n-1})} 
+ \frac{64}{3} R(I_{n-1}\cup I_n) \\
&\hspace{140pt} +\left(\frac{32}{3}+2\right) C_0 \nor{\del u}{\dot{H}^1(\R^2)}.
\end{aligned}
\end{equation}
Thus, 
\begin{equation}
\begin{aligned}
& \nor{\del \eh}{L^\infty(I_{n-1})} + \nor{\del \fh}{L^\infty(I_{n-1})}
 \le \sum_{j=1}^3 \left(
\nor{\del \eh_j}{L^\infty(I_{n-1})} + \nor{\del \fh_j}{L^\infty(I_{n-1})}
 \right)\\
& \le \frac{1}{4} \left(\nor{\del \eh}{L^\infty(I_{n-1})} + \nor{\del \fh}{L^\infty(I_{n-1})}\right) \\
& \hspace{90pt} 
+ 8\cdot \frac{32}{3}R(I_{n-1}\cup I_n) +9\cdot \frac{16}{3}C_0\nor{\del u}{\dot{H}^1}, 
\end{aligned}
\end{equation}
which implies
\begin{equation}
%\begin{aligned}
 \nor{\del \eh}{L^\infty(I_{n-1})} + \nor{\del \fh}{L^\infty(I_{n-1})} \le \left( \frac{32}{3} \right)^{2} +12 \cdot \frac{16}{3}C_0 \nor{\del u}{\dot{H}^1}. 
%\end{aligned}
\end{equation}
We repeat the above argument, then for $k=1,\cdots n$, we have 
\begin{equation}
\begin{aligned}
& \nor{\del \eh}{L^\infty(I_{k})} + \nor{\del \fh}{L^\infty(I_{k})}\\
& \hspace{30pt} \le \left( \frac{32}{3} \right)^{n+k-1} R([b_{k-1} ,\infty )) + \frac{16}{3}\cdot \left( 12 \right)^{n-k} C_0 \nor{\del u}{\dot{H}^1}. 
\end{aligned}
\end{equation}
Therefore,
\begin{equation}\label{a65}
\begin{aligned}
& \nor{\del \eh}{L^\infty (0,\infty)} = \max_{k=1,\cdots ,n} \nor{\del \eh}{L^\infty (I_k)}\\
& \le 
\left( \frac{32}{3} \right)^{n} R((0 ,\infty )) + \frac{16}{3}\cdot \left( 12 \right)^{n-1} C_0 \nor{\del u}{\dot{H}^1}
%&= \left( \frac{32}{3} \right)^{n} \max_{j=1,2} \left\{ \nor{|\rd_r \del v|\cdot |v_j^{(1)}|}{L^1(0,\infty)} + \nor{|\rd_r v^{(2)}|\cdot |\del v_j|}{L^1(0,\infty)} \right\} \\
%&\hspace{180pt} + \frac{16}{3}\cdot \left( 12 \right)^{n-1} C_0 \nor{\del u}{\dot{H}^1}\\
\le C(M) \nor{\del u}{\dot{H}^1} .
\end{aligned}
\end{equation}
(\ref{a65}) also holds when $\eh$ is replaced by $\fh$.
%Similarly, we have
%\begin{equation}
%\nor{\del \fh}{L^\infty (0,\infty)} \le C(M) \nor{\del u}{\dot{H}^1}.
%\end{equation}
As a consequence of the $L^\infty$ bound, 
using the relation $v =\eh \times \fh$, we obtain
\begin{equation}
\nor{\del v_3}{L^\infty} \le C(M) \nor{\del u}{\dot{H}^1},  
\end{equation}
which implies
\begin{equation}\label{a68}
\nor{\del u}{L^\infty} \le C(M) \nor{\del u}{\dot{H}^1}.
\end{equation}
Thus, we have
\begin{equation}\label{a69}
\begin{aligned}
& \nor{\rd_r \del \eh_j}{L^2_e}  \\
&\le \left\|
- \left[ \del \eh\cdot \rd_r v^{(1)} \right] v^{(1)}_j  
-\left[ \eh^{(2)}\cdot \rd_r \del v \right] v^{(1)}_j -\left[\eh^{(2)}\cdot \rd_r v^{(2)}\right] \del v_j
\right\|_{L^2_e} \\
&\le \nor{\del \eh}{L^\infty} \nor{\rd_r v^{(1)}}{L^2_e} 
+ \nor{\rd_r \del v}{L^2_e} + \nor{\rd_r v^{(2)}}{L^2_e} 
\nor{\del v_j}{L^\infty} \\
&\le C(M) \nor{\del u}{\dot{H}^1} .
\end{aligned}
\end{equation}
On the other hand, since we have $\frac{\eh_3}{r} = \fh_1 \frac{u_2}{r} - \fh_2 \frac{u_1}{r}$ by (\ref{a52}), similar argument yields 
\begin{equation}
\nor{\frac{\del \eh_3}{r}}{L^2_e} \le C(M) \nor{\del u}{\dot{H}^1},
\end{equation}
which is the desired conclusion. 
%(iii) is immediate from (ii). 
\hfill$\square$
\bigskip

Here, let us introduce a new function space. 
Set Hilbert space $H^1_e := L^2_e \cap \dot{H}^1_e$ with inner product
\begin{equation}
\inp{f}{g}_{H^1_e} := \inp{f}{g}_{L^2_e} + \inp{\rd_r f}{\rd_r g}_{L^2_e} +\inp{\frac{f}{r}}{\frac{g}{r}}_{L^2_e}.
\end{equation} 
Then, we define 
\begin{equation}
H^{-1}_e := (H^1_e)^{*}
\end{equation}
\begin{equation}
\nor{f}{H^{-1}_e} := \sup_{g\in H^1_e,\ \nor{g}{H^1_e}=1} \inp{f}{g}_{H^{-1}_e, H^1_e}.
\end{equation}
%This space 
%corresponds to one derivative of $L^2_e$ function. 
This space is introduced so that $\rd_r f$ and $\frac{f}{r}$ are accommodated for $f\in L^2_e$. 
%In \cite{BIKT2}, 
Such kind of space is also considered in \cite{BIKT2}, which, however, uses $\dot{H}^{-1}_e$ 
defined as the dual space of $\dot{H}^1_e$.
%Bejenaru et al. introduces the space $\dot{H}^{-1}_e$, which is the dual space of $\dot{H}^1_e$. 
The reason why we take account of the $L^2_e$-finiteness in our definition is 
the compatibility with the space $H^{-1}(\R^2)$, which will be seen in Lemma \ref{L7}. \par
Although $H^{-1}_e$ penetrates the space of distribution, some calculations such as 
differentiation and multiplication 
are permitted under certain restrictions. 

%%%%%%%%%%%%%%%%%%%%%%%%%%%%%%%%%%%%%%%%%%%%%%%%%%%%%%%%%%%%%%%%%%
%%%%%%%%%%%%%               Lemma 2               %%%%%%%%%%%%%%%%
%%%%%%%%%%%%%%%%%%%%%%%%%%%%%%%%%%%%%%%%%%%%%%%%%%%%%%%%%%%%%%%%%%
\begin{lem}\label{L2}
\begin{enumerate}[(i)]
\item For $f\in L^2_e$, we can define $\rd_r f \in H^{-1}_e$ and $\frac{f}{r} \in H^{-1}_e$ by 
\begin{equation}\label{a74}
\inp{\rd_r f}{\vph}_{H^{-1}_e, H^1_e} := - \inp{f}{\rd_r \vph + \frac{\vph}{r}}_{L^2_e},
\end{equation}
\begin{equation}\label{a75}
\inp{\frac{f}{r}}{\vph}_{H^{-1}_e, H^1_e} := - \inp{f}{\frac{\vph}{r}}_{L^2_e}.
\end{equation}
Moreover, we have
\begin{equation}\label{a76}
\nor{\rd_r f}{H^{-1}_e}\le C \nor{f}{L^2_e},\quad \nor{\frac{f}{r}}{H^{-1}_e}\le \nor{f}{L^2_e}.
\end{equation}
\item Let $g\in L^\infty ((0,\infty))$ with $\rd_r g\in L^2_e$. Then, for $f\in H^{-1}_e$, we can define 
$gf\in H^{-1}_e$ by 
\begin{equation}\label{a77}
\inp{gf}{\vph}_{H^{-1}_e, H^1_e} := \inp{f}{g\vph}_{H^{-1}_e, H^1_e}.
\end{equation}
Moreover, we have
\begin{equation}
\nor{gf}{H^{-1}_e} \le C \nor{f}{H^{-1}_e} \left( \nor{g}{L^\infty} + \nor{\rd_r g}{L^2_e} \right).
\end{equation}
\item For $f\in L^2_e$ and $g\in L^\infty ((0,\infty))$ with $\rd_r g\in L^2_e$, the Leibniz rule holds for $fg$. 
Namely, $\rd_r (fg) = \rd_rf\cdot g + f\cdot\rd_r g$.
\end{enumerate}
\end{lem}

\proof 
(i) is immediate. For $g\in L^\infty$ with $\rd_r g\in L^2_e$ and $\vph\in H^1_e$, 
\begin{equation}
\nor{g\vph}{H^1_e} \le C \left( \nor{g}{L^\infty} +\nor{\rd_r g}{L^2_e} \right) \nor{\vph}{H^1_e},
\end{equation}
which leads to (ii). (iii) follows easily from (i) and (ii).\qedhere

%\begin{rmk}
%The calculation laws (\ref{a74}), (\ref{a75}), and (\ref{a77}) coincide with those in distribution $\mathcal{D}' ((0,\infty))$.
%\end{rmk}
%

%%%%%%%%%%%%%%%%%%%%%%%%%%%%%%%%%%%%%%%%%%%%%%%%%%%%%%%%%%%%%%%%%%
%%%%%%%%%%%%%               Lemma 3               %%%%%%%%%%%%%%%%
%%%%%%%%%%%%%%%%%%%%%%%%%%%%%%%%%%%%%%%%%%%%%%%%%%%%%%%%%%%%%%%%%%
\begin{lem}\label{L3}
In addition to the conditions in Lemma \ref{L1} (iii), 
we further assume that $u(t,x)$ is weakly differentiable with respect to $t$ 
(i.e., $\rd_t u\in L^1_\loc (I\times\R^2)$), and 
that $\rd_t u \in C(I;L^2(\R^2))$. Then, 
\begin{enumerate}[(i)]
\item $u(t,x) - u(t_0,x) \in C^1 (I;L^2(\R^2))$ for any fixed $t_0\in I$.
\item $\rd_r u(t)\in C^1 (I;H^{-1}_e)$ and 
\begin{equation}
\rd_t \rd_r u(t) = \rd_r \rd_t u(t).
\end{equation}
\item %If $\sup_{t\in I} \nor{u}{\dot{H}^1} <\infty$, then 
$\eh (t) - \eh(t_0) \in C^1 (I;L^2_e)$ for any fixed $t_0\in I$. 
Furthermore, if 
$\sup_{t\in I} \nor{u}{\dot{H}^1} <\infty$, then
\begin{equation}
\nor{\rd_t \eh}{L^\infty (I;L^2_e)} \le C \nor{\rd_t u}{L^\infty (I; L^2)}
\end{equation}
%\begin{equation}
%\nor{\rd_t\eh (t) -\rd_t\eh(s)}{L^2_e} \le C \nor{\rd_t v(t) - \rd_t v(s)}{L^2_e},
%\end{equation}
for some $C=C(\sup_{t\in I} \nor{u}{\dot{H}^1})>0$.
\end{enumerate}
\end{lem}

\proof 
(i) By assumptions, the map $t\mapsto u(t,x)$ is absolutely continuous for almost every $x\in\R^2$, 
and its differentiation $\frac{\rd u}{\rd t} (t,x)$ coincides with the weak derivative of $u$ with respect to $t$ (See \cite{GT}, Problem 7.8, for example). 
Therefore, for $s,t\in I$,
\begin{equation}\label{a82}
u(t,x) - u(s,x) = \int_s^t \rd_t u(t',x) dt' \quad\text{ for almost every } x\in\R^2,
\end{equation}
which implies $u(t)- u(s) \in C^1 (I ;L^2(\R^2))$. \par
(ii) By pointwise relation (\ref{a82}), we have
\begin{equation}\label{a83}
v(t,r) - v(s,r) = \int_s^t \rd_t v(t',r) dt' \quad\text{ for almost every } r\in (0,\infty).
\end{equation}
We now interpret (\ref{a83}) as a relation in $L^2_e$. (The right hand side is considered as the Bochner integral.) 
Since $\rd_r : L^2_e\to H^{-1}_e$ is bounded, we have
\begin{equation}
\rd_r v(t) - \rd_r v(s) = \int_s^t \rd_r \rd_t v(t') dt', %\quad\text{ for almost every } x\in\R^2 ,
\end{equation}
which implies (ii).\par
(iii) %First, we show $\eh(t,r) -\eh(s,r)\in L^2_e$ for $s,t\in I$.
We may replace $I$ by an interval contained compactly in $I$, and thus 
may assume $M:=\sup_{t\in I}\nor{u}{\dot{H}^1}<\infty$. 
To take care of integrability, we perform the proof via approximations. 
%We take $I'\Subset I$ (the interior of $I$). 
Let $I'\Subset I$ be arbitrary interval. Applying Lemma \ref{L11}, 
we take a sequence $\{u_n(t)\}_{n=1}^\infty\subset C^\infty (I'\times\R^2)$ satisfying (A), 
(B), and (C)' in Lemma \ref{L11}. 
%Let $\{u_n(t)\}_{n=1}^\infty$ be a sequence as in Lemma \ref{L11},  
And then let $\eh_n (t)$ be the function in Lemma \ref{L1} (i) corresponding to $u_n(t)$. 
The elementary ODE theory yields that 
$\eh_n (t,r)\in C^\infty (I'\times (0,\infty))$, and that 
there exists $R_n>0$ such that if $r\ge R_n$, then $\eh_n(t,r) = {}^t(h_3(r) ,0, -h_1(r) )$ for all $t\in I'$. 
%The former fact is the consequence of elementary ODE theory, and the latter comes from the uniqueness. 
Furthermore, (\ref{a41}) implies
\begin{equation}\label{a85}
\sup_{t\in I'} \nor{\eh (t) - \eh_n(t)}{\dot{H}^1_C} \le C(M) \sup_{t\in I'} \nor{u(t)-u_n(t)}{\dot{H}^1} \to 0 
\ \ \text{ as } n\to\infty .
\end{equation}
%where $M= \sup_{t\in I} \nor{u}{\dot{H}^1}$. 
In particular, there exists $C_1(M)>0$ such that 
$\sup_{t\in I'} \nor{\rd_r u_n (t)}{L^2_e} + \sup_{t\in I'} \nor{\rd_r \eh_n (t)}{L^2_e} \le C_1(M)$ for sufficiently large $n$. \par
Here, we use the following inequality:
%%%%%%%%%%%%%%%%%%%%%%%%%%%%%%%%%%%%%%%%%%%%%%%%%%%%%%%%%%%%%%%%%%
%%%%%%%%%%%%%               Lemma 4               %%%%%%%%%%%%%%%%
%%%%%%%%%%%%%%%%%%%%%%%%%%%%%%%%%%%%%%%%%%%%%%%%%%%%%%%%%%%%%%%%%%
\begin{lem}\label{L4}
For $R_0>0$ and $f: (0,\infty) \to \C$, the following holds: 
\begin{equation}
\nor{\int_r^\infty f(r') dr'}{L^2_e(R_0,\infty)} \le \nor{f}{L^1_e(R_0,\infty)}
\end{equation}
\end{lem}
For the proof, see Lemma 4.1 in \cite{GK}.\bigskip\par
To the moment, we write $\eh, u$ instead of $\eh_n, u_n$ for abbreviation, respectively. 
We take $t,s\in I'$ and use the notation $\del\cdot$ as the difference between the values at $t$ and at $s$  
($\del\eh = \eh (t)-\eh(s)$, for example). 
Then 
%\begin{equation}\label{a87}
%\del \eh = 
%\int_r^\infty - (\del \eh \cdot \rd_r v(t)) v(t) + 
%(\rd_r \eh (s) \cdot \del v) v(t) + (\rd_r \eh (s) \cdot v(s)) \del v\ dr'.
%\end{equation}
\begin{equation}\label{a87}
\begin{aligned}
\del \eh &= -
\int_r^\infty (\del \eh \cdot \rd_r v(t)) v(t) + 
(\eh (s) \cdot \del \rd_r v) v(t) + ( \eh (s) \cdot \rd_r v(s)) \del v\ dr' \\
&= [ (\eh (s) \cdot \del v) v(t) ](r)\\
&\hspace{12pt} + 
\int_r^\infty \left[ -(\del \eh \cdot \rd_r v(t)) v(t) 
+ (\rd_r \eh (s) \cdot \del v) v(t) \right. \\
&\hspace{65pt} \left.+ (\eh (s) \cdot \del v) \rd_r v(t) 
- ( \eh (s) \cdot \rd_r v(s)) \del v \right] \ dr'
\end{aligned}
\end{equation}
where we have used integration by part. By Lemma \ref{L4}, for $R>0$, we have 
\begin{equation}\label{a88}
\begin{aligned}
\nor{\del \eh}{L^2_e (R,\infty)} \le & \nor{|\del \eh|\cdot |\rd_r v(t)|}{L^1_e (R,\infty)} \\
&+
\left( 1 +\nor{\rd_r \eh (s)}{L^2_e} + \nor{\rd_r v(s)}{L^2_e} + \nor{\rd_r v(t)}{L^2_e} \right)\nor{\del v}{L^2_e}.
\end{aligned}
\end{equation}
Here, we take a partition $0=b_0<b_1<\cdots <b_{k-1} <b_k=\infty$, $k=k(M)$, such that 
$\nor{\rd_r v(t)}{L^2_e (I_k)} \le \frac{1}{2}$. 
(Note that $k$ is independent of $n$.) \par
When $R=b_{k-1}$, (\ref{a88}) and H\"{o}lder's inequality give 
\begin{equation}
\nor{\del \eh}{L^2_e (I_k)} \le \frac{1}{2} \nor{\del \eh}{L^2_e (I_k)} + (1+3C_1 (M)) \nor{\del v}{L^2_e}, 
\end{equation}
which implies 
\begin{equation}
\nor{\del \eh}{L^2_e (I_k)} \le 2(1+3C_1 (M)) \nor{\del v}{L^2_e}.
\end{equation}
When $R=b_{k-2}$, (\ref{a88}) implies 
\begin{equation}
\nor{\del\eh}{L^2_e(I_{k-1}\cup I_k)} \le \sum_{j=k-1}^k \frac{1}{2} \nor{\del \eh}{L^2_e (I_j)} +(1+3C_1(M)) \nor{\del v}{L^2_e},
\end{equation}
which implies
\begin{equation}
\nor{\del \eh}{L^2_e (I_{k-1})} \le 4(1+3C_1 (M)) \nor{\del v}{L^2_e}.
\end{equation}
Repeating this argument and undoing the abbreviation, we obtain 
\begin{equation}\label{a93}
\nor{\eh_n (t) -\eh_n(s)}{L^2_e} \le C_2(M) \nor{\del v_n}{L^2_e} \le C_2(M) \int_{s}^t \nor{\rd_t v_n(t')}{L^2_e} dt',
\end{equation}
for some constant $C_2(M)$.
%, which is bounded with respect to $n$. By (\ref{a85}) and the 
%property (C)' in Lemma \ref{L11}, we have
It follows from (\ref{a85}) and (\ref{a93}) that 
$\eh (t) - \eh(s)$ in $L^2_e$ %for all $t,s\in I'$, and 
and $\eh_n(t)-\eh(s) \rightharpoonup \eh (t) - \eh(s)$ in $L^2_e$. 
Moreover, (\ref{a93}) gives
\begin{equation}
\begin{aligned}
\nor{\eh (t) - \eh(s)}{L^2_e} &\le \liminf_{n\to\infty} \nor{\eh_n (t) - \eh_n(s)}{L^2_e} \\
& \le C_2(M) \liminf_{n\infty} \int_{s}^t \nor{\rd_t v_n(t')}{L^2_e} dt' \\
& = C_2(M) \int_{s}^t \nor{\rd_t v(t')}{L^2_e} dt'.
\end{aligned} 
\end{equation}
Since $I'\Subset I$ is arbitrary, this holds for all $t,s\in I$. 
This implies $\eh(t) - \eh(t_0) \in W^{1,\infty} (I,L^2_e)$ for any fixed $t_0\in I$ 
(See Theorem 1.4.40 in \cite{CH}, for example), 
and $\nor{\rd_t \eh}{L^\infty (I;L^2_e)} \le C(M) \nor{\rd_t v}{L^\infty (I;L^2_e)}$.\par
%%%%%%%%%%%%%%%%%%%%%%%%%%%%%%%
Next, we check the continuity of $\rd_t \eh$. 
%Let $t\in I$ and $h\in\R$ with $t+h\in I$. 
%Substituting $t,s$ in (\ref{a87}) for $t+h$, $t$, 
Since it has turned out that $\eh (t) -\eh (s)\in L^2_e$ for all $t,s\in I$, 
it follows that $v$ and $\eh$ satisfy (\ref{a87}). 
Dividing (\ref{a87}) by $t-s$ and 
taking the limit as $s\to t$, we obtain
\begin{equation}\label{a95}
\begin{aligned}
\rd_t \eh (t) 
&= [ (\eh (t) \cdot \rd_t v(t)) v(t) ](r)\\
&\hspace{10pt} + 
\int_r^\infty \left[ -(\rd_t \eh(t) \cdot \rd_r v(t)) v(t) 
+ (\rd_r \eh (t) \cdot \rd_t v(t)) v(t) \right. \\
&\hspace{40pt} \left.+ (\eh (t) \cdot \rd_t v(t)) \rd_r v(t) 
- ( \eh (t) \cdot \rd_r v(t)) \rd_t v(t) \right] \ dr'
\end{aligned}
\end{equation}
%where $\Del_h \eh = \frac{\eh(t+h)-\eh(t)}{h}$ and so on. 
for almost all $t\in I$ and for all $n\in\bbn$. 
It also holds that $v_n$ and $\eh_n$ satisfy 
(\ref{a95}) when $v$, $\eh$ are replaced by $v_n$, $\eh_n$, respectively. 
By taking the difference, the same argument provides 
\begin{equation}\label{d460}
\nor{\rd_t\eh_n (t) -\rd_t\eh_{n'}(t)}{L^2_e} \le C(M) \nor{\rd_t v_n(t) - \rd_t v_{n'}(t)}{L^2_e},
\end{equation}
\begin{equation}\label{d461}
\nor{\rd_t\eh (t) -\rd_t\eh_n(t)}{L^2_e} \le C(M) \nor{\rd_t v(t) - \rd_t v_n(t)}{L^2_e},
\end{equation}
for almost all $t$ and for sufficiently large $n,n'$. (Here, $C(M)$ is a constant independent of $t$ and $n$.) 
Hence by the property (C)' in Lemma \ref{L11}, we have $\eh_t \in C(I;L^2_e)$.\qedhere
\bigskip

We now turn to consider a solution to (\ref{1}); 
$u(t)\in C(I;\Sig_m\cap\dot{H}^2)$, where $I\subset\R$ is an open interval. 
Since (\ref{a68}) implies $u(t)\in C(I;L^\infty(\R^2))$, 
it follows from (\ref{1}) that $\rd_t u(t) \in C(I; L^2(\R^2))$, which enables us to apply Lemma \ref{L3}.\par
Here, note that a map $u\in \Sig_m$ belongs to $\dot{H}^2$ if and only if $\Del u\in L^2_e$, hence
\begin{equation}
H_m v_j, \quad H_0 v_3 \in L^2_e \quad (j=1,2 )
\end{equation}
where
\begin{equation}\label{a99}
H_k :=\rd_{rr} +\frac{1}{r}\rd_r -\frac{k^2}{r^2} 
\quad\text{for } k\in \bbn.
\end{equation}
It is known that there are equivalences 
\begin{equation}\label{a100}
\nor{H_k f}{L^2_e} \sim \nor{\rd_{rr} f}{L^2_e} + \nor{\frac{\rd_r f}{r}}{L^2_e}
+ \nor{\frac{f}{r^2}}{L^2_e} \quad\text{for } k\ge 2,
\end{equation}
\begin{equation}\label{a101}
\nor{H_1 f}{L^2_e} \sim \nor{\rd_{rr} f}{L^2_e} + \nor{\frac{\rd_r f}{r}- \frac{f}{r^2}}{L^2_e} ,
\end{equation}
%\begin{equation}
%\nor{H_0 f}{L^2_e} \sim \nor{\rd_{rr} f}{L^2_e} + \nor{\frac{\rd_r f}{r}}{L^2_e} .
%\end{equation}
We can show these either by direct calculations or by using Hankel transform, for which we refer to \cite{BIKT2}.
%We refer to \cite{BIKT2} for the proof. 
Note that (\ref{a100}) does not hold for $k=1$. Thus, we have to be careful in the $1$-equivariant case. 
Denoting $J =J^{v(r,t)} =v(r,t)\times\cdot $, we define 
\begin{equation}
W(r,t) := \rd_r v - \frac{m}{r} JRv
\end{equation}
\begin{equation}
q(r,t) \equiv q_1(r,t) +i q_2(r,t) := %\left( \rd_r v - \frac{m}{r} JRv \right) 
W \cdot 
\left( \eh +iJ\eh \right),
\end{equation}
\begin{equation}
p(r,t) \equiv p_1(r,t) +i p_2(r,t) := \rd_t v  \cdot \left( \eh +iJ\eh \right),
\end{equation}
\begin{equation}
\nu (r,t) \equiv \nu_1(r,t) +i \nu_2(r,t) := JRv \cdot \left( \eh +iJ\eh \right).
\end{equation}
Besides these notations, we define 
\begin{equation}
\al (r,t) := \rd_t\eh \cdot J\eh = D_t \eh \cdot J\eh,
\end{equation}
or equivalently $D_t \eh = \al J\eh$.\par
The following lemma is concerned with the regularity of each quantity. %The proof is parallel to 
%that in \cite{BIKT2}. 
\begin{lem}\label{L5}
\begin{enumerate}[(i)]
\item $W\in C(I;H^1_e)\cap C^1(I;H^{-1}_e)$.
\item $q\in C(I;H^1_e)\cap C^1(I;H^{-1}_e)$.
\item $p\in C(I;L^2_e)$.
\item $\nu \in C(I; \dot{H}^1_e)$.
\item $\al \in C(I;L^2_e)$.
\end{enumerate}
\end{lem}

\proof 
(iii), (iv), and (v) are immediate from Lemma \ref{L3} and (\ref{a68}). 
Now we prove (i) by following \cite{BIKT2} partially. 
$W\in C(I; L^2_e)$ follows from the definition. 
%We note that $u$ is in $C(I;\dot{H}^2)$ if and only if $\Del u\in L^2_e$, namely
%\begin{equation}\label{a102}
%\rd_{rr} v_j +\frac{1}{r}\rd_r v_j -\frac{m^2}{r^2} v_j,\quad 
%\rd_{rr} v_3 +\frac{1}{r}\rd_r v_3 \in L^2_e \quad (j=1,2).
%\end{equation} 
%When $m\ge 2$, 
When we write  
\begin{equation}\label{a107}
W=
\begin{pmatrix}
W_1\\ W_2\\ W_3
\end{pmatrix}
=
\begin{pmatrix}
\rd_r v_1 +\frac{m}{r}v_3 v_1 \\
\rd_r v_2 + \frac{m}{r}v_3v_2 \\
\rd_r v_3 -\frac{m}{r}(v_1^2+v_2^2)
\end{pmatrix}
,
\end{equation}
we have
\begin{equation} 
\rd_r W_j = v_{jrr} + m v_3 \left( \frac{v_{jr}}{r} - \frac{v_j}{r^2} \right) +\frac{m}{r} v_{3r} v_j \quad \text{for } j=1,2,
\end{equation}
\begin{equation} 
\rd_r W_3 = v_{3rr} +\frac{m}{r^2}(v_1^2+v_2^2) +\frac{2m}{r} v_3v_{3r}.
\end{equation}
%From (\ref{a68}), (\ref{a99}), (\ref{a100}) and (\ref{a101}), we have $\rd_r W \in C(I;L^2_e)$. 
Thus, to derive $\rd_r W \in C(I;L^2_e)$, it suffices to show that 
$\frac{v_1^2+v_2^2}{r^2} \in C(I;L^2_e)$. % even when $m=1$. 
We first note that $\frac{v_1v_{1r} + v_2v_{2r}}{r} = -\frac{v_3 v_{3r}}{r} \in C(I;L^2_e)$. 
Hence, 
\begin{equation}
\frac{v_1^2+v_2^2}{r^2} = \frac{v_1v_{1r} + v_2v_{2r}}{r} - v_1 \left( \frac{v_{1r}}{r} - \frac{v_1}{r^2}  \right) - v_2 \left( \frac{v_{2r}}{r} - \frac{v_2}{r^2} \right) ,
\end{equation}
which is in $C(I; L^2_e)$. 
This also implies $\frac{W_3}{r} \in C(I;L^2_e)$. For $j=1,2$, we can write
\begin{equation} 
\frac{W_j}{r} = \frac{\rd_r v_j}{r} -\frac{v_j}{r^2} +\frac{1+v_3}{r^2}v_j +\frac{m-1}{r^2}v_3v_j .
\end{equation}
In order to derive $\frac{W_j}{r} \in C(I; L^2_e)$, it suffices to show that 
$\frac{1+v_3}{r^2} \in C(I;L^2_e)$. 
Since $v_3 (0, t)=-1$ for all $t\in I$, there exists $R_0$ such that $v_3 (r,t) <0$ for all $r\in (0,R_0]$ 
and for all $t\in I$. Then,
\begin{equation}\label{a112}
\begin{aligned}
\frac{1+v_3}{r^2} &= \frac{1+v_3}{r^2} (\chi_{(0,R_0]} +\chi_{(R_0, \infty )}) \\
&\le \frac{1-v_3^2}{r^2} + \frac{1}{r^2}\chi_{(R_0, \infty )} 
= \frac{v_1^2 +v_2^2}{r^2}+ \frac{1}{r^2}\chi_{(R_0, \infty )},
\end{aligned}
\end{equation}
which implies $\frac{1+v_3}{r} \in L^2_e$ for all $t\in I$. (Here $\chi_A$ represents the characteristic function for some set $A$.) 
The continuity with respect to $t$ also follows easily from the expression in (\ref{a112}). %, taking into account that $\frac{1}{1-v_3}\chi_{(0,R_0]} \in C(I;L^\infty)$. 
Hence, $W\in C(I; H^1_e)$. \par
To derive $W\in C^1(I; H^{-1}_e)$, it suffices to show that $\frac{v_3v_j}{r}$ and $\frac{v_1^2+v_2^2}{r} \in C^1(I;H^{-1}_e)$ for $j=1,2$. 
By Lemma \ref{L2}, it is reduced to showing that $v_j(t)v_k(t) - v_j(t_0)v_k(t_0) \in C^1(I;L^2)$ 
for some fixed $t_0\in I$ and $j,k=1,2,3$, which follows immediately from Lemma \ref{L3}. 
Hence (i) is achieved. \par
We now prove (ii). $q\in C(I;H^1_e)$ is easily derived from (i) and Lemma \ref{L1}. %To see $q\in C^1(I;H^{-1}_e)$, we observe the difference of the value of $q$ between at $t+h$ and at $t$, where $t\in I$ and $h>0$ with $t+h\in I$. By definition, we have
Let $t\in I$ and $h>0$ with $t+h\in I$. By definition, we have
\begin{equation}
\frac{q_1(t+h) - q_1(t)}{h} = 
\frac{W(t+h) -W(t)}{h} \cdot \eh (t+h) + W(t)\cdot \frac{\eh (t+h) -\eh (t)}{h} .
\end{equation}
From (i) and Lemmas \ref{L1} and \ref{L2}, 
the first term tends to $\rd_t W (t) \cdot \eh (t)$. 
The second term converge to $W(t) \cdot \rd_t \eh (t)$ by (i) and Lemma \ref{L3}. 
%Since the argument is similar if $q_1$ is replaced by $q_2$, the desired conclusion follows.
The argument for $q_2$ is similar. Hence $q\in C^1(I;H^{-1}_e)$.\qedhere\bigskip\par
%\vspace{-3pt}
%We have just accomplished the preparation for the derivation of (\ref{13}). 
Let us move on to the derivation of (\ref{13}). 
The outline is based on \cite{CSU}. 
The original equation (\ref{1}) gives, in $L^2_e$ relation, 
\begin{equation}
\begin{aligned}
\rd_t v &= v \times \left( \rd_{rr} v + \frac{1}{r}\rd_r v + \frac{m^2}{r^2}R^2v \right)\\
%& = J \left( D_r W +\frac{1}{r}W + \frac{m}{r} JR W \right)\\
& = J \left( D_r W +\frac{1}{r}W - \frac{m}{r} v_3 W \right) .
\end{aligned}
\end{equation}
This implies 
\begin{equation}\label{a115}
p = i \left( q_r +\frac{q}{r} - \frac{m}{r} v_3 q \right) .
\end{equation}
in $L^2_e$ relation. Next, by Lemmas \ref{L2} and \ref{L3}, we have
\begin{equation}
\rd_t W = \rd_r\rd_t v + \frac{m}{r} \left( \rd_t v_3\, v + v_3\, \rd_t v  \right)
\end{equation}
in $H^{-1}_e$ relation. Hence, considering Lemma \ref{L2} (ii) and (iii), we have
\begin{equation}
%\begin{aligned}
\rd_t W\cdot (\eh + iJ\eh) = \left( \rd_r\rd_t v + \frac{m}{r} v_3 \rd_t v \right) \cdot (\eh +iJ\eh ) 
%&=\rd_r \rd_t v\cdot (\eh +iJ\eh ) + \frac{m}{r}v_3 \rd_t v\cdot (\eh +iJ\eh )\\
=\rd_r p +\frac{m}{r} v_3 p.
%\end{aligned}
\end{equation}
On the other hand, we obtain
\begin{equation}
\rd_t q = \rd_t W\cdot (\eh +iJ\eh ) + W\cdot D_t (\eh +iJ\eh ) 
= \rd_t W\cdot (\eh +iJ\eh ) - i \al q .
\end{equation}
Therefore, 
\begin{equation}\label{a119}
\rd_t q +i\al q = \rd_r p + \frac{m}{r}v_3 p
\end{equation}
in $H^{-1}_e$ relation. 
In turn, direct calculations give 
\begin{equation}\label{a120}
\begin{aligned}
D_r D_t \eh &= \rd_r D_t \eh + (D_t \eh\cdot \rd_r v) v \\
&= \rd_r \rd_t \eh + (\rd_r \eh \cdot \rd_t v) v + (\eh \cdot \rd_r\rd_t v) v +
(\eh\cdot \rd_t v )\rd_r v + (\rd_t \eh \cdot \rd_r v)v
\end{aligned}
\end{equation}
in $H^{-1}_e$ relation. Similarly, we have
\begin{equation}\label{a121}
D_t D_r \eh = (\rd_t\eh\cdot \rd_r v)v + (\eh \cdot \rd_t\rd_r v) v 
+ (\eh\cdot \rd_r v )\rd_t v +(\rd_r \eh \cdot \rd_t v)v
\end{equation}
in $H^{-1}_e$ relation. 
Subtracting (\ref{a121}) from (\ref{a120}), we obtain 
\begin{equation}\label{a122}
D_r D_t \eh - D_t D_r \eh = (\eh\cdot \rd_t v )\rd_r v - (\eh\cdot \rd_r v )\rd_t v,
\end{equation}
where we use the $H^{-1}_e$ identity $\rd_r\rd_t\eh = \rd_t\rd_r\eh $, 
which is shown in the same manner 
as Lemma \ref{L3} (ii). %The same relation holds for $J\eh$. 
On the other hand, we have $D_r \eh =0$ and 
$
D_rD_t \eh = D_r (\al J\eh) = \al_r J\eh
$
 in $H^{-1}_e$. Substituting these for (\ref{a122}), we have
\begin{equation}\label{a123}
\al_r J\eh = (\eh\cdot \rd_t v )\rd_r v - (\eh\cdot \rd_r v )\rd_t v.
\end{equation}
Note that both sides of (\ref{a123}) belong to the class $L^1_e$, while we 
have carried out the calculation in $H^{-1}_e$. 
In particular, it makes sense to consider the inner product of (\ref{a123}) and $J\eh$ in $\R^3$ 
pointwise, which leads to
\begin{equation}\label{a124}
%\begin{aligned}
\al_r = (\eh\cdot \rd_t v ) (\rd_r v \cdot J\eh) - (\eh\cdot \rd_r v ) (\rd_t v \cdot J\eh) 
= \im \left[ \ovl{p} \left( q+\frac{m}{r}\nu \right) \right] .
%\end{aligned}
\end{equation}
By (\ref{a115}), 
\begin{equation}\label{a125}
\al_r = - \re \left[ \left( \ovl{q}+\frac{m}{r} \ovl{\nu} \right) \left( q_r + \frac{1-mv_3}{r}q \right) \right] .
\end{equation}
%Since the right hand side is in $L^1_e$, 
Keeping in mind that (\ref{a125}) is a relation in $L^1_e$, we integrate (\ref{a124}) 
with respect to $r$ and obtain
\begin{equation}
\al (r) = \re \int_r^\infty \left( \ovl{q}+\frac{m}{r} \ovl{\nu} \right) \left( q_r + \frac{1-mv_3}{r}q \right) dr ,
\end{equation}
which coincides with $N(q)$ defined in (\ref{14}). 
On the other hand, by substituting (\ref{a115}) for (\ref{a119}), we have 
\begin{equation}\label{a127}
\rd_t q + i\al q = i\left\{ q_{rr} + \frac{q_r}{r} - \frac{(1-mv_3)^2}{r^2}q - \frac{m}{r}v_{3r} q \right\}
\end{equation}
in $H^{-1}_e$ relation. This is the prototype of (\ref{13}), which represents the relation of the radial components of $\tq$. \par
Here, we need to justify the multiplication of elements in $H^{-1}_e$ by $e^{(m+1)\theta R}$  
in order to derive the relation of $\tq = e^{i(m+1)\theta} q$. 
For $f\in L^2_e$, we define $S(f)\in L^2(\R^2)$ by $S(f) (x) := e^{i(m+1)\theta} f(r)$, 
where $(r,\theta)$ is the polar coordinates of $x$. 
It is trivial that $\nor{S(f)}{L^2(\R^2)} = \sqrt{2\pi}\nor{f}{L^2_e}$. 
% and $\nor{S(f)}{\dot{H}^1(\R^2)} \sim \nor{f}{\dot{H}^1_e}$. 
%We also define the bounded operator $L:H^1_e \to H^{-1}_e$ by $L:= \rd_{rr} +\frac{1}{r}\rd_r - \frac{(m+1)^2}{r^2}$. 
%Here, recall that we have defined the operator $H_{m+1} = \rd_{rr} +\frac{1}{r}\rd_r - \frac{(m+1)^2}{r^2}$. 
Note that $H_{m+1}: H^1_e \to H^{-1}_e$ is bounded.
%%%%%%%%%%%%%%%%%%%%%%%%%%%%%%%%%%%%%%%%%%%%%%%%%%%%%%%%%%%%%%%%%%
%%%%%%%%%%%%%               Lemma 6               %%%%%%%%%%%%%%%%
%%%%%%%%%%%%%%%%%%%%%%%%%%%%%%%%%%%%%%%%%%%%%%%%%%%%%%%%%%%%%%%%%%
\begin{lem}\label{L6}
\begin{enumerate}[(i)]
\item $f \in H^1_e$ if and only if $S(f) \in H^1(\R^2)$, and we have 
$\nor{f}{H^1_e} \sim \nor{S(f)}{H^1(\R^2)}$.
\item $H_{m+1}-1 : H^1_e \to H^{-1}_e$ is invertible.
\item If $f\in L^2_e$, then $S ((H_{m+1}-1)^{-1} f ) \in H^2 (\R^2)$.
\item For $f\in L^2_e$, we have $(\Del -1) S ((H_{m+1}-1)^{-1} f) =S(f)$.
\end{enumerate}
\end{lem}

\proof
(i) Let $f\in H^1_e$. 
We say that a set $F\subset \R^2$ is \textit{fan-shaped} if 
there exist $\theta_0 \in [0,2\pi)$, $\gam \in (0,\pi)$, and two positive numbers $0<r_0<r_1 <\infty$ 
such that $F$ can be written as 
\begin{equation}
F= F_{\theta_0 , \gam}^{r_0,r_1} := \{ (r\cos \theta , r\sin \theta) | r\in (r_0,r_1), \theta - \theta_0 \in (-\gam, \gam) \}.
\end{equation}
For a fan-shaped set $F$, the polar coordinates transformation $\psi:(x,y)\mapsto (r,\theta)$ is 
a $C^\infty$-diffeomorphism. 
%Since $f$ is weak differentiable on $(0,\infty)$ (namely, $f\in W^1(0,\infty)$), we have $S(f) \in W^1 (F)$ and $\rd_x f = \frac{x}{\sqrt{x^2 +y^2}} S (\rd_r f ) - \frac{y}{x^2+y^2} \rd_\theta $.
Since $S(f)\circ \psi^{-1}$ is weakly differentiable on $\psi^{-1}(F)=(r_0, r_1)\times (\theta_0-\gam, \theta_0+\gam)$, 
we have $S(f) \in W^1 (F)$ and 
\begin{equation}\label{a129}
 \rd_x S(f) = \frac{x}{\sqrt{x^2 +y^2}} S (\rd_r f ) - 
i(m+1) \frac{y}{x^2+y^2} 
 S(f) ,
\end{equation}
\begin{equation}\label{a130}
 \rd_y S(f) = \frac{y}{\sqrt{x^2 +y^2}} S (\rd_r f ) + 
i(m+1) \frac{x}{x^2+y^2} 
 S(f) 
\end{equation}
on $F$ (see \cite{GT} for example). Note that this expression does not depend on the choice of $F$. % and 
%the right hand sides can be defined on $(x,y) \in \R^2\bk {0}$.
We define $g_x, g_y :\R^2\to \C$ by
 the right hand sides of (\ref{a129}) and (\ref{a130}) for $(x,y)\in \R^2\bk \{ 0\}$, respectively, 
and $g_x(0)=g_y(0)=0$. 
Obviously, $g_x,g_y \in L^2(\R^2)$, and thus 
$S(f)$, $g_x$ and $g_y$ belong to the space of tempered distributions $\cas'=\cas'(\R^2)$. 
It suffices to show that $\rd_x S(f) =g_x$ and $\rd_y S(f)=g_y$. 
We only observe the former equality. %, since the latter is derived similarly. \par
From the above argument, 
for $\vph \in C^\infty_0(\R^2\bk \{0\})$, % with $\supp \vph \subset \R^2\bk \{ 0\}$, 
we have 
\begin{equation}\label{a131}
\inp{S(f)}{\rd_x \vph}_{\cas',\cas} = -\inp{g_x}{\vph}_{\cas',\cas},\hspace{4pt} 
\inp{S(f)}{\rd_y \vph}_{\cas',\cas} = -\inp{g_y}{\vph}_{\cas',\cas}.
\end{equation} 
Indeed, if we take a finite cover of $\supp \vph$ which consists of fan-shaped sets, and take a 
partition of unity subordinate to it, 
then (\ref{a131}) follows from (\ref{a129}) and (\ref{a130}). 
It further follows that (\ref{a131}) holds for $\vph \in \cas (\R^2)$ with $\R^2\bk (\supp \vph)$ containing a neighborhood of origin. %, then (\ref{a131}) holds for this $\vph$. 
Indeed, let $\eta_0$ be a $C_0^\infty (\R^2)$-function which is $1$ for $|x|\le 1$, $0$ for $|x|\ge 2$, 
and $0\le \eta (x) \le 1$ for all $x\in\R^2$. And let $\eta_j := \eta_0 (\cdot / 2^j) - \eta_0 (\cdot /2^{j+1})$ 
for $j\in \bbn$. Then, for the above $\vph$, we have
%\begin{equation}
$
\left( \sum_{j=0}^{J} \eta_j \right) \vph \to \vph$ in $\cas$ as $J\to \infty
$, 
%\end{equation}
which implies (\ref{a131}). Hence, it follows that $\caf [\rd_x S(f) - g_x ]$ is 
a polynomial, where $\caf$ is the Fourier transform. 
However, since $(1+|\xi|)^{-1}\caf [\rd_x S(f) - g_x ](\xi) \in L^2 (\R^2)$, it must be $0$, which is the desired conclusion. 
%\begin{equation}
%S(f) = g_x + \sum_{\al\in\bbz^2_{\ge 0} , |\al |\le k  } a_\al \rd^\al \del_0
%\end{equation}
%for some $k\in \bbz_{\ge 0}$ and $a_\al \in \C$. 
\par
%%%
%%%
%%%
(ii) We define bilinear form $F: H^1_e \times H^1_e \to \R$ by 
\begin{equation*}
F(f, g):= -\inp{(H_{m+1}-1)f}{g}_{H^{-1}_e ,H^1_e}.
\end{equation*}
Then, $F$ is bounded and coercive, and thus %, and we have
%\begin{equation}
%F(f,f) = \nor{f_r}{L^2_e}^2 + (m+1)^2 \nor{\frac{f}{r}}{L^2_e}^2 +\nor{f}{L^2_e}^2 \ge \nor{f}{H^1_e}^2,
%\end{equation}
%which illustrates the coercivity of $F$. 
%Thus 
$H_{m+1}-1$ is invertible by Lax-Milgram's theorem. \par
(iii), (iv) Let $\psi := (H_{m+1}-1)^{-1} f$. Then we have 
\begin{equation}
\psi_{rr} + \frac{1}{r}\psi_r -\frac{(m+1)^2}{r^2} \psi = f+ \psi \quad\text{in } \mathcal{D}'((0,\infty )).
\end{equation}
Since $f+\psi \in L^2_{\loc} ((0,\infty))$, 
it holds that $\psi \in H^2_{\loc} ((0,\infty))$ (see Theorem 8.8 in \cite{GT} for example). 
%$\psi: (r,\theta) \to x$, $\psi^{-1}:x\to (r,\theta)$. 
%\begin{equation}
%D\psi^{-1} =
%\begin{pmatrix}
%\frac{x}{r} & \frac{y}{r} \\
%-\frac{y}{r^2} & \frac{x}{r^2}
%\end{pmatrix}
%\end{equation}
This implies that $S(\psi) \in W^1(F)$ for every fan-shaped set $F$, and
\begin{equation}
\Del S(\psi) = S(\psi_{rr} + \frac{1}{r}\psi_r -\frac{(m+1)^2}{r^2} \psi ) = S(f) +S(\psi)
\end{equation}
in $F$. By the same argument as (i), it follows that $\caf [\Del S(\psi)- S(f) -S(\psi)]$ is a polynomial. 
From (i), we have $S(\psi)\in H^1(\R^2)$. Thus, $(1+|\xi|)^{-1} \caf [\Del S(\psi)- S(f) -S(\psi)](\xi) \in L^2(\R^2)$, which implies $\Del S(\psi) = S(f)+S(\psi)$ in $\cas'$. 
In particular, $\Del S(\psi)\in L^2(\R^2)$, thus (iii) and (iv) follows. \qedhere
\bigskip\par

Considering Lemma \ref{L6}, for $f\in H^{-1}_e$, we define 
\begin{equation}
S(f) := (\Del - 1) S( (H_{m+1}-1)^{-1} f).
\end{equation}

\begin{lem}\label{L7}
\begin{enumerate}[(i)]
\item $S$ is a bounded linear operator from $H^{-1}_e$ to $H^{-1}(\R^2)$. 
\item For $f\in H^1_e$, we have $\Del S(f) = S(H_{m+1}f)$.
\end{enumerate}
\end{lem}

Lemma \ref{L7} can be expressed by the following commutative diagram:
$$
\begin{CD}
H^1_e @>S >> H^1(\R^2)\\
@VH_{m+1}VV   @VV\Del V\\
H^{-1}_e @>> S> H^{-1}(\R^2)\\
\end{CD} 
$$
This indicates that our definition of $H^{-1}_e$ is a natural radial counterpart of $H^{-1}(\R^2)$. 
\bigskip\\
\textit{Proof of Lemma \ref{L7}.}
(i) is immediate from Lemma \ref{L6}. Let $f\in H^1_e$ and set $\psi := (H_{m+1}-1)^{-1}H_{m+1}f$. 
Then, $\psi = f + (H_{m+1}-1)^{-1}f$. Hence, 
\begin{equation}
\begin{aligned}
S(H_{m+1}f) &= (\Del -1) S(\psi ) \\
&= (\Del -1) S(f) + (\Del -1) S((H_{m+1}-1)^{-1}f) = \Del S(f),
\end{aligned}
\end{equation}
where we use Lemma \ref{L6} (iv). Thus (ii) holds. \hfill $\square$
\bigskip\par
Recall that $\tq (t) =S(q (t))$. 
By Lemmas 5, 6, and 7, $\tq \in C(I; H^1(\R^2))\cap C^1(I;H^{-1}(\R^2))$. 
Operating $S$ on both sides of (\ref{a127}), we obtain 
%\begin{equation}
%i\tq_t = -\Del \tq -\frac{m(1+v_3) (2+m(1-v_3))}{r^2}\tq + \frac{m}{r}v_{3r}\tq + \tq N(q) .
%\end{equation}
%Thus, we 
(\ref{13}).%, and Proposition \ref{P3} is proved.

%%%%%%%%%%%%%%%%%%%%%%%%%%%%%%%%%%%%%%%%%%%%%%%%%%%%%%%%%%%%%%%%%%%
%%%%%%%%%%%%%%%
%%%%%%%%%%%%             SECTION 5
%%%%%%%%%%%%%%%
%%%%%%%%%%%%%%%%%%%%%%%%%%%%%%%%%%%%%%%%%%%%%%%%%%%%%%%%%%%%%%%%%%%

\section{Derivation of Estimates}
\subsection{Proof of Proposition \ref{P8}}
In this section, 
we provide a proof of (\ref{c31}) in Proposition \ref{P8}. 
The proof is an extension of the work in \cite{GKT}, %which consider the regularity up to first 
%derivatives. 
where the a priori estimates for (\ref{13}) %of up to first derivatives of $\tq$ 
are established for up to first spatial derivatives.\par
%However, when $m=1$, there is an obstacle which comes from technical difficulties.
Let $u(t)\in L^\infty (I;\Sig_m\cap \dot{H}^3)$ be a solution to (\ref{1}), 
where $I=(\tau , \tau +\sig)$ for some $\tau>0$ and $\sig>0$. 
Note that $u(t)$ is automatically in $C(I;\Sig_m)$ from the identity (\ref{1}). 
We have checked the regularity of $\tq (t)$ in the previous section 
when $u(t)\in C(I;\Sig_m)\cap L^\infty (I;\dot{H}^2)$. 
In the present case, $\tq (t)$ has additional regularity as follows.

\begin{lem}\label{L8}
$\tq (t) \in L^\infty (I;H^2)\cap W^{1,\infty} (I;L^2)$.
\end{lem}

\proof
We first note that if $u\in \Sig_m\cap \dot{H}^3$, then 
\begin{equation}\label{a139}
\rd_r H_m v_j,\quad \frac{1}{r} H_m v_j ,\quad \rd_r H_0 v_3 \quad  \in L^2_e
\end{equation}
for $j=1,2$, and 
\begin{equation}\label{a140}
\nor{u}{\dot{H}^3} \sim \sum_{j=1}^2 \left( \nor{\rd_r H_m v_j}{L^2_e} + \nor{\frac{1}{r} H_m v_j}{L^2_e} \right) + 
\nor{\rd_r H_0 v_3}{L^2_e} .
\end{equation}
This is immediate from the equivalence $\nor{u}{\dot{H}^3} \sim \nor{\Del u}{\dot{H}^1}$. 
For $W$ in (\ref{a107}) and for $j=1,2$, direct calculations yield
\begin{equation}
\begin{aligned}
H_{m+1} W_j = 
\rd_r &H_m W_j + m \frac{v_3}{r} H_m v_j -2m^2 \frac{1+v_3}{r^3} v_j\\
 &-2m\frac{1+v_3}{r^2} v_{jr} 
+2m v_{3r} \left( \frac{v_{jr}}{r} -\frac{v_j}{r^2} \right) + \frac{m}{r} v_j H_0 v_3
\end{aligned}
\end{equation}
Then, the above representation implies $H_{m+1} W_j \in L^\infty(I;L^2_e)$. Indeed, 
if we regard $v_{jr}$, $\frac{v_j}{r}$, and $H_0 v_3$ as radial symmetric functions in $L^2(\R^2)$, 
Sobolev embedding implies  
all of these three quantities are in $L^\infty(I;L^p_e)$ for $p\in [2,\infty)$. 
Thus it follows from (\ref{a112}) that $\frac{1+v_3}{r^2} \in L^\infty(I;L^p_e)$ for $p\in [2,\infty)$. 
Summarizing these up, we obtain $H_{m+1} W_j \in C(I;L^2_e)$, and hence $(H_{m+1}W_j)\cdot \eh \in L^\infty(I;L^2_e)$. \par
We next show $(H_{m+1} W_3)\cdot \eh_3 \in L^\infty(I;L^2_e)$. 
Direct calculations yield
\begin{equation}
\begin{aligned}
(H_{m+1} F_3)\cdot \eh_3 
= &\rd_r H_0 v_3 \cdot \eh_3 + m(m^2+2m) \frac{1-v_3^2}{r^2} \frac{\eh_3}{r} 
- 2m \frac{v_{3r}}{r} \frac{\eh_3}{r} \\
&+ 2m v_3 v_{3rr} \frac{\eh_3}{r} 
+ 2m v_{3r}^2 \frac{\eh_{3}}{r} -(m^2+2m) \frac{v_{3r}}{r}\frac{\eh_3}{r} .
\end{aligned}
\end{equation}
Thus, it suffices to show that $\frac{\eh_3}{r}$ and $\frac{v_{3r}}{r}$ are in $L^\infty(I;L^p_e)$ 
for $p\in [2,\infty)$. 
The former follows from (\ref{a52}). Since $v_{3rr} +\frac{v_{3r}}{r} \in L^\infty(I;L^p_e)$ 
by the Sobolev inequality, it suffices to show $v_{rr}\in L^\infty(I;L^p_e)$. 
However, this is immediate from change of coordinates; 
$u_{rr} = \frac{x^2}{r^2} u_{xx} + \frac{y^2}{r^2} u_{yy} + \frac{2xy}{r^2}u_{xy}$. 
Hence, $(H_{m+1}W)\cdot \eh \in L^\infty(I;L^2_e)$. 
Note that 
$H_{m+1} q = (H_{m+1}W)\cdot \eh +2 W_r\cdot \eh_r$. 
Since both $W_r$ and $\eh_r$ are in $L^\infty (I;L^4_e)$, 
we have $\tq \in L^\infty (I;H^2)$ from Lemma 7 (iii). \par
$\tq (t) \in W^{1,\infty} (I;L^2)$ follows from the fact that $\tq$ satisfies (\ref{13}) and by Hardy's inequality (see (\ref{c156}) below). \qedhere
\bigskip\par

We move on to the derivation of the estimate (\ref{c31}). 
In the proof, we sometimes write the spaces of radial component like $L^p_e$ as $L^p$ for abbreviation 
if there is no ambiguity. 
Let $x_i$ indicate the $i$-th spatial coordinate of $\R^2$ for $i=1,2$. 
By operating $\rd_{x_i} \rd_{x_j}$ on the equation (\ref{13}) for $i,j=1,2$, we obtain 
\begin{equation}
i U_t + \Del U = \sum_{k=1}^{9} A_k,
\end{equation}
where 
\begin{align*}
& U = \rd_{x_i} \rd_{x_j} \tq,  \\
& A_1 = \frac{m(1+v_3)(mv_3-m-2)}{r^2} U, \\
& A_2 = \left( \frac{m(1+v_3)(mv_3-m-2)}{r^2} \right)_{x_i}\hspace{-2pt} \tq_{x_j} +\left( \frac{m(1+v_3)(mv_3-m-2)}{r^2} \right)_{x_j}\hspace{-2pt} \tq_{x_i}, \\
%& A_3 = \left( \frac{m(1+v_3)(mv_3-m-2)}{r^2} \right)_{x_j} \tq_{x_i}, \\
& A_3 = \left( \frac{m(1+v_3)(mv_3-m-2)}{r^2} \right)_{x_ix_j} \tq ,
\end{align*}
\begin{align*}
& A_4 = \frac{m}{r}v_{3r} U,
&& A_5 = \left( \frac{m}{r}v_{3r} \right)_{x_i} \tq_{x_j} + \left( \frac{m}{r}v_{3r} \right)_{x_j} \tq_{x_i},\\
%& A_7 = \left( \frac{m}{r}v_{3r} \right)_{x_j} \tq_{x_i},
& A_6 = \left( \frac{m}{r}v_{3r} \right)_{x_ix_j} \tq ,
&& A_7 = N(q) U,\\
& A_{8} = (N(q))_{x_i} \tq_{x_j} + (N(q))_{x_j} \tq_{x_i},
%& A_{11} = (N(q))_{x_j} \tq_{x_i},
&& A_{9} = (N(q))_{x_ix_j} \tq . 
\end{align*}
By Strichartz estimates, we have
\begin{equation}
\nor{U}{\str (I)} \le C \left( \nor{U(\tau)}{L^2_x} + \sum_{k=1}^{9} \nor{A_k}{L^{\sy}_{t,x}} \right) .
\end{equation}
To derive the bound for each $A_k$, we make here 
some preparations based on the papers \cite{GKT}, \cite{GKT2}. 
%concern about 
\begin{itemize}
%\item $h_1, 1+h_3$, and its derivatives have good decaying properties. For example, 
%$\nor{\frac{h_1}{r}}{L^\infty}, \nor{\frac{1+h_3}{r^2}}{L^\infty} <\infty. $ Moreover, since
%\begin{equation*}
%$
%(J \bj )_r = - \frac{m}{r}h_1 h
%$ and  
%$
%h_r = \frac{m}{r}h_1 (J \bj ) ,
%$ 
%\end{equation*}
%we have $\nor{(1 +r) (J\bj)_r}{L^\infty}$, $\nor{(1+r) h_r}{L^\infty} <\infty$. 
%Similarly, we have $\nor{(r+r^2)h_{rr}}{L^\infty} <\infty$.
\item If $\del$ is sufficiently small, then 
\begin{equation}
\nor{z}{\dot{H}^1_e} \sim \sqrt{\cae(u) -4\pi m}.% \lesssim \nor{z}{\dot{H}^1_e}
\end{equation}
%from the choice of $(s,\al)$ in Proposition \ref{P4}. 
In particular, $\nor{z}{L^\infty}$ is sufficiently small, and hence we have 
\begin{equation}\label{a146}
|\gam| \lesssim |z|^2,\quad |\gam_r| \lesssim |z||z_r|.%,\quad |\gam_{rr}|\lesssim |z_r|^2 +|z||z_{rr}|. 
\end{equation}
\item For simplicity, we sometimes write $a_1\eh + a_2 J^v\eh$, 
$ a_1 \bj + a_2 J^{h}\bj$ as $a\eh$, $a\bj$ for $a=a_1+ia_2\in\C$, respectively. Under this convention, 
$q$ and $z$ satisfy the relation 
\begin{equation}\label{a147}
s e^{-\al R} (q\eh )(sr) = (L_0 z)\bj + (\gam h)_r +\frac{2m}{r} h_3 \gam h + \frac{m}{r} \xi_3 \xi,
\end{equation}
where %$\eta := e^{m\theta R} (z\bj)$, $L f := \rd_r f + \frac{m}{r}(f_3 h + h_3 f)$ for $f:(0,\infty)\to \R^3$
$L_0 z := z_r +\frac{m}{r} h_3 z$
, and $\xi(r) := e^{-\al R}v(sr) - h(r) = z\bj +\gam h$.
\item To obtain bounds of quantities related to $z$ by those of $q$, 
the following lemma effectively works:%plays an important role to obtain bounds for the norms in terms of $q$: 
%%%%%%%%%%%
% lemma 9 %
%%%%%%%%%%%
\begin{lem}\label{L9}
Let $g(r):(0,\infty)\to \C$ be a measurable function satisfying $g\in L^\infty((0,\infty))$ and 
$g_r\in L^1_\loc ((0,\infty))$. 
And assume that $m\in \bbn$, $p\in [2, \infty)$, and $a\in \R$ satisfy $m-a+\frac{2}{p}>0$ 
and $pa>2$. 
Then, if $\frac{g_r}{r^{a-1}} - \frac{m}{r^a}g \in L^p_e$, %for some $p\in [2,\infty)$ and $m\ge a$, 
we have $\frac{g_r}{r^{a-1}},\frac{1}{r^a}g \in L^p_e$ and
\begin{equation}
\nor{\frac{g_r}{r^{a-1}}}{L^p_e} + \nor{\frac{g}{r^a}}{L^p_e} \le C \nor{\frac{g_r}{r^{a-1}} - \frac{m}{r^a}g}{L^p_e}
\end{equation}
for some constant $C=C(m,p,a)$.
\end{lem}
This is a simple extension of Lemma 3.6 in \cite{CSU} and Lemma 4.2 in \cite{GKT}. %We prove this lemma in Section 5.3. 
The proof of this lemma follows that of Lemma 4.2 in \cite{GKT}, 
replacing $\rd_r$, $\frac{\cdot}{r}$ by $\frac{1}{r^{a-1}}\rd_r$, $\frac{\cdot}{r^{a}}$, respectively. 
\item The following estimate for $z$ is obtained in \cite{GKT}, Lemma 4.8. 
%\begin{lem}
Namely, there exists $\del_0>0$ such that for $u\in \Sig_m$ with $\del= \sqrt{\cae (u)-4\pi m} < \del_0$ 
and for $p\in [2,\infty)$, 
we have
\begin{equation}\label{a149}
\nor{z_r}{L^p_e} +\nor{\frac{z}{r}}{L^p_e} \lesssim s^{1-2/p}\nor{q}{L^p_e} +\nor{q}{L^2_e}
\end{equation}
%\begin{equation}
%\begin{aligned}\label{a150}
%\nor{z_{rr}}{L^p_e} +\nor{\frac{z_r}{r} -\frac{mz}{r^2}}{L^p_e} 
%\lesssim & s^{2-2/p} \left(\nor{q_r}{L^p_e} +\nor{\frac{q}{r}}{L^2_e} + \nor{q}{L^{2p}_e}^2\right) \\
%&+s^{1-2/p}\nor{q}{L^p_e} +\nor{q}{L^2_e}^2 +\nor{q}{L^2_e}.
%\end{aligned}
%\end{equation}
%\end{lem}
%some estimates in terms of $z$ are established. 
%(\ref{a149}) and (\ref{a150}) with $p=2$ are established in \cite{GKT}, and 
%our new ingredient is (\ref{a150}) for general $p$, 
%which we prove in Section ??? below. 
In \cite{GKT}, the estimate for $\nor{z_{rr}}{L^2}$ is also obtained and used to 
find the bound for $\nor{v_{3rr}}{L^2}$. 
However, we slightly modify their method, and 
we only need (\ref{a149}) to estimate the terms which $z$ concerns. %estimates up to first derivatives of $z$. 
%This enables us to deal with the terms containing $v_{3rrr}$. 
\item We observe here several simple estimates which are also seen in \cite{GKT}. 
Let $u\in \Sig_m$ with $\del = \sqrt{\cae(u) -4\pi m}<\del_0$. 
Since $\nor{q}{L^2_e} =\pi^{-1}\del \lesssim 1$, %and since $\del_0$ is sufficiently small, we may have $\nor{q}{L^2_e} \le 1$. %if $\del_0$ is sufficiently small. 
%we have especially $\nor{q}{L^2_e} \le 1$ by using the smallness of $\del_0$.
%Thus, 
(\ref{a146}) and (\ref{a149}) provide
\begin{equation}\label{a150}
\begin{aligned}
\nor{\frac{1+v_3}{r^2}}{L^4} &= s^{-3/2} \nor{\frac{1+h_3 +z_2 h_1 +\gam h_3}{r^2}}{L^4}\\
&\le s^{-3/2} \left( \nor{\frac{1+h_3}{r^2}}{L^4} + \nor{\frac{z_2 h_1}{r^2}}{L^4} +\nor{\frac{\gam h_3}{r^2}}{L^4} \right) \\
&\lesssim s^{-3/2} \left(1 +\nor{\frac{z}{r}}{L^4} +\nor{\frac{z}{r}}{L^8}^2 \right)\\
&\lesssim s^{-3/2} + s^{-1} \nor{q}{L^4} +\nor{q}{L^8}^2.
\end{aligned}
\end{equation}
Similarly, we have
\begin{equation}\label{a151}
\begin{aligned}
\nor{\frac{v_{3r}}{r}}{L^4} \hspace{-1pt} 
&= \hspace{-1pt} \nor{\frac{1}{sr} \hspace{-2pt}
\left[ 
\frac{mh_1^2}{r} + h_1z_{2r} - \frac{mh_1h_3z_2}{r} + \gam_r h_3 + \frac{m\gam h_1^2}{r^2}
\right] \hspace{-2pt} \left( \frac{r}{s} \right)
}
{L^4}\\
&\lesssim s^{-3/2} \left( 1 + \nor{z_r}{L^4} + \nor{\frac{z}{r}}{L^4} + \nor{z_r}{L^8}^2 + \nor{\frac{z}{r}}{L^8}^2 \right)\\
&\lesssim s^{-3/2} + s^{-1} \nor{q}{L^4} + \nor{q}{L^8}^2,
\end{aligned}
\end{equation}
%\begin{equation}
%\nor{\frac{v_{3r}}{r}}{L^4} 
%\lesssim 
%\nor{\frac{W_3}{r}}{L^4} + \nor{\frac{1-v_3^2}{r^2}}{L^4} 
%\lesssim 
%\nor{\frac{q}{r}}{L^4} + \nor{}{}
%\end{equation}
\begin{equation}
\nor{v_{3r}}{L^4} \lesssim s^{-1/2} + \nor{q}{L^4},\qquad \nor{v_{3r}}{L^8} \lesssim s^{-3/4} + \nor{q}{L^8}.
\end{equation}
%\begin{equation}
%\nor{v_{3r}}{L^8} \lesssim s^{-3/4} + \nor{q}{L^8}.
%\end{equation}
Here, the notation in the first line of (\ref{a151}) means 
the composite of the functions in the square brackets and $\frac{r}{s}$. 
%In addition, (\ref{a150}) yields 
%\begin{equation}
%\begin{aligned}
%\nor{v_{3rr}}{L^4}
%&= 
% \left\| \frac{1}{s^2}
%\left[ 
%\left(\frac{mh_1^2}{r}\right)_r 
% + 
%\left( h_1z_{2r}\right)_r 
% - 
%\left( \frac{mh_1h_3z_2}{r}\right)_r \right. \right. \\
%& \qquad\qquad\qquad\qquad \left. \left. + 
%\left( \gam_r h_3\right)_r  + 
%\left( \frac{m\gam h_1^2}{r^2} \right)_r 
%\right] \left( \frac{r}{s} \right) 
%\right\|_{L^4}\\
%&\lesssim s^{-3/2} 
%\left( 
%1 +
%\nor{z_r}{L^4} +
%\nor{\frac{z}{r}}{L^4} + \nor{z_r}{L^8}^2 + \nor{z_{rr}}{L^4}
%\right)\\
%&\lesssim 
%s^{-3/2} + s^{-1} \nor{q}{L^4} + \nor{q}{L^8}^2 + \nor{q_r}{L^4} + \nor{\frac{q}{r}}{L^4}.
%\end{aligned}
%\end{equation}
Since 
\begin{equation}
v_{3rr} = \rd_r W_3 - m\frac{1-v_3^2}{r^2} -2mv_3\frac{v_{3r}}{r},
\end{equation}
we have 
\begin{equation}
\begin{aligned}
\nor{v_{3rr}}{L^4} &\lesssim
 \nor{\rd_r (q\eh)}{L^4} +\nor{\frac{1+v_3}{r^2}}{L^4} +\nor{\frac{v_{3r}}{r}}{L^4} \\
%&\lesssim \nor{q_r}{L^4} + \nor{q v_r}{L^4} +\nor{\frac{1+v_3}{r^2}}{L^4} +\nor{\frac{v_{3r}}{r}}{L^4}\\
&\lesssim \nor{q_r}{L^4} + \nor{\frac{q}{r}}{L^4} + \nor{q}{L^8}^2 
+ s^{-1} \nor{q}{L^4} + s^{-\Na}.
\end{aligned}
\end{equation}
%where we use $v_r = W + \frac{m}{r} JRv$ in the last inequality.
%\item We further observe the bound of $v_{3rrr}$. 
\item The nonlocal terms are treated in the following manner. 
By Hardy's inequality, for $p\in [1,\infty)$ and for $f\in L^p_e$, we have
\begin{equation}\label{c156}
\nor{\int_{r}^\infty \frac{f(r')}{r'} dr'}{L^p_e} \lesssim \nor{f}{L^p_e}.
\end{equation}
Since $|\nu|=\sqrt{1-v_3^2}$, 
(\ref{c28}) and (\ref{a150}) implies
\begin{equation}
\begin{aligned}
\nor{N(q)}{L^4} \lesssim \nor{q}{L^8}^2 +\nor{\frac{\nu}{r}q}{L^4} 
&\lesssim \nor{\frac{1-v_3^2}{r^2}}{L^4} + \nor{q}{L^8}^2 \\
&\lesssim s^{-3/2} + s^{-1}\nor{q}{L^4} + \nor{q}{L^8}^2.
\end{aligned}
\end{equation}
%where we use (\ref{a150}).
\end{itemize}

Now, we derive the bound for each $A_k$. 
\begin{equation}
\begin{aligned}
\nor{A_1}{\La} 
&\lesssim \nor{\frac{1+v_3}{r^2}}{L^{\sy}_t L^4_x} \nor{U}{L^\infty_t L^2_x} \\
&\lesssim
\left( \us^{-\Na} \sig^{\Nc} + \us^{-1}\sig^{\Nb} \nor{q}{L^4_tL^4_x} +\nor{q}{L^{\sh}_t L^8_x}^2  \right) \nor{U}{L_t^\infty L^2_x}\\
& \lesssim
\left( \us^{-\Na} \sig^{\Nc} +\nor{q}{L^{\sh}_t L^8_x}^2 
+ \nor{q}{L^4_tL^4_x}^3 \right) \nor{U}{L_t^\infty L^2_x}
\end{aligned}
\end{equation}
\begin{equation}
\begin{aligned}
\nor{A_2}{\La} %+ \nor{A_3}{\La}\\
&\lesssim 
\nor{\left( \left|\frac{v_{3r}}{r^2}\right| + \left|\frac{1+v_3}{r^3}\right|  \right) \left( 
\left| q_r \right| + \left| \frac{q}{r} \right| \right)}{\La}\\
& \lesssim 
\left( \us^{-\Na} \sig^{\Nc} +\nor{q}{\Lb}^2 + \nor{q}{L^4_tL^4_x}^3 \right)  \nor{\left|\frac{q_r}{r}\right| + \left|\frac{q}{r^2}\right|}{L_t^\infty L^2_x} .
\end{aligned}
\end{equation}
\begin{equation}
\begin{aligned}
&\nor{A_3}{L^{\Nc}_tL^{\sy}_x}  \\
&\lesssim
\left(
\nor{\left|\frac{1+v_3}{r^2}\right| + \left|\frac{v_{3r}}{r}\right| + \left|v_{3rr}\right|}{L^{\sy}_t L^4_x} 
%+ \nor{\frac{v_{3r}}{r}}{L^{\sy}_t L^4_x}
+ \nor{v_{3r}}{L^{\sh}_tL^8_x}^2 
%+ \nor{v_{3rr}}{L^{\sy}_t L^4_x}
\right)
%&\hspace{230pt} \times 
\nor{\frac{q}{r^2}}{L^\infty_t L^2_x}\\
&\lesssim
\left(
\us^{-\Na} \sig^{\Nc} + \nor{q}{\Lb}^2 +\nor{q}{\Ld}^3 \right. \\
&\hspace{110pt} \left. +\sig^{\Nb} (\nor{q_r}{\Ld} + \nor{\frac{q}{r}}{\Ld})   
\right) \nor{\frac{q}{r^2}}{L^\infty_t L^2_x}
\end{aligned}
\end{equation}
\begin{equation}
\begin{aligned}
\nor{A_4}{\La}
&\lesssim \nor{\frac{1+v_3}{r^2}}{L^{\sy}_t L^4_x} \nor{U}{L^\infty_t L^2_x}\\
&\lesssim 
\left( \us^{-3/2} \sig^{3/4} +\nor{q}{L^{\sh}_t L^8_x}^2 +\nor{q}{L^4_tL^4_x}^3 \right) \nor{U}{L_t^\infty L^2_x}.
\end{aligned}
\end{equation}
\begin{equation}
\begin{aligned}
\nor{A_5}{\La} %+\nor{A_7}{\La}\\
&\lesssim 
\nor{\left( \left| v_{3rr} \right| + \left| \frac{v_{3r}}{r} \right| \right) 
\left( \left| \frac{q_r}{r} \right| + \left| \frac{q}{r^2} \right| \right) }{\La} \\
& \lesssim 
\left(
\us^{-\Na} \sig^{\Nc} + \nor{q}{\Lb}^2 +\nor{q}{\Ld}^3 \right. \\
& \hspace{40pt} \left. +\sig^{\Nb} (\nor{q_r}{\Ld} + \nor{\frac{q}{r}}{\Ld})   
\right) \nor{\frac{q}{r^2}}{L^\infty_t L^2_x}
\end{aligned}
\end{equation}
\begin{equation}
\begin{aligned}
\nor{A_7}{\La}
&\lesssim \nor{N(q)}{L^{\sy}_tL^4_x} \nor{U}{\Lc} \\
&\lesssim \left( \us^{-3/2} \sig^{3/4} +\nor{q}{L^{\sh}_t L^8_x}^2 +\nor{q}{\Ld}^3 \right) \nor{U}{L_t^\infty L^2_x}.
\end{aligned}
\end{equation}
\begin{equation}
\begin{aligned}
&\nor{A_8}{\La}\\ %+ \nor{A_{11}}{\La} \\
&\lesssim \nor{\left( q + \frac{m}{r}\nu \right) \left( q_r + \frac{1-mv_3}{q} \right) 
\left( \left| q_r \right| + \left| \frac{q}{r} \right| \right)}{L^{\sy}_tL^4_x}  \\
&\lesssim \sig^{\Nb} 
\hspace{-3pt}
\left( \nor{q_r}{\Ld} \hspace{-1pt} + \nor{\frac{q}{r}}{\Ld} \right) 
\hspace{-3pt}
\left( \nor{q_{rr}}{\Lc} \hspace{-1pt} + \nor{\frac{q_r}{r}}{\Lc} \hspace{-1pt} + \nor{\frac{q}{r^2}}{\Lc} \right) .
\end{aligned}
\end{equation}
\begin{equation}
\begin{aligned}
&\nor{A_9}{\La} \lesssim 
\nor{N(q)_{rr}\tq}{\La} + \nor{\frac{1}{r} N(q)_r \tq}{\La} \\
&\lesssim
\nor{\left( \left| q_r \right| +\left| \frac{\nu_r}{r} \right| + \left| \frac{\nu}{r^2} \right| \right) 
\left(
\left| q_r \right| + \left| \frac{q}{r} \right| 
\right)
\left| q\right|
}{\La}\\
& 
\qquad +\nor{
\left(
\left| q \right| + \left| \frac{\nu}{r} \right|
\right)
\left(
\left| q_{rr} \right| +\left| \frac{q_r}{r} \right| + \left| \frac{q}{r^2} \right| + \left| \frac{v_{3r}}{r} q \right|
\right)
\left| q \right|
}{\La} \\
&\qquad 
+\nor{\left( q + \frac{m}{r}\nu \right) \left( q_r + \frac{1-mv_3}{q} \right) 
 \left| \frac{q}{r} \right| }{L^{\sy}_tL^4_x} \\
&\lesssim
\left(
\us^{-\Na} \sig^{\Nc}  
+ \nor{q}{\Lb}^2 + \nor{q}{\Ld}^3 +\sig^{\Nb} (\nor{q_r}{\Ld} + \nor{\frac{q}{r}}{\Ld})   
\right) \\
& \hspace{95pt} \times \left( \nor{q_{rr}}{L^\infty_t L^2_x} + \nor{\frac{q_r}{r}}{L^\infty_t L^2_x} + \nor{\frac{q}{r^2}}{L^\infty_t L^2_x} \right)
\end{aligned}
\end{equation}
The remaining term is $A_6$. First, we have
\begin{equation}
\nor{A_6}{\La}
\lesssim \nor{\frac{v_{3rrr}}{r}q}{\La} + \nor{\frac{v_{3rr}}{r^2}q}{\La} +\nor{\frac{v_{3r}}{r^3}q}{\La} .
\end{equation}
The second and third terms are estimated in the same manner above. 
Hence, it suffices to find the bound of the first term. 
By direct calculations, 
\begin{equation}
\begin{aligned}
v_{3rrr} = \rd_{rr} W_3 -& \frac{2m}{r} v_3 \rd_r W_3 - \frac{2m}{r}W_3^2 
+\frac{4m}{r^2} (2mv_3^2 + v_3 -m) W_3 \\
& + \frac{m}{r^3} (1-v_3^2) ( 6mv_3^2 +6mv_3 +2 -2m^2 ).
\end{aligned}
\end{equation}
Hence 
\begin{equation}
\begin{aligned}
\nor{\frac{v_{3rrr}}{r} q}{\La} 
\lesssim 
&\nor{(q\eh )_{rr} \frac{q}{r}}{\La} + \nor{(q\eh )_{r} \frac{q}{r^2}}{\La} \\
&\hspace{-5pt} +\nor{\frac{q^3}{r^2}}{\La}
+\nor{\frac{q^2}{r^3}}{\La} + \nor{\frac{1+v_3}{r^2} \frac{q}{r^2}}{\La}.
\end{aligned}
\end{equation}
Each term is estimated as follows:
\begin{itemize}
\item The fifth term is treated in the same manner, namely,
\begin{equation}
\begin{aligned}
&\nor{\frac{1+v_3}{r^2} \frac{q}{r^2}}{\La} \lesssim 
\nor{\frac{1+v_3}{r^2}}{L^{\sy}_tL^4_x} \nor{\frac{q}{r^2}}{\Lc} \\
& \lesssim 
\left( \us^{-3/2} \sig^{3/4} +\nor{q}{L^{\sh}_t L^8_x}^2 +\nor{q}{\Ld}^3  \right) \nor{\frac{q}{r^2}}{L_t^\infty L^2_x}.
\end{aligned}
\end{equation}
\item The third and forth terms are bounded by
\begin{equation}
%\begin{aligned}
%&\nor{\frac{q^3}{r^2}}{\La} +\nor{\frac{q^2}{r^3}}{\La} \\
\left( \nor{q}{\Lb}^2 + \sig^{\Nb}\nor{\frac{q}{r}}{\La} \right) \nor{\frac{q}{r^2}}{L_t^\infty L^2_x}.
%\end{aligned}
\end{equation}
\item For the second term,
\begin{equation}
\begin{aligned}
&\nor{(q\eh )_r\frac{q}{r^2}}{\La} \le \nor{q_r \frac{q}{r^2}}{\La} + \nor{q v_r \frac{q}{r^2}}{\La}  \\
&\lesssim \nor{q_r \frac{q}{r^2}}{\La} + \nor{\frac{q^3}{r^2}}{\La} +\nor{\frac{q^2}{r^3}}{\La} \\
&\lesssim \left( \nor{q}{\Lb}^2 + \sig^{\Nb}
\left( \nor{q_r}{\La} + \nor{\frac{q}{r}}{\La}\right) \right) \nor{\frac{q}{r^2}}{L_t^\infty L^2_x}.
\end{aligned}
\end{equation}
\item We note that
\begin{equation}\label{a170}
\begin{aligned}
(q\eh)_{rr} 
=q_{rr}\eh - 2 \left( q_{1r} \eh_r + q_{2r}(J\eh)_r \right) 
&-\left( q_1 (v_r\cdot \eh ) + q_2 (v_r\cdot J\eh ) \right) v_r \\
&-\left( q_1 (v_{rr}\cdot \eh )+ q_2 (v_{rr}\cdot J\eh ) \right)v .
\end{aligned}
\end{equation}
Hence, the first term is bounded by
\begin{equation}
\begin{aligned}
%&\nor{(q\eh )_{rr}\frac{q}{r}}{\La} \\
& \left(
\us^{-\Na} \sig^{\Nc}  
+ \nor{q}{\Lb}^2 + \nor{q}{\Ld}^3 %\right. \\
%& \hspace{135pt} \left. 
+\sig^{\Nb} (\nor{q_r}{\Ld} + \nor{\frac{q}{r}}{\Ld})   
\right) \\
& \hspace{80pt} \times \left( \nor{q_{rr}}{L^\infty_t L^2_x} + \nor{\frac{q_r}{r}}{L^\infty_t L^2_x} + \nor{\frac{q}{r^2}}{L^\infty_t L^2_x} \right) .
\end{aligned}
\end{equation}
\end{itemize}
Therefore, we have
\begin{equation}
\begin{aligned}
&\nor{A_6}{\La} \\
&\lesssim 
\left(
\us^{-\Na} \sig^{\Nc}  
+ \nor{q}{\Lb}^2 + \nor{q}{\Ld}^3 +\sig^{\Nb} (\nor{q_r}{\Ld} + \nor{\frac{q}{r}}{\Ld})   
\right) \\
& \hspace{95pt} \times \left( \nor{q_{rr}}{L^\infty_t L^2_x} + \nor{\frac{q_r}{r}}{L^\infty_t L^2_x} + \nor{\frac{q}{r^2}}{L^\infty_t L^2_x} \right) .
\end{aligned}
\end{equation}
Applying the equivalence (\ref{a100}) to $\tq$, 
we obtain (\ref{c31}).

\subsection{Proof of Proposition \ref{P9}}

In this section, we prove Proposition \ref{P9}. 
Let $u\in \Sig_m\cap \dot{H}^3$ with $\del = \sqrt{\cae(u) -4\pi m}< \del_0$. 
We first note that 
%\begin{equation}
%\begin{aligned}
%\nor{u}{\dot{H}^3}
%= \nor{\Del u}{\dot{H}^1}
%\lesssim \nor{\rd_r \Del u}{L^2} + \nor{\frac{1}{r}\rd_\theta \Del u}{L^2}
%\lesssim 
%\nor{\rd_r H_m v_j}{}
%\end{aligned}
%\end{equation}
since we have the equivalence (\ref{a140}), 
it suffices to find the bounds for 
%\begin{equation}
$
\rd_r H_m v_j$, $\frac{1}{r} H_m v_j$ and $\rd_r H_0 v_3
$ 
%\end{equation}
for $j=1,2$. %Without loss of generality, we may consider only the case $j=1$. 
%We first consider $\rd_r H_m v_1$ and $\frac{1}{r} H_m v_1$, while the case when $j=2$ is similar. 
By straightforward computations, we have
\begin{equation}
\begin{aligned}
\rd_r H_m v_1 = 
&\rd_{rr} W_1 + (1-mv_3) \frac{1}{r}\rd_r W_1 
- 2m \frac{1}{r}W_1 W_3 + 2m^2 v_3v_1 \frac{1}{r^2} W_3 \\
& +(-3m^2(1-v_3^2) +mv_3 -1) \frac{1}{r^2}W_1
-mv_1 \frac{1}{r}\rd_r W_3 \\
&+ (m+3m^2 v_3) \frac{1}{r^2} W_3
+\frac{4m^2}{r^3}v_1(1-v_3^2) +\frac{6m^3}{r^3} v_3 v_1 (1-v_3^2),
\end{aligned}
\end{equation}
\begin{equation}
\frac{1}{r} H_m v_1 = \frac{1}{r} \rd_r W_1 + (1-mv_3)\frac{1}{r^2}W_1 
-mv_1\frac{1}{r^2}W_3 - \frac{2m^2}{r^3} v_1 (1-v_3^2),
\end{equation}
\begin{equation}
\begin{aligned}
\rd_r H_0 v_3 
=
& \rd_{rr} W_3 + (-2mv_3 +1)\frac{1}{r}\rd_r W_3 
-2m \frac{1}{r} W_3^2 \\
&+(8m^2v_3 +2mv_3 -4m^2 -1) \frac{1}{r^2}W_3 \\
&+\frac{1-v_3^2}{r^3} (6m^3v_3^2 +4m^2v_3 -2m^3).
\end{aligned}
\end{equation}
Hence, we obtain
\begin{equation}
\begin{aligned}
\nor{u}{\dot{H}^3}
&\lesssim 
\nor{\rd_{rr} W}{L^2} + \nor{\frac{1}{r} \rd_r W}{L^2} + \nor{\frac{1}{r^2}W}{L^2} 
+ \nor{\frac{1}{r}W^2}{L^2} \\
 & \quad + \sum_{j=1}^2 \nor{\frac{v_1 (1-v_3^2)}{r^3}}{L^2} 
+ \nor{\frac{1-v_3^2}{r^3} (6m^3v_3^2 +4m^2v_3 -2m^3)}{L^2}\\
&=: \sum_{k=1}^6 B_k.
\end{aligned}
\end{equation}
By (\ref{a170}), we have
\begin{equation}\label{a178}
B_1 \lesssim \nor{q_{rr}}{L^2} + \nor{q_r v_r}{L^2} +\nor{qv_{rr}}{L^2} +\nor{q \left| v_r\right|^2}{L^2}.
\end{equation}
Direct calculations yield 
\begin{equation}
v_{1rr}  = \rd_r W_1 - mv_3 \frac{1}{r} W_1 -mv_1 \frac{1}{r} W_3 +\frac{1}{r^2}v_1 (-m^2+2m^2v_3^2 + mv_3),
\end{equation}
\begin{equation}
v_{3rr}  = \rd_r W_3 - 2mv_3 \frac{1}{r} W_3 -\frac{1-v_3^2}{r^2} (m+2m^2v_3),
\end{equation}
hence 
\begin{equation}\label{a181}
\begin{aligned}
\nor{qv_{rr}}{L^2} &\lesssim \nor{q \rd_r W}{L^2} +\nor{q \frac{1}{r}W}{L^2} + \nor{\frac{q}{r^2}}{L^2}\\
&\lesssim \nor{q q_r}{L^2} + \nor{q^2v_r}{L^2} +\nor{\frac{q^2}{r}}{L^2} +\nor{\frac{q}{r^2}}{L^2}.
\end{aligned}
\end{equation}
Each term is estimated by
\begin{equation}
\nor{q q_r}{L^2} \le \nor{q}{L^2} \nor{q_r}{L^\infty} \lesssim \nor{q_{rr}}{L^2} + \nor{\frac{q_r}{r}}{L^2},
\end{equation}
\begin{equation}
\nor{\frac{q^2}{r}}{L^2} \le \nor{q}{L^2} \nor{\frac{q}{r}}{L^\infty} \lesssim \nor{\frac{q_r}{r}}{L^2} +\nor{\frac{q}{r^2}}{L^2},
\end{equation}
\begin{equation}
%\begin{aligned}
\nor{q^2 v_r}{L^2} \lesssim \nor{q^2 \left| v_r-\frac{m}{r} JRv \right|}{L^2} + \nor{\frac{q^2}{r}}{L^2}%\\ % \lesssim \nor{q_{rr}}{L^2} + \nor{\frac{q_r}{r}}{L^2}
%&\lesssim \nor{q^3}{L^2} + \nor{q}{L^2} \nor{\frac{q}{r}}{L^\infty}
\lesssim \nor{q}{L^6}^3 + \nor{\frac{q_r}{r}}{L^2} +\nor{\frac{q}{r^2}}{L^2}.
%\end{aligned}
\end{equation}
Therefore, (\ref{a181}) implies
\begin{equation}
\nor{qv_{rr}}{L^2} 
\lesssim \nor{q_{rr}}{L^2} + \nor{\frac{q_r}{r}}{L^2} + \nor{\frac{q}{r^2}}{L^2} + \nor{q}{L^6}^3.
\end{equation}
Moreover, we obtain
\begin{equation}
\nor{q_r v_r}{L^2} \lesssim \nor{q_r \left| v_r-\frac{m}{r} JRv \right|}{L^2} + \nor{\frac{q_r}{r}}{L^2}
\lesssim \nor{q_{rr}}{L^2} + \nor{\frac{q_r}{r}}{L^2},
\end{equation}
\begin{equation}
\nor{q \left| v_r\right|^2}{L^2} \lesssim \nor{q \left| v_r-\frac{m}{r} JRv \right|^2}{L^2} + \nor{\frac{q}{r^2}}{L^2} \lesssim \nor{q}{L^6}^3 + \nor{\frac{q}{r^2}}{L^2}.
\end{equation}
Hence, (\ref{a178}) provides 
\begin{equation}
B_1 \lesssim \nor{q_{rr}}{L^2} + \nor{\frac{q_r}{r}}{L^2} + \nor{\frac{q}{r^2}}{L^2} + \nor{q}{L^6}^3.
\end{equation}
We can similarly obtain
\begin{equation}
%\begin{aligned}
B_2 \lesssim \nor{\frac{q_r}{r}}{L^2} + \nor{\frac{q}{r} v_r}{L^2} 
%\lesssim \nor{\frac{q_r}{r}}{L^2} + \nor{\frac{q^2}{r}}{L^2} + \nor{\frac{q}{r^2}}{L^2}
\lesssim \nor{\frac{q_r}{r}}{L^2} + \nor{\frac{q}{r^2}}{L^2},
%\end{aligned}
\end{equation}
\begin{equation}
B_3 \le \nor{\frac{q}{r^2}}{L^2},\qquad B_4 \le \nor{\frac{q_r}{r}}{L^2} + \nor{\frac{q}{r^2}}{L^2}.
\end{equation}
%\begin{equation}
%B_4 \le \nor{\frac{q_r}{r}}{L^2} + \nor{\frac{q}{r^2}}{L^2}.
%\end{equation}
In order to estimate $B_5$, we first observe that 
\begin{equation}
%\begin{aligned}
v_1^2 +v_2^2 = \left[ \left( (1+\gam ) h_1 -z_2h_3 \right)^2 + z_1^2 \right] \left( \frac{r}{s} \right) 
\lesssim \left[ |z|^2 + h_1^2 \right] \left( \frac{r}{s} \right)
%\end{aligned}
\end{equation}
where we use (\ref{a146}). Hence 
\begin{equation}
B_5 \lesssim \nor{\frac{\left( v_1^2 + v_2^2\right)^{\Na}}{r^3}}{L^2}
\lesssim s^{-2} \nor{\frac{|z|^3 +h_1^3}{r^3}}{L^2} 
\lesssim \nor{q}{L^6}^3 + s^{-2}
\end{equation}
where we use (\ref{a149}) in the last inequality. \par
It remains to control $B_6$. 
To this end, we require dividing the case into $m=1$ and $m=2$. 
When $m=1$, we have the factorization 
\begin{equation}\label{a194}
\frac{1-v_3^2}{r^3} (6v_3^2 +4v_3 -2) = \frac{2}{r^3} (1-v_3^2) (1+v_3) (3v_3-1).
\end{equation}
Moreover, we have
\begin{equation}
1+v_3 \lesssim \left[ (1+h_3) + |z| \right] \left( \frac{r}{s} \right) .
\end{equation}
Hence 
\begin{equation}
B_6 \lesssim 
s^{-2} \nor{\frac{|z|^3 + (1+h_3)^3 + h_1^3}{r^3}}{L^2} \lesssim \nor{q}{L^6}^3 + s^{-2}.
\end{equation}
In the case when $m\ge 2$, we no longer have such a factorization as (\ref{a194}). 
Instead, we are able to use Lemma \ref{L9} in this case. 
First, we have
\begin{equation}
B_6 \lesssim \nor{v_1^2 +v_2^2}{L^2} \lesssim 
s^{-2} \nor{\frac{|z|^2 + h_1^2}{r^3}}{L^2} \lesssim s^{-2} \left( \nor{\frac{z^2}{r^3}}{L^2} +1 \right) .
\end{equation}
Here, we apply Lemma \ref{L9} with $a= \frac{3}{2}$ and $p=4$. Then we obtain
\begin{equation}
\begin{aligned}
\nor{\frac{z^2}{r^3}}{L^2} = \nor{\frac{z_r}{r^{\Nb}}}{L^4}^2 
&\lesssim \nor{\frac{1}{r} \left( z_r -\frac{mz}{r} \right)^2}{L^2} \\
& \lesssim \nor{\frac{1}{r} \left( z_r +\frac{m}{r} h_3z \right)^2}{L^2} 
+ \nor{\frac{1}{r} \left( \frac{m}{r} (1 +h_3) z \right)^2}{L^2} \\
& =: E_1 + E_2.
\end{aligned}
\end{equation}
%Since $L_0 z = e^{m\theta R} [ (z_r +\frac{m}{r} z) \bj ]$, 
If we use the relation (\ref{a147}), then we have
\begin{equation}
\begin{aligned}
E_1 &= \nor{\frac{1}{r} \left| 
s e^{-\al R} (q\eh )(sr) - (\gam h)_r -\frac{2m}{r} h_3 \gam h - \frac{m}{r} \xi_3 \xi
 \right|^2}{L^2}\\
&\lesssim \nor{\frac{1}{r} s^2 \left[ \left| q\eh \right|^2 \right] (sr)}{L^2} 
+ \nor{ \frac{1}{r} \left| (\gam h)_r \right|^2 }{L^2} \\
& \hspace{60pt} + \nor{ \frac{1}{r} \left| \frac{2m}{r} h_3 \gam h \right|^2 }{L^2}
+ \nor{ \frac{1}{r} \left| \frac{m}{r} \xi_3 \xi \right|^2 }{L^2}\\
&\lesssim s^{2} \left( \nor{\frac{q_r}{r}}{L^2} +\nor{\frac{q}{r^2}}{L^2} +\nor{q}{L^6}^3 \right).
\end{aligned}
\end{equation}
%Noting that $|\xi| \lesssim |z|$, we obtain
%\begin{equation}
%E_1 \lesssim s^{2} \left( \nor{\frac{q_r}{r}}{L^2} +\nor{\frac{q}{r^2}}{L^2} +\nor{q}{L^6}^3 \right) 
%+1.
%\end{equation}
$E_2$ is estimated by
\begin{equation}
E_2 \lesssim \nor{\frac{(1+h_3)^3}{r^3} z^2}{L^2} \lesssim 1.
\end{equation}
Thus
\begin{equation}
B_6 \lesssim \nor{\frac{q_r}{r}}{L^2} +\nor{\frac{q}{r^2}}{L^2} +\nor{q}{L^6}^3 + s^{-2}.
\end{equation}
Hence we complete the proof.

\section{Proof of Proposition \ref{P4}}

\textit{Step 1}. 
We begin with showing the following claim:\bigskip\par
\vspace{-5pt}
There exists $\del_1 >0$ and $C_1>0$ such that for $u\in \Sig_m$ 
with $\nor{u-Q}{\dot{H}^1}$\\ $<\del_1$, there is a pair $(s,\al) = (s_0(u), \al_0(u))\in \R^+\times \bbt^1$ 
such that 
\begin{itemize}
\item $\inp{h_1}{z}_{\dot{H}^1} =0$,
\item $|s_0(u)-1| +|\al_0 (u)| \le C_1 \nor{u-Q}{\dot{H}^1}$.
\end{itemize}
Moreover, if $(s,\al)$ satisfies $\inp{h_1}{z}_{\dot{H}^1}=0$ and $|s-1| + |\al| \le C_1\del_1$, 
then $(s,\al)$ coincides with $(s_0 (u),\al_0 (u))$. \bigskip\par
\vspace{-3pt}
First, we introduce the function space
\begin{equation}
Y:= \left\{ e^{m\theta R} v(r) \left|  v(r):(0,\infty)\to \R^3, v_1,v_2\in \dot{H}^1_e, v_3\in L^\infty , \rd_r v_3\in L^2_e  \right. \right\}
\end{equation}
with the norm
\begin{equation}
\nor{u}{Y} := \nor{v_1}{\dot{H}^1_e} + \nor{v_2}{\dot{H}^1_e} +\nor{v_3}{L^\infty} +\nor{\rd_r v_3}{L^2_e}.
\end{equation}
Then, we can easily check that $(Y,\nor{\cdot}{Y})$ is a Banach space. 
Define $F={}^t (F_1,F_2): Y\times \R^+ \times\R \to \R^2$ by
\begin{equation}
\begin{pmatrix}
F_1(u,s,\al) \\
F_2(u,s,\al)
\end{pmatrix}
:= 
\begin{pmatrix}
\inp{h_1}{z_1}_{\dot{H}^1_e} \\
\inp{h_2}{z_2}_{\dot{H}^1_e}
\end{pmatrix}
=
\begin{pmatrix}
-\int_0^\infty e^{-\al R} (v(sr)\cdot \bj) (H_m h_1) rdr\\
-\int_0^\infty e^{-\al R} (v(sr)\cdot J^h\bj) (H_m h_1) rdr
\end{pmatrix}
\end{equation}
where $H_m$ is defined in (\ref{a99}). 
Then, $F$ is $C^1$, and we have
\begin{equation}
F(Q,1,0) =0,\quad 
\begin{pmatrix}
\rd_s F_1 (Q,1,0) & \rd_\al F_1 (Q,1,0)\\
\rd_s F_2 (Q,1,0) & \rd_\al F_2 (Q,1,0)
\end{pmatrix}
=
\nor{h_1}{\dot{H}^1_e}^2
\begin{pmatrix}
0 & -1\\
m & 0
\end{pmatrix}
.
\end{equation}
Thus, by the implicit function theorem, there exist neighborhoods $V\subset Y$ of $Q$, and $W\subset \R^+ \times\R$ of $(1,0)$, respectively, and 
%$V\times W \subset Y\times (\R_{>0} \times \R)$, and 
a function $(s_0(u), \al_0(u)) :V\to W$ such that 
\begin{itemize}
\item $(s_0(u),\al_0(u))$ is $C^1$ on $V$.
\item $(s_0 (Q),\al_0 (Q)) = (1,0)$.
\item For $(u,s,\al) \in V\times W$, 
\begin{equation}\label{a212}
F(u,s,\al) = 0 \Longleftrightarrow (s,\al) = (s(u),\al(u)).
\end{equation}
\item For any $u\in V$. we have
\begin{equation}
\det 
\begin{pmatrix}
\rd_s F_1 (u,s_0(u),\al_0(u) ) & \rd_\al F_1 (u,s_0(u), \al_0(u))\\
\rd_s F_2 (u,s_0(u),\al_0(u)) & \rd_\al F_2 (u,s_0(u),\al_0(u))
\end{pmatrix}
\neq 0
\end{equation}
and 
\begin{equation}\label{a214}
\begin{pmatrix}
\inp{d_u s_0(u)}{\del u}_{Y^*,Y} \\
\inp{d_u \al_0(u)}{\del u}_{Y^*,Y}
\end{pmatrix}
= -
\begin{pmatrix}
\rd_s F_1 & \rd_\al F_1 \\
\rd_s F_2 & \rd_\al F_2 
\end{pmatrix}
^{-1}
\hspace{-3pt}
\begin{pmatrix}
\inp{d_u F_1}{\del u}_{Y^*,Y} \\
\inp{d_u F_2}{\del u}_{Y^*,Y}
\end{pmatrix}
\end{equation}
for $\del u \in Y$.
\end{itemize}
%We may assume that $V$ is convex, say, a ball centered at $Q$. 
%Here, we take a constant $C>0$ such that 
%\begin{equation}
%\left| 
%\begin{pmatrix}
%\rd_s F_1 & \rd_\al F_1 \\
%\rd_s F_2 & \rd_\al F_2 
%\end{pmatrix}
%^{-1}
%\right| \le \frac{1}{2}C^{\Nb}, \quad
%\nor{d_u F_j (u,s_0(u), \al_0(u))}{Y^*} \le \frac{1}{2}C^{\Nb}
%\end{equation}
%for all $u\in V$ and $j=1,2$.  
From (\ref{a214}), there is a constant $C$ such that 
\begin{equation}
\nor{d_u s(u)}{Y^*} \le C,\quad \nor{d_u \al (u)}{Y^*} \le C
\end{equation}
for all $u\in V$. (If necessary, we replace $V$ by a smaller neighborhood.)\par
Now, let $u\in \Sig_m$. 
Then $u\in Y$, and from Lemma \ref{L1}, we have 
\begin{equation}
\nor{u-Q}{Y} \le C(\nor{u}{\dot{H}^1}) \nor{u-Q}{\dot{H}^1}.
\end{equation}
In particular, there exists $\del_1 >0$ such that $u\in V$ if $\nor{u-Q}{Y} <\del_1$.
Then, $Q+t(u-Q)\in V$ for all $t\in [0,1]$, and
\begin{equation}
\begin{aligned}
|s(u)-1| =|s(u)-s(Q)| &= \left| \int_0^1 \frac{d}{dt} s(Q+t(u-Q)) dt \right| \\
&= \left| \int_0^1\inp{d_u s(Q+t(u-Q))}{u-Q}_{Y^*,Y} dt \right| \\
& \le C\nor{u-Q}{Y} \le \frac{1}{2}C_1 \nor{u-Q}{\dot{H}^1} 
\end{aligned}
\end{equation}
for some $C_1 >0$, and similarly 
\begin{equation}
|\al(u)| =|\al (u)-\al(Q)| \le \frac{1}{2}C_1 \nor{u-Q}{\dot{H}^1}. 
\end{equation}
Moreover, when $\del_1$ is sufficiently small, we have 
$(s,\al)\in W$ if $|s-1| +|\al| \le C_1\del_1$. 
Hence, by (\ref{a212}), the claim of Step 1 follows. \bigskip\\
%%%%%%%%%%%%%%%%%%%%%%%%%%%%%%%%%%%%%%%%%%%%%%%%%%%%%%%%%%%%%%%%%%%%%%%%%%%%%%%%%%%%%%%%%%%%%%%
%%%%%%%%%%%%%%%%%%%%%%%%%%%%%%%%%%%%%%%%%%%%%%%%%%%%%%%%%%%%%%%%%%%%%%%%%%%%%%%%%%%%%%%%%%%%%%%
\textit{Step 2}. We now prove the existence of $s(u)$, $\al(u)$ as in Proposition \ref{P4} 
(i) and (ii). 
Suppose $\del_0 >0$ and $u\in\Sig_m$ with $\del <\del_0$. 
%By Proposition \ref{P1}, 
%there exist $\del_2>0$ and $C_2>0$ such that we can take the unique $(s_*(u),\al_*(u))$ satisfying (\ref{1}) 
To avoid the ambiguity, we change the notation $\del_0$ in Proposition \ref{P1} into $\del_2$. 
Then, there exists $C_2>0$ such that if $\del <\del_2$, then %there exists $C_2$ such that 
\begin{equation}
\nor{u- e^{\al_*(u) R} Q(\frac{\cdot}{s_*(u)})}{\dot{H}^1} \le C_2 \del.
\end{equation}
Here we set $\tilde{u}(x) := e^{-\al_*(u) R} u(s_*(u) x)$. 
By the scale invariance, we have $\nor{\tu -Q}{\dot{H}^1} \le C_2\del$. 
Thus, if we choose $\del_0 = \min \{ C_2^{-1} \del_1 ,\del_2\}$, 
we can apply the result in Step 1 to $\tu$. Namely, 
there exists $(s_0(\tu), \al_0(\tu))$ such that $F(u,s_0(\tu) ,\al_0(\tu))=0$ and that 
$\left| s_0(\tu) -1\right| + \left| \al_0(\tu) \right| \le C_1 \nor{\tu -Q}{\dot{H}^1} \le C_1C_2\del$. 
Here we set $C_3:=C_1C_2$ and 
\begin{equation}
s(u) := s_* (u) s_0(\tu) ,\quad \al(u) := \al_* (u) + \al_0 (\tu).
\end{equation}
Then we have
\begin{equation}
F(u ,s(u),\al(u)) = F(\tu, s_0(\tu), \al_0(\tu)) =0,
\end{equation}
\begin{equation}\label{a223}
\left| \frac{s(u)}{s_*(u)} -1\right| + \left| \al(u) - \al_*(u) \right| \le C_3 \del.
\end{equation}
Moreover, if $(s,\al)$ satisfies (\ref{a223}) in which $(s(u),\al(u))$ is replaced by $(s,u)$, then 
we obtain $(\frac{s}{s_*}, \al-\al_*)=(s_0(\tu),\al_0(\tu)) $ since 
$F(\tu, \frac{s}{s_*}, \al-\al_*) =0$. 
Therefore, we achieve (i) and (ii) in Proposition 4.\bigskip\par 
%%%%%%%%%%%%%%%%%%%%%%%%%%%%%%%%%%%%%%%%%%%%%%%%%%%%%%%%%%%%%%%%%%%%%%%%%%%%%%%%%%%%%%%%%%%%%%%
%%%%%%%%%%%%%%%%%%%%%%%%%%%%%%%%%%%%%%%%%%%%%%%%%%%%%%%%%%%%%%%%%%%%%%%%%%%%%%%%%%%%%%%%%%%%%%%
\textit{Step 3}. 
Finally, we show the regularity property of $(s(u),\al(u))$ as in Proposition \ref{P4} (iii). 
We only consider the case $u(t)\in C(I;\Sig_m)\cap C^1 (I;L^2(\R^2))$. 
(The case $u(t)\in C(I;\Sig_m)\cap W^{1,\infty}(I;L^2(\R^2))$ can be derived by small modifications.) \par
Fix $t_0\in I$. It suffices to show that $s(u(t))$, $\al(u(t))$ are $C^1$ 
in some neighborhood of $t=t_0$. 
Define linear transform $\mct$ on $Y$ by $\mct u:= e^{-\al_* (u(t_0))R} u(s_*(u(t_0))\cdot )$ for $u\in Y$. 
Then, if $t-t_0$ is sufficiently small, we have
\begin{equation}
s(u(t)) = s_0 (\mct u(t)) s_*(u(t_0)) ,\quad \al (u(t)) = \al (\mct u(t)) + \al_* (u(t_0))
\end{equation}
by the uniqueness proved in Step 1. 
Hence, we may assume $s_*(u(t_0)) =1$, $\al_*(u(t_0))=0$, and thus $\nor{u(t_0)-Q}{\dot{H}^1} < \del_1$. Under this simplification, 
we can use the fact that $s(u)$ and $\al(u)$ 
is $C^1$-function on $Y$ by Step 1. 
By continuity, there exists $\sig_0 >0$ such that $\nor{u(t)- Q}{\dot{H}^1} <\del_1$ if $|t-t_0| \le\sig_0 $. 
Note that for fixed $s$ and $\al$, $F(\cdot, s,\al)$ is bounded linear functional on $L^2_e$.\par 
It suffices to show that $s(u(t))$ and $\al(u(t))$ is $C^1$ on $(t_0-\sig_0 , t_0+\sig_0)$ and
\begin{equation}\label{a225}
\begin{pmatrix}
\frac{d}{dt} s(u(t))\\
\frac{d}{dt} \al (u(t))
\end{pmatrix}
 = 
-
B(u(t)) 
\begin{pmatrix}
F_1(\rd_tu(t), s(u(t)), \al (u(t)) )  \\
F_2(\rd_tu(t), s(u(t)), \al (u(t)) )
\end{pmatrix}
,
\end{equation}
where 
\begin{equation}
\begin{aligned}
B(u(t)) &\equiv 
\begin{pmatrix}
B_{11}(u) & B_{12}(u)\\
B_{21}(u) & B_{22}(u)
\end{pmatrix}
\\
&:=
\begin{pmatrix}
\rd_s F_1 (u, s(u), \al (u)) & \rd_\al F_1 (u(t), s(u), \al (u))\\
\rd_s F_2 (u, s(u), \al (u)) & \rd_\al F_2 (u(t), s(u), \al (u))
\end{pmatrix}
^{-1}
\end{aligned}
\end{equation}
for $t\in (t_0-\sig_0 , t_0+\sig_0 )$. Let $\sig \in \R$ with $|t+\sig -t_0|<\sig_0$. 
Then, 
\begin{equation}
\begin{aligned}
&\frac{1}{\sig} \left[ s(u(t +\sig)) -s(u(t)) \right] + \sum_{j=1}^2 B_{1j}F_j(\rd_t u,s(u),\al(u) ) \\
&=\frac{1}{\sig} \int_{0}^{1} \inp{d_us(u(t) + \xi (u(t+\sig) - u(t) ))}{\ u(t+\sig) -u(t)}_{Y^* ,Y} d\xi\\
&\hspace{140pt} +\sum_{j=1}^2 B_{1j}F_j(\rd_t u,s(u),\al(u) )\\
&=: K_1 + K_2,
\end{aligned}
\end{equation}
where
\begin{equation}
K_1:= \int_{0}^{1} \inp{d_us(c(\xi))- d_u s(c(0))}{\ \Del_\sig u}_{Y^* ,Y} d\xi
\end{equation}
\begin{equation}
K_2:= \inp{d_us(u)}{\Del_\sig u}_{Y^* ,Y} + \sum_{j=1}^2 B_{1j}F_j(\rd_t u,s(u),\al(u) )
\end{equation}
\begin{equation}
c(\xi ):= u(t) +\xi \left( u(t+\sig) -u(t) \right) ,\quad 
\Del_\sig u = \frac{1}{\sig} \left( u(t+\sig) -u(t) \right).
\end{equation}
Using (\ref{a214}), we have
\begin{equation}
%\begin{aligned}
K_2= 
- \sum_{j=1}^2 B_{1j} \left[ F_j \left(\Del_\sig u - \rd_t u(t) ,\ s(u(t)) ,\ \al(u(t)) \right) \right] 
\to 0 %\quad \text{ as } \sig\to 0,
%\end{aligned}
\end{equation}
as $\sig\to 0$, since $u(t)\in C^1 (I;L^2)$. On the other hand, we can write 
\begin{equation}
\begin{aligned}
&\inp{d_us(c(\xi))- d_u s(c(0))}{\ \Del_\sig u}_{Y^* ,Y}\\
& = -\sum_{j=1}^2 \left[ B_{1j} (c(\xi ) ) F_j (\Del_\sig u, s(c(\xi )), \al (c(\xi ) ) ) \right. \\
&\hspace{90pt} \left. -B_{1j} (c(0)) F_j (\Del_\sig u, s(c(0 )), \al (c(0) ) )  \right]\\
&= - (L_1 +L_2 +L_3), 
\end{aligned}
\end{equation}
where
\begin{equation}
L_1 := \sum_{j=1}^2 \left[ B_{1j} (c(\xi)) - B_{1j} (c(0)) \right] F_j 
\left( \Del_\sig u , s(c(\xi)), \al(c(\xi)) \right),
\end{equation}
\begin{equation}
\begin{aligned}
L_2 &:= \sum_{j=1}^2 B_{1j} (c(0)) \left[ F_j 
\left( \Del_\sig u , s(c(\xi)), \al(c(\xi)) \right) \right. \\
& \hspace{115pt}\left. -
F_j 
\left( \Del_\sig u , s(c(0)), \al(c(\xi)) \right) 
\right],
\end{aligned}
\end{equation}
\begin{equation}
\begin{aligned}
L_3 &:= \sum_{j=1}^2 B_{1j} (c(0)) \left[ F_j 
\left( \Del_\sig u , s(c(0)), \al(c(\xi)) \right) \right. \\
& \hspace{115pt}\left. -
F_j 
\left( \Del_\sig u , s(c(0)), \al(c(0)) \right) 
\right].
\end{aligned}
\end{equation}
Each $|L_k|$ is bounded by some constants independent of $\sig$ and $\xi$, 
and converges to $0$ as $\sig\to 0$. 
Hence, by Lebesgue's dominant convergence theorem, we obtain
$K_1 \to 0$ as $\sig \to 0$. 
The same argument can be applied to $\al$, and hence 
we achieve (\ref{a225}). \par
The continuity of (\ref{a225}) can be shown in the same manner as above. 
Hence, $s(u(t))$, $\al(u(t))$ is $C^1$.

%%%%%%%%%%%%%%%%%%%%%%%%%%%%%%%
%%%%%%%   Section 7   %%%%%%%%%
%%%%%%%%%%%%%%%%%%%%%%%%%%%%%%%

\section{Some Technical Lemmas}
\subsection{Continuity of Reconstruction in Higher Regularity}
In this section, we show some technical lemmas used in the previous sections. 
We first make a further observation concerned with the 
continuous dependence of the map $(q ,s,\al)\mapsto u$ in Proposition \ref{P5}. 

\begin{lem}\label{L10}
The reconstruction map $(q,s,\al)\mapsto u$ in Proposition \ref{P5} %$H^1_e \times \R^+ \times \bbt^1 \ni (q,s,\al) \mapsto u \in \Sig_m$ is 
is continuous from $H^1_e \times \R^+ \times \bbt^1$ to $\Sig_m\cap \dot{H}^2$.
\end{lem}

\proof
We first fix $(q,s,\al) \in H^1_e\times \R^+\times \bbt^1$. 
It suffices to show that
$\nor{ u(q',s',\al') - u(q,s,\al) }{\dot{H}^2}\to 0 $ as $(q',s',\al')\to (q,s,\al)$ in 
$H^1_e\times \R^+\times \bbt^1$. 
(Here $u(q,s,\al)\in \Sig_m$ denotes the map reconstructed from $(q,s,\al)$.) 
In the proof, we write the difference $q'-q$ as $\del q$, and also adopt this convention 
to other quantities. 
\par
We write $u=e^{m\theta R }v(r) := u(q,s,\al)$, $u' = e^{m\theta R} v'(r) :=u(q',s',\al')$. 
Taking account of the equivalence 
\begin{equation}\label{d71}
\nor{\del u}{\dot{H}^2} \sim \sum_{j=1}^2 \nor{\del H_m v_j}{L^2_e} + \nor{\del H_0 v_3}{L^2_e},
\end{equation}
we only need to see the right hand side of (\ref{d71}). %and  
%$\nor{\del H_m v_1}{L^2_e}$
Direct computations yield
\begin{equation}\label{b240}
\del H_m v_j = \del \rd_r W_j  
 - \del \left( \frac{m}{r} v_1 W_3 \right) + 
\del \left( \frac{1-mv_3}{r} W_j \right) 
- \del \left( \hspace{-1pt} \frac{2m^2}{r^2} v_j (1-v_3^2) \hspace{-1pt} \right),
\end{equation}
\begin{equation}\label{b241}
\del H_0 v_3 = \del \rd_r W_3 + \del \left( \frac{1-2mv_3}{r} W_3 \right) 
- \del \left( \frac{2m^2}{r^2} v_3 (1-v_3^2) \right)
\end{equation}
for $j=1,2$. Here, we have
\begin{equation}
\begin{aligned}
&\nor{\del \rd_r W}{L^2_e} =\nor{\del (q_r \eh + q \eh_r )}{L^2_e}  \\
&\le \nor{\del q_r}{L^2_e} +\nor{q_r (\del\eh )}{L^2_e} + \nor{(\del q) \eh'_r }{L^2} + \nor{q(\del\eh_r)}{L^2} \\
&\le \nor{\del q_r}{L^2_e} + \nor{q_r}{L^2_e} \nor{\del \eh}{L^\infty} 
+ \nor{\del q}{L^\infty} \nor{\eh'_r}{L^2_e} +\nor{q}{L^\infty} \nor{\del \eh_r}{L^2_e}\\
&\lesssim \left( 1+ \nor{u'}{\dot{H}^1} \right) \nor{\del q}{H^1_e} +\nor{q}{H^1_e} \nor{\del u}{\dot{H}^1} 
\end{aligned}
\end{equation}
where we use Lemma \ref{L1} and (\ref{a68}). Moreover, for $j=1,2,3$, we have
\begin{equation}
\begin{aligned}
\nor{\del \left( v_j \frac{W}{r} \right) }{L^2_e} &= \nor{(\del v_j)\frac{W'}{r} + v_j \frac{\del W}{r} }{L^2_e}\\
&\le \left( \nor{q}{H^1_e} +\nor{q'}{H^1_e} \right) \nor{\del u }{\dot{H}^1} + \nor{\del q}{H^1_e}.
\end{aligned}
\end{equation}
Therefore, %all terms except last terms in the right hand sides of (\ref{b240}) and (\ref{b241}) converge 
%to $0$ as $(q',s',\al')\to (q,s,\al)$ in $H^1_e\times \R^+\times \bbt^1$. \par
in both (\ref{b240}) and (\ref{b241}), all terms but the last term converge to $0$ as 
$(q',s',\al')\to (q,s,\al)$ in $H^1_e\times \R^+\times \bbt^1$. 
In order to derive the convergence of the last terms, it suffices to show that 
$
\del \left( \frac{1-v_3^2}{r^2} \right) \to 0
$ 
in $L^2_e$. Now, we have
\begin{equation}\label{b244}
\begin{aligned}
\nor{\del\frac{1-v_3^2}{r^2}}{L^2_e} \le  
&\nor{\left[ \frac{1-v_3'^2}{r^2} \right] (r) - 
\frac{s'^2}{s^2} \left[ \frac{1-v_3'^2}{r^2} \right] \left(\frac{s'}{s} r\right)}{L^2_e}\\
&\hspace{30pt} +\nor{\frac{s'^2}{s^2} \left[ \frac{1-v_3'^2}{r^2} \right] \left(\frac{s'}{s} r\right) 
- \frac{1-v_3^2}{r^2} (r) }{L^2_e}.
\end{aligned}
\end{equation}
The first term in (\ref{b244}) tends to $0$ by (\ref{a150}). 
Since 
\begin{equation}
%\begin{aligned}
\left[ 1-v_3^2 \right] (sr) %&= 1- \left( (1+\gam) h_3 + z_2 h_1\right)^2 
= h_1^2 - (2\gam + \gam^2) h_3^2 - 2(1+\gam) z_2 h_1h_3 -z_2^2 h_1^2,
%\end{aligned}
\end{equation}
the second term in (\ref{b244}) can be written as 
\begin{equation}\label{b246}
\begin{aligned}
&\frac{1}{s} \nor{ \del \left( \frac{
 h_1^2 - (2\gam + \gam^2) h_3^2 -2(1+\gam) z_2h_1h_3 -z_2^2 h_1^2
}{r^2} \right) }{L^2_e} \\
%&- \left.
%\frac{1}{s^2} \left[ \frac{
%h_1^2 - (2\gam + \gam^2) h_3^2 -2(1+\gam) z_2h_1h_3 -z_2^2 h_1^2
%}{r^2} \right] \left( \frac{r}{s} \right)
% \right\|_{L^2_e} \to 0
&\lesssim \frac{1}{s} \left( \nor{ \frac{\del\gam}{r^2}}{L^2_e} 
+ \nor{ \frac{\del (\gam ^2)}{r^2}}{L^2_e}
+ \nor{ \frac{ (\del \gam )z_2' + (1+\gam ) \del z_2  }{r} }{L^2_e} 
+ \nor{ \frac{\del (z_2^2)}{r^2}}{L^2_e} \right) .
\end{aligned}
\end{equation}
%as $(q',s',\al') \to (q,s,\al)$. 
Since $|\del\gam| \lesssim (|z| +|z'|) |\del z|$, (\ref{b246}) is bounded by
\begin{equation}\label{b247}
\frac{C}{s} \left( \left( \nor{\frac{z}{r}}{L^4_e} + \nor{\frac{z'}{r}}{L^4_e} \right) \nor{\frac{\del z}{r}}{L^4_e} 
+ \nor{\frac{\del z}{r}}{L^2_e} \right)
\end{equation}
for some constant $C$. Since 
%\begin{equation}
%\nor{\frac{\del z}{r}}{L^2_e}= \nor{\frac{1}{r} e^{-\al R} v()}{}
%\end{equation}
$z= \left( e^{-\al R} v(sr) \right)\cdot \left( \bj + iJ^h\bj \right)$, 
we can easily show that the second term in (\ref{b247}) converges to $0$ 
as $(q',s',\al')\to (q,s,\al)$. 
Hence, it suffices to prove $\nor{\frac{\del z}{r}}{L^4_e} \to 0$. \par
Using the relation (\ref{a147}), we obtain 
\begin{equation}\label{b248}
\begin{aligned}
\nor{\del z_r - \frac{m}{r} \del z}{L^4_e} &\le 
\nor{\del (L_0 z) \bj}{L^4_e} + \nor{m \frac{1+h_3}{r} \del z}{L^4}\\
&\lesssim \nor{\del \left( s e^{-\al R} [q\eh ] (sr) \right)}{L^4} + 
\nor{\del G_0(z)}{L^4_e} + \nor{\del z}{\infty}
\end{aligned}
\end{equation}
where $G_0 (z) := (\gam h)_r +\frac{2m}{r} h_3 \gam h + \frac{m}{r} \xi_3 \xi$. 
%In the same manner as above, 
Similarly, we have 
\begin{equation}
\begin{aligned}
\nor{\del G_0 (z)}{L^4_e} &\lesssim \left( \nor{z}{L^\infty} + \nor{z'}{L^\infty} \right) 
\left( \nor{\del z_r}{L^4_e} + \nor{\frac{\del z}{r}}{L^4_e} \right) \\
&\lesssim \left( \nor{z}{L^\infty} + \nor{z'}{L^\infty} \right) 
\nor{\del z_r - \frac{m}{r} \del z}{L^4_e}
\end{aligned}
\end{equation}
where we apply Lemma \ref{L9} in the last inequality. 
By the smallness of $\nor{z}{L^\infty}$ and $\nor{z'}{L^\infty}$, 
(\ref{b248}) yields 
\begin{equation}
\nor{\del z_r - \frac{m}{r} \del z}{L^4_e} 
\lesssim \nor{\del \left( s e^{-\al R} [q\eh ] (sr) \right)}{L^4} + \nor{\del z}{\infty},
\end{equation}
and the right hand side tends to $0$ as $(q',s',\al')\to (q,s,\al)$. 
By using Lemma \ref{L9} again, we obtain $\nor{\frac{\del z}{r}}{L^4_e}\to 0$. 
Hence the proof is accomplished. \qedhere

\subsection{Approximation}
Finally, we show that each function $u(t)= e^{m\theta R}v(t,r) \in C(I;\dot{H}^1(\R^2;\bbs^2))$ can be approximated by smooth functions as follows. 

\begin{lem}\label{L11}
(i) Let $I$ be an open interval, and suppose that $u(t) = e^{m\theta R}v(t,r) \in C(I; \dot{H}^1(\R^2;\bbs^2))$ with $v(\infty) = \vec{k}$. 
Then, for every $I'\Subset I$, 
there exist $u_n(t,x) = e^{m\theta R} v_n(t,r)$, $n\in \bbn$ such that 
\begin{enumerate}[(A)]
\item $u_n(t,x) \in C^\infty (I'\times\R^2)$. 
\item For all $n\in \bbn$, there exists $R_n>0$ such that $v_n(t,r) = h(r)$ 
for all $t\in I'$ and $r\ge R_n$. 
\item $\sup_{t\in I'} \nor{u_n(t) - u(t)}{\dot{H}^1} \to 0 $ as $n\to \infty$.
\end{enumerate}
(ii) Moreover, if $\rd_t u \in C(I;L^2)$, then there exist $u_n(t,x) = e^{m\theta R} v_n(t,r)$, $n\in \bbn$ satisfying (A), (B), and the following (C)':
\begin{enumerate}
\item[(C)'] $\sup_{t\in I'} \nor{u_n(t) - u(t)}{\dot{H}^1}
+ \sup_{t\in I'} \nor{\rd_t u_n(t) - \rd_t u(t)}{L^2(\R^2)} \to 0 $ as $n\to \infty$.
\end{enumerate}
\end{lem}

\begin{rmk}
%Lemma \ref{L11} (ii) also holds when $\rd_t u \in L^{\infty} (I;L^2(\R^2))$, where the 
%supremum in (C)' is interpreted as the essential meaning. 
Such kind of approximation is originally considered in \cite{BT} in an implicit way. 
\end{rmk}

\noindent\textit{Proof of Lemma \ref{L11}.} 
The proof consists of two steps. \\
\textit{Step 1}. 
We first show the existence of sequence $\{u_n \}_{n=1}^\infty$ satisfying (A) and (C) in the case (i), 
or (A) and (C)' in the case (ii). \par
%We write $I=[a,b]$. We extend $u(t)$ to whole $\R$ by defining 
%$u(t)= u(a)$ if $t\le a$, and $u(t)=u(b)$ when $t\ge b$. 
We take an interval $I'\Subset I$ and set $\tu (t,x) := u(t,x) - Q(x)$, where $Q(x)=e^{m\theta R} h(r)$. 
We also take radially symmetric mollifiers $\eta_1 \in C_0^\infty (\R)$ and $\eta_2 \in C_0^\infty (\R^2)$, 
and then define 
$\eta (t,x) := \eta_1 (t) \eta_2(x) \in C_0^\infty(\R\times\R^2)$. \par
Here, for $(t,x)\in I'\times \R^2$ and $\vep>0$, we define 
\begin{equation}\label{a241}
\begin{aligned}
u_{\vep} (t,x) &:= Q(x) + \eta_{\vep} * \tu (t,x) \\
&= Q(x) + \int_{I\times\R^2} \eta_{\vep} (t-s, x-y) \tu (s,y)\ ds dy
\end{aligned}
\end{equation}
where $\eta_\vep (t,x) := \vep^{-3} \eta (\vep^{-1} t, \vep^{-2}x)$. (Note that (\ref{a241}) 
is well-defined if $\vep$ is sufficiently small.) 
Obviously, $u_\vep \in C^\infty (I'\times \R^2)$. 
%From the continuity of $I\ni t\to v(t) \in \dot{H}^1$, 
Since $\tu(t,x) \to 0$ as $|x|\to \infty$ uniformly on any 
subinterval of $I$, it follows that $\eta_\ep * \tu \to \tu$ as $\vep\to 0$ in $L^\infty (I'\times \R^2)$.
%
%We now claim that $\eta_\ep * \tu \to \tu$ as $\vep\to 0$ in $L^\infty (I'\times \R^2)$. 
%To this end, it suffices to show that 
%%$\tu$ is uniformly continuous on $I'\times \R^2$. (The rest )
%$\tu(t,x) \to 0$ as $|x|\to \infty$ uniformly in $t\in I''$, where $I''\Subset I$ is any 
%interval containing $I'$. 
%However, this follows immediately from (\ref{c47}) since the map $I''\ni t\to v(t) \in \dot{H}^1$ 
%is uniformly continuous. \par
%In particular, we have $u_\vep \to u$ as $\vep \to 0$ in $L^\infty (I'\times \R^2)$. 
In particular, the above claim implies that for sufficiently small $\vep$, 
%we have $|u_\vep | \ge \frac{1}{2}$. 
%Here, 
we can define 
\begin{equation}
U_n := |u_{\frac{1}{n}}|^{-1} u_{\frac{1}{n}} %\frac{u_{\frac{1}{n}}}{|u_{\frac{1}{n}}|}
\end{equation}
for sufficiently large $n\in \bbn$. 
Then, $U_n(t,\cdot)$ is $m$-equivariant for each $t$, 
since we can check that 
\begin{equation}
\tu_\vep (t,R(\theta)x)
=e^{m\theta R} \tu_{\vep} (t,x),\quad
R(\theta)=
\begin{pmatrix}
\cos \theta & -\sin \theta\\
\sin \theta & \cos \theta
\end{pmatrix}
\end{equation}
for all $\theta \in \bbt^1$. 
We can verify that $\{U_n\}_n$ satisfies the desired properties, hence Step 1 is established.
\bigskip\\
\textit{Step 2}. 
We show that $u_n$ obtained in Step 1 can be chosen with satisfying (B). 
By Step 1, we may assume $u(t) \in C(I; \dot{H}^1)\cap C^\infty (I\times\R^2)$. 
%for $I'$ as in Step 1. 
We take a cut-off function $\phi\in C^\infty (\R)$ which satisfies 
$0\le \phi (s) \le 1$ for all $s\in\R$, $\phi (s) = 1$ if $|s|\le 1$, 
and $\phi (s) =0$ when $|s|\ge 2$. 
And then we set 
$\phi_n(s) := \phi (n^{-1}s)$ for $n\in\bbn$. 
Define 
\begin{equation}
u_n(t,x) := Q(x) + \phi_n(r) \tu (t,x).
\end{equation}
Since $\tu(t,x) \to 0$ as $|x|\to \infty$ uniformly on %any 
%subinterval of $I$
$I'$, we can define
%
%Since  
%$\sup_{t\in I'} |\tu (t,x)| \to 0$ as $|x|\to \infty$ as before, 
%if $n$ is sufficiently large, 
%$|\phi_n(r) \tu (t,x)| \le \frac{1}{2} $ for all $(t,x)\in I'\times \R^2$. 
%For such $n$, we can define 
\begin{equation}
U_n := |u_{n}|^{-1} u_{n}. %\frac{u_{\frac{1}{n}}}{|u_{\frac{1}{n}}|}.
\end{equation}
for sufficiently large $n$. 
It is clear that $U_n$ is $m$-equivariant and satisfies (A) and (B). 
Thus, it suffices to show (C) or (C)'. 
We can easily check that 
$\sup_{t\in I'} \nor{u_n-u}{L^\infty} \to 0$ as $n\to \infty$. 
Moreover, denoting one of $\rd_{x_1}$, $\rd_{x_2}$, and $\rd_t$ by $\rd$, we have
\begin{equation}
\nor{\rd u_n -\rd u}{L^2_x} 
\le \frac{1}{n} \nor{\phi'}{L^\infty} \nor{\tu (t) \chi_{\{n\le |x| \le 2n \}} }{L^2_x} 
+ \nor{\rd \tu \chi_{\{|x|\ge n\}} }{ L^2_x}
\end{equation}
for each $t$. Thus
\begin{equation}
\begin{aligned}
&\sup_{t\in I'} \nor{\rd u_n -\rd u}{L^2} \\
&\hspace{20pt}\lesssim  \sup_{t\in I'} \nor{\tu (t) \chi_{\{n\le |x| \le 2n \}}}{L^\infty} 
+ \sup_{t\in I'} \nor{\tu (t) \chi_{\{|x|\ge n\}}}{L^2} \to 0
\end{aligned}
\end{equation}
as $n\to \infty$. Hence, we have $\sup_{t\in I'} \nor{\rd U_n -\rd u}{L^2} \to 0$ 
as $n\to\infty$ in the same manner as Step 1.  
Hence, we complete the proof. \hfill $\square$\bigskip\\
\textit{Acknowledgments} 
The author would like to thank Stephen Gustafson for giving him useful suggestions on 
the calculation in the proof of Lemma 4.1, and for discussing 
the problems the author kept in mind, which leads to the subjects of the present paper.
The author also would like to thank Kenji Nakanishi for useful comments, 
which are reflected on this paper.

Ikkei Shimizu\par 
{\sc Department of Mathematics\par 
Graduate School of Science\par 
Kyoto University\par 
Oiwakecho, Kitashirakawa, Sakyoku, Kyoto, 
606-8502, Japan}\par
E-mail: \url{ishimizu@math.kyoto-u.ac.jp}


\begin{thebibliography}{99}
%\bibitem{B} I. Bejenaru, \textit{Global results for Schr\"odinger maps in dimensions $n\ge 3$}, 
%Comm. Partial Differential Equations, \textbf{33} (2008), 451--477 %(X)
%\bibitem{B2} I. Bejenaru, \textit{On Schr\"odinger maps}, Amer. J. Math. 
%\textbf{130} (2008), 1033-1065. %(?)
%\bibitem{BIK} I. Bejenaru, A. D. Ionescu and C. E. Kenig, 
%\textit{Global existence and uniqueness of Schr\"odinger maps in dimensions $d\ge 4$}, 
%Adv. Math. \textbf{215} (2007), 263--291.
\bibitem{BIKT} I. Bejenaru, A. D. Ionescu, C. E. Kenig and D. Tataru, \textit{Global \scr maps in dimensions $d\ge 2$ : Small data in the critical Sobolev spaces}, 
Annals of Mathematics \textbf{173} (2011), 1443--1506.
\bibitem{BIKT2} I. Bejenaru, A. D. Ionescu, C. E. Kenig and D. Tataru, \textit{Equivariant Schr\"{o}dinger maps in two spatial dimensions}, 
Duke Math. J. \textbf{162} (2013), 1967-2025.
\bibitem{BT} I. Bejenaru and D. Tataru, Near Soliton Evolution for Equivariant Schr\"odinger Maps 
in Two Spatial Dimensions, Mem. Amer. Math. Soc. \textbf{228}, Amer. Math. Soc., Providence, 2014.
\bibitem{CH} T. Cazenave and A. Haraux, \textit{An Introduction to Semilinear Evolution Equations}, 
Oxford Lecture Series in Mathematixs and Its Applications, vol. 13, The Clarendon Press, Oxford, 1998.
  \bibitem{CSU} N. H. Chang, J. Shatah and K. Uhlenbeck, \textit{Schr\"{o}dinger maps}, 
Comm. Pure Appl. Math. \textbf{53} (2000), 590--602.
\bibitem{D} L.-S. Da Rios, \textit{On the motion of an unbounded fluid with a vortex filament of any shape}, Rend. Circ. Mat. Palermo \textbf{22} (1906), 117--135. %(Mada Check Shitenai!)
\bibitem{DW} W. Ding and Y. Wang, \textit{Local Schr\"odinger flow into K\"ahler manifolds}, 
Sci. China Ser. A-Math. \textbf{44} (2001), 1446--1464. %(O)
%\bibitem{Dod} B. Dodson, \textit{Global well-posedness for the Schr\"odinger map problem with
% small Besov norm}, 2017, arXiv:1701.08195v1. %(O) %(O)
%  \bibitem{4} A. M. Kosevich, B. A. Ivanov, A. S. Kovalev, \textit{Magnetic solitons}, Phys. Rep. 
% \textbf{194}, (1990), 117-238
\bibitem{GKT} S. Gustafson, K. Kang and T.-P. Tsai, Schr\"odinger flow near harmonic maps. Comm. Pure Appl. Math. \textbf{60} (2007), 463--499.
\bibitem{GKT2} S. Gustafson, K. Kang and T.-P. Tsai, 
\textit{Asymptotic stability of harmonic maps under the Schr\"odinger flow}, Duke Math. J. \textbf{145} (2008), 
537--583.
\bibitem{GK} S. Gustafson and E. Koo, \textit{Global well-posedness for $2D$ radial \scr maps into the sphere}, 2011, arXiv:1105.5659v1.
\bibitem{GNT} S. Gustafson, K. Nakanishi and T.-P. Tsai, \textit{Asymptotic stability, concentration, 
and oscillation in harmonic map heat-flow, Landau-Lifshitz, and Schr\"odinger maps on $\R^2$}, 
Comm. Math. Phys. \textbf{300} (2010), 205--242.
\bibitem{GT} D. Gilbarg and N. Trudinger, \textit{Elliptic Partial Differential Equations of Second Order}, Reprinted of the 1998 ed., Springer-Verlag, Berlin, 2001.
%\bibitem{IK} A. D. Ionescu and C. E. Kenig, \textit{Low-regularity Schr\"odinger maps}, Differential Integral Equations, \textbf{19} (2006), 1271--1300.
%\bibitem{IK2} A. D. Ionescu and C. E. Kenig, \textit{Low-regularity Schr\"odinger maps, II: 
%Global well-posedness in dimensions $d\ge 3$}, Comm. Math. Phys. \textbf{271} (2007), 523--559.
\bibitem{KIK} A. M. Kosevich, B. A. Ivanov and A. S. Kovalev, \textit{Magnetic Solitons}, 
Phys. Rep. \textbf{194} (1990), 117--238.
\bibitem{KLPST} C. E. Kenig, T. Lamm, D. Pollack, G. Staffilani and T. Toro, \textit{
The Cauchy problem for Schr\"odinger flows into K\"ahler manifolds}, Discrete Contin. Dyn. Syst. 
\textbf{27} (2010), 389--439. 
\bibitem{M} H. McGahagan, \textit{An approximation scheme for Schr\"odinger maps}, Comm. Partial Differential Equations \textbf{32} (2007), 375--400.
\bibitem{MRR} F. Merle, P. Rapha\"el and I. Rodnianski, \textit{
Blowup dynamics for smooth data equivariant solutions to the critical Schr\"odinger map problem}, 
Invent. Math. \textbf{193} (2013), 249--365.
\bibitem{P} G. Perelman, \textit{Blow up dynamics for equivariant critical Schr\"odinger maps}, 
Comm. Math. Phys. \textbf{330} (2014), 69--105.
\bibitem{SSB} P.-L. Sulem, C. Sulem and C. Bardos, \textit{
On the continuous limit for a system of classical spins}, Comm. Math. Phys. \textbf{107} (1986), 431--454.  
\end{thebibliography}
\end{document}